\pdfoutput=1
\RequirePackage{ifpdf}
\ifpdf 
\documentclass[pdftex]{sigma}
\else
\documentclass{sigma}
\fi

\numberwithin{equation}{section}

\newtheorem{Theorem}{Theorem}[section]
\newtheorem*{Theorem*}{Theorem}
\newtheorem{Corollary}[Theorem]{Corollary}

\newtheorem{Proposition}[Theorem]{Proposition}

\theoremstyle{definition}
\newtheorem{Definition}[Theorem]{Definition}

\newtheorem{Remark}[Theorem]{Remark}

\usepackage{tikz}
\usetikzlibrary{arrows.meta,positioning}

\begin{document}

\allowdisplaybreaks

\newcommand{\arXivNumber}{2503.11214}

\renewcommand{\PaperNumber}{040}

\FirstPageHeading

\ShortArticleName{Reformulation of $q$-Middle Convolution and Applications}

\ArticleName{Reformulation of $\boldsymbol{q}$-Middle Convolution\\ and Applications}

\Author{Yumi ARAI and Kouichi TAKEMURA}
\AuthorNameForHeading{Y.~Arai and K.~Takemura}
\Address{Department of Mathematics, Ochanomizu University, \\ 2-1-1 Otsuka, Bunkyo-ku, Tokyo 112-8610, Japan}
\Email{\mail{araiyumi.math@gmail.com}, \mail{takemura.kouichi@ocha.ac.jp}}

\ArticleDates{Received October 23, 2025, in final form April 11, 2026; Published online April 26, 2026}

\Abstract{We reformulate the $q$-convolution and the $q$-middle convolution introduced by Sakai and Yamaguchi, and we introduce $q$-analogues of the addition which is related to the gauge-transformation. A merit of the reformulation is the additivity on composition of two $q$-middle convolutions. We obtain sufficient conditions that the Jackson integrals associated with the $q$-convolution converge and satisfy the $q$-difference equation associated with the $q$-convolution. We present several third-order linear $q$-difference equations and solutions of them by using the $q$-middle convolution and the $q$-analogues of the addition.}


\Classification{33D15; 39A13}

\section{Introduction}

Differential equations and difference equations play important roles in investigating special functions, and the middle convolution has been applied for the study of special functions and differential equations.
Katz introduced the middle convolution in his book \cite{Katz} to study rigid local systems.
Dettweiler and Reiter \cite{DR1,DR2} reformulated the middle convolution for the Fuchsian equations in terms of the linear algebra.
The Fuchsian equations is the system of linear differential equations written as
\begin{equation}
\frac{{\rm d}Y}{{\rm d}x}=\left( \frac{A_0}{x-t_0}+\frac{A_1}{x-t_1} + \dots + \frac{A_N}{x-t_N} \right) Y, \label{eq:original}
\end{equation}
where $Y$ is a column vector with $m$ entries and $A_0, A_1, \dots ,A_N$ are constant square matrices of size $m$.
We review briefly the definition of the convolution and the middle convolution for equation \eqref{eq:original} (or the tuple of the matrices $(A_0, \dots ,A_N)$).
Let $\lambda \in \mathbb{C}$ and $F_i$ $(i=0,\dots ,N)$ be the square matrix of size $(N+1) m $ of the form
\begin{equation}
F_{0} =
\begin{pmatrix}
A_{0} + \lambda I_{m} & A_1 & \cdots & A_N \\
0 & 0 &\cdots & 0 \\
\vdots & \vdots & & \vdots \\
0 & 0 & \cdots & 0
\end{pmatrix}
, \qquad
F_{1} =
\begin{pmatrix}
0 & 0 & \cdots & 0 \\
A_{0} & A_1 + \lambda I_{m} & \cdots & A_N \\
\vdots & \vdots & & \vdots \\
0 & 0 & \cdots & 0
\end{pmatrix}
 ,  \ \ \dots ,
\label{eq:Fi}
\end{equation}
where $I_m $ is the identity matrix of size $m$.
Then, the correspondence of the tuple of matrices $(A_0, \dots ,A_N ) \mapsto (F_0 , \dots ,F_N )$ (or the correspondence of the associated Fuchsian system) is called the convolution.
It is shown that the following subspaces $\mathcal{K}$, $\mathcal{L}$ of $\mathbb{C}^{(N+1) m}$ are preserved by the action of $F_i$ $(i=0,\dots ,r)$.
\begin{equation*}
\mathcal{K} = \left(
\begin{matrix}
\operatorname{ker}A_{0} \\
\vdots \\
\operatorname{ker}A_{N}
\end{matrix}
\right), \qquad
\mathcal{L} = \operatorname{ker}(F_{0} + F_{1} + \cdots + F_{N}).
\end{equation*}
We denote the linear transformation induced from the action of $F_{i}$ on the quotient space $\mathbb{C}^{(N+1) m}/(\mathcal{K} + \mathcal{L})$ by $\overline{F}_{i}$.
The correspondence of the tuple of matrices $(A_0, \dots ,A_N ) \mapsto \bigl(\overline{F}_0 , \dots ,\allowbreak\overline{F}_N \bigr)$ (or the correspondence of the associated Fuchsian system) is called the middle convolution, and we denote it by $mc_{\lambda }$.
See also Haraoka \cite[Section~7.5]{Har} for a detailed description.

It is known that the middle convolution has additivity.
Namely, $mc_{0 } $ is an identity map and we have
$
mc_{\lambda} \circ mc_{\mu} = mc_{\lambda + \mu} $,
under some conditions (see the conditions $(*)$ and $(**)$ in Definition \ref{def:condqmc}).

A theory of the $q$-deformation of the middle convolution was constructed by Yamaguchi and Sakai.
A tentative theory was presented in master's thesis of Yamaguchi \cite{Yg} which was supervised by Sakai, and a full paper was published later by their joint paper \cite{SY}.
They handled the system of $q$-difference equations described as
\begin{equation}
 Y(qx) = \left( B_{\infty} + \frac{B_1}{1 - x/ b_1} +\dots + \frac{B_N}{1 - x/b_N } \right) Y(x),
\label{eq:YqxSY0}
\end{equation}
where $Y(x)$ is a column vector with $m$ entries and $B_{\infty} , B_1, \dots ,B_N$ are constant square matrices of size $(N+1)m $.
The construction of the $q$-middle convolution is similar to the case of the Fuchsian system of differential equations.
They investigated several properties of the $q$-middle convolution and related topics.
On the composition of the $q$-middle convolution, they obtained
\begin{equation}
\overline{\Psi }_{\mu } \circ \overline{\Psi }_{\lambda } \simeq \overline{\Psi }_{\log ( q^{\lambda } + q^{\mu } -1)/ \log q }.  \label{eq:Psimulam0}
\end{equation}
See the appendix for the notation.

The theory by Sakai and Yamaguchi was applied to studies of some $q$-difference equations in~\mbox{\cite{AT, STT}}.
Recall that the convolution in equation \eqref{eq:Fi} is related with Euler's integral transformation.
Let $Y(x)$ be a solution of equation \eqref{eq:original}.
Set
\begin{equation*}
W_{j}(x) = \frac{Y(x)}{x - t_{j}}, \qquad
W(x) =
\left(
\begin{matrix}
W_{0}(x) \\
\vdots \\
W_{N}(x)
\end{matrix}
\right).
\end{equation*}
We apply Euler's integral transformation for each entry of $W(x)$, i.e., we set
\begin{equation*}
\widetilde{Y} (x) = \int_{\Delta}W(s)(x - s)^{\lambda}{\rm d}s ,
\end{equation*}
where $\Delta $ is an appropriate cycle in $\mathbb{C} $ with the variable $s$.
Dettweiler and Reiter \cite{DR1} established that the function $\widetilde{Y} (x)$ satisfies the following Fuchsian system of differential equations:
\[
\frac{{\rm d}\widetilde{Y}}{{\rm d}x} = \left( \frac{F_0}{x-t_0}+\frac{F_1}{x-t_1}+\dots + \frac{F_N}{x-t_N} \right)\widetilde{Y} ,
\]
where $F_0, \dots , F_N$ were defined in equation \eqref{eq:Fi}.
Sakai and Yamaguchi \cite{SY} introduced the $q$\nobreakdash-deformation of the convolution by considering a $q$-deformation of the integral transformation of Dettweiler and Reiter.
See the appendix for a review of them.
In \cite{STT}, the theory of Sakai and Yamaguchi was applied to the $q$-Heun equation.
Consequently, a $q$-integral transformation of the $q$-Heun equation was suggested.
Note that it was discussed in \cite{Tki} from another approach based on the kernel function identity, and the $q$-integral transformation was extended to the variants of the $q$-Heun equations.
In our previous paper \cite{AT}, the $q$-convolution of Sakai and Yamaguchi was reconsidered.
Namely, one parameter which corresponds to the integral cycle was added to the $q$-integral transformation of Sakai and Yamaguchi, and convergence of the $q$-integral transformation was considered.
As an application, $q$-integral representations of solutions to the $q$-hypergeometric equation and its variants were obtained.

In this paper, we reformulate the $q$-convolution of Sakai and Yamaguchi.
Our definition of the $q$-convolution is simpler than that of Sakai and Yamaguchi, and our discussion takes care of convergence.
We rewrite equation \eqref{eq:YqxSY0} to
\begin{equation}
\frac{Y(q x) - Y(x )}{-x} = \Bigg[ \sum^{N}_{i = 0}\frac{B_{i}}{x -b_{i}} \Bigg] Y(x) ,
\label{eq:YqxBYx0}
\end{equation}
where $b_0=0$, and define the $q$-convolution and the $q$-middle convolution for the tuple $(B_0, B_1,\allowbreak \dots , B_N)$.
The $q$-difference operator in equation \eqref{eq:YqxBYx0} has already appeared in the work of Kakei and Kikuchi \cite{KK}.
A merit of our definition is the additivity of two $q$-middle convolutions.
Namely, we have
\begin{equation}
mc^q _{\lambda} \circ mc^q _{\mu} = mc^q _{\lambda + \mu} ,
\label{eq:mcaddi0}
\end{equation}
instead of equation \eqref{eq:Psimulam0} (see Section~\ref{sec:qmc} for the definition).
The $q$-convolution in this paper is also related to a $q$-integral transformation, and the $q$-middle convolution induces the $q$-integral transformation for solutions to $q$-difference equations.

We pay attention to a $q$-analogue of the addition for the Fuchsian equations, which is induced from a gauge-transformation of solutions.
Let $Y$ be a solution to the system of differential equation given in equation \eqref{eq:original}, and let $(a_0, a_1, \dots , a_N)$ be a tuple of complex numbers.
Then, the function $Z= (x-t_0)^{a_0} (x-t_1)^{a_1} \cdots (x-t_N)^{a_N} Y$ satisfies the following system of differential equations:
\begin{equation*}
\frac{{\rm d}Z}{{\rm d}x}=\left( \frac{A_0+a_0 I_m}{x-t_0}+\frac{A_1+a_1 I_m}{x-t_1}+\dots + \frac{A_N + a_N I_m}{x-t_N} \right) Z.
\end{equation*}
The correspondence of the tuple of the matrices $(A_0, A_1, \dots ,A_N ) \mapsto (A_0 + a_0 I_m , A_1 +a_1 I_m ,\dots , \allowbreak A_N +a_N I_m )$ is called the addition.

We introduce an analogue of the addition for systems of $q$-difference equations.
In particular, we investigate the transformations of the tuples $(B_0, B_1, \dots , B_N) $ induced from the transformations of the solutions by
\begin{equation} \label{eq:gaugetr}
Y(x) \mapsto x^{\mu } Y(x)\qquad \text{and} \qquad Y(x) \mapsto \frac{(x/\gamma ;q )_{\infty }}{(x/\delta ;q )_{\infty }} Y(x) ,
\end{equation}
where
$
(\alpha ;q ) _{\infty } = (1-\alpha )(1-q \alpha )\bigl(1-q^2 \alpha \bigr) \cdots $.
The transformation \eqref{eq:gaugetr}, which is related to the $q$-analogue of addition, also appeared in \cite{KK}.
See Section~\ref{sec:add} for details of the $q$-analogue of addition.

In our previous paper \cite{AT}, we obtained $q$-integral representations of solutions to the $q$-hypergeometric equation and the variants of the $q$-hypergeometric equation by applying the $q$-middle convolution.
A standard form of the $q$-hypergeometric equation is given as
\begin{equation} \label{eq:qhg-sta-sngl}
(x-q)g(x/q)+(abx-c)g(qx)-\{(a+b)x-q-c \}g(x)=0.
\end{equation}
In this paper, we obtain several $q$-difference equations by applying $q$-analogues of the addition and the reformulated $q$-middle convolution repeatedly.
Consequently, we have formal $q$-integral representations of solutions.
The reformulated $q$-middle convolution would be better for discussing a sequence of $q$-difference equations.
We discuss sufficient conditions that the $q$-integral representations converge and they satisfy the attached $q$-difference equation, which is used to generate actual solutions.
The generalized $q$-hypergeometric equation of order $3$ is written in the form
\begin{gather*}
(a_3 x-b_3)g\bigl(q ^3 x\bigr) + (a_2 x-b_2)g\bigl(q ^2 x\bigr) +(a_1 x-b_1)g(q x) + (a_0 x-b_0)g( x) =0,\\ a_3 b_3 a_0 b_0 \neq 0,
\end{gather*}
and we obtain double $q$-integral representations of solutions of this equation.
We also obtain other $q$-difference equations of order $3$, which seems to be novel (e.g., see equation \eqref{eq:21-111-111-g1}).

This paper is organized as follows.
In Section~\ref{sec:add}, we discuss $q$-analogues of the addition.
In Section~\ref{sec:qmc}, we reformulate the $q$-convolution and discuss convergence of the $q$-integral transformation which is associated with the $q$-convolution.
In Section~\ref{sec:qhypgeomgen}, we obtain the $q$-hypergeometric equation and some generalizations by applying the $q$-middle convolutions and the $q$-analogues of the addition.
We obtain $q$-integral representations of solutions formally, and discuss convergence in some cases.
In Section~\ref{sec:compqmc}, we discuss composition of the $q$-middle convolutions and establish the formula in equation \eqref{eq:mcaddi0}.
In Section~\ref{sec:rmk}, we give concluding remarks.
In the appendix, we~review the $q$-convolution and the $q$-middle convolution of Sakai and Yamaguch.

Throughout this paper, we assume that the complex number $q$ satisfies $0 <|q|<1$.

\section[q-analogues of addition]{$\boldsymbol{q}$-analogues of addition} \label{sec:add}
To reformulate the $q$-convolution in a simpler form, we replace the expression of the linear $q$-difference equation which is associated with the $q$-convolution.
Sakai and Yamaguchi used the linear $q$-difference equation
\begin{equation}
 Y(qx)= \Biggl( B_{\infty}+\sum_{i=1}^N \frac{B_i}{1-x/b_i} \Biggr) Y(x).
\label{eq:YqxSY}
\end{equation}
Set $B_0 = I_m - B_\infty - B_{1} - \dots -B_N$ and $b_0=0$.
Then, equation \eqref{eq:YqxSY}
 is equivalent to
\begin{equation}
\frac{Y(q x) - Y(x )}{-x} = \Biggl[ \frac{B_{0}}{x} + \sum^{N}_{i = 1}\frac{B_{i}}{x -b_{i}} \Biggr] Y(x) = \Biggl[ \sum^{N}_{i = 0}\frac{B_{i}}{x -b_{i}} \Biggr] Y(x).
\label{eq:YqxBYx2}
\end{equation}
Note that we can regard the left-hand side of equation \eqref{eq:YqxBYx2} as the $q$-derivative of $Y(x)$ by dividing by $1-q$.

We introduce a $q$-analogue of the addition.
If the function $Y(x) \in \mathbb{C}^m$ satisfies $Y(qx) =B(x)Y(x) $, then the function
\begin{equation}
Z(x)= x^{\mu } Y(x)
\label{eq:ZY}
\end{equation}
satisfies $Z(qx) = q^{\mu } B(x)Z(x) $.
In particular, if $Y(x)$ is a solution to equation \eqref{eq:YqxBYx2}, then the function $Z(x)$ in equation \eqref{eq:ZY} satisfies
\begin{equation*}
\frac{Z(q x) - Z(x )}{-x} = \Biggl[ \frac{(1-q^{\mu} )I_m + q^{\mu }B_0}{x} + \sum^{N}_{i = 1}\frac{q^{\mu } B_{i}}{x -b_{i}} \Biggr] Z(x).
\end{equation*}
Thus, we define the addition ${\rm add}_{\mu } $ by
\[
{\rm add}_{\mu }\colon \ (B_0, B_1, \ldots, B_N) \mapsto ((1-q^{\mu} )I_m + q^{\mu }B_0 , q^{\mu }B_1, \ldots, q^{\mu }B_N).
\]
We introduce another $q$-analogue of the addition which is induced by the gauge-transformation
\begin{equation}
Y_g (x) = \frac{( x/ b'_1;q)_{\infty }}{( x/b_1;q)_{\infty}} Y(x).
\label{eq:anotheradd}
\end{equation}
We have
\begin{equation*}
Y_g (qx) = \frac{( qx/ b'_1;q)_{\infty }}{( qx/b_1;q)_{\infty}} Y(qx) = \frac{(1-x/b_1)}{(1-x/b'_1)} \frac{( x/ b'_1;q)_{\infty }}{( x/b_1;q)_{\infty}} Y(qx).
\end{equation*}
If the function $Y(x) \in \mathbb{C}^m$ satisfies $Y(qx) =B(x)Y(x) $, then the function $Y_g (x)$ satisfies
\begin{equation*}
Y_g (qx)
= \frac{(1-x/b_1)}{(1-x/b'_1)} \frac{( x/ b'_1;q)_{\infty }}{( x/b_1;q)_{\infty}} B(x) Y(x) = \frac{(1-x/b_1)}{(1-x/b'_1)} B(x) Y_g (x).
\end{equation*}
Note that equation \eqref{eq:YqxBYx2} is equivalent to
\begin{equation}
Y(q x) = B(x) Y(x) , \qquad B(x)= \Biggl[ I -B_{0} + \sum^{N}_{i = 1} B_{i} \frac{-x }{x -b_{i}} \Biggr].
\label{eq:YqxBxYxBx}
\end{equation}
By rewriting $B(x) (1-x/b_1)/(1-x/b'_1) $ as the form of the second equation in \eqref{eq:YqxBxYxBx}, we obtain
\begin{align*}
Y_g (qx) ={}& \left[ I -B_{0} + \left\{ \left(1- \frac{b'_1 }{b_1} \right) ( I - B_{0} ) + \frac{b'_1}{b_1} B_{1} - \sum^{N}_{i = 2} \frac{1 - b'_1 / b_1 }{ 1 - b_i/ b'_1 } B_{i} \right\} \frac{-x }{x -b'_1} \right.\\
&\left.+ \sum^{N}_{i = 2} \frac{1 - b_i /b_1}{1 - b_i /b'_1 } B_{i} \frac{-x }{x -b_{i}} \right] Y_g (x).
\end{align*}
Therefore, the gauge-transformation \eqref{eq:anotheradd} induces the change of the poles $(b_1, b_2, \dots ,b_N) \mapsto (b'_1, b_2, \dots ,b_N)$ and the tuple of the matrices $(B_0, B_1, \dots ,B_N) \mapsto (B'_0, B'_1, \dots ,B'_N)$, where
\begin{align*}
& B'_0=B_0, B'_1 = \left(1- \frac{b'_1 }{b_1} \right) ( I - B_{0} ) + \frac{b'_1}{b_1} B_{1} - \sum^{N}_{i = 2} \frac{1 - b'_1 / b_1 }{ 1 - b_i/ b'_1 } B_{i} , \\
& B'_j = \frac{1 - b_j /b_1}{1 - b_j /b'_1 } B_{j},\qquad j=2, \dots ,N.
\end{align*}
We may also regard this correspondence as a $q$-analogue of addition.

\section[q-middle convolution and convergence]{$\boldsymbol{q}$-middle convolution and convergence} \label{sec:qmc}

We reformulate the $q$-convolution by modifying the definition of Sakai and Yamaguchi.
See the appendix for the $q$-convolution of Sakai and Yamaguchi.
\begin{Definition}[reformulation of $q$-convolution] \label{def:tqc}
Let $ ( B_{0}, B_{1} ,\dots ,B_N ) $ be the tuple of $m\times m $ matrices and $\lambda \in \mathbb{C} $.
We define the tuple $( G_0, G_1, \dots , G_N )$ of $(N+1)m \times (N+1)m$ matrices as follows:
\begin{gather}
 G_0 = \left(
 \begin{matrix}
 q^{-\lambda } B_0 + \bigl(1-q^{-\lambda}\bigr)I_m & q^{-\lambda }B_1 & \cdots & \cdots & q^{-\lambda } B_N \\
 0 & 0 & \cdots & \cdots & 0 \\
 \vdots & \vdots & {} & {} & \vdots \\
 \end{matrix}
\right), \nonumber\\
 \qquad \qquad \qquad \vdots \nonumber\\
 G_i = \left(
 \begin{matrix}
 {} & {} & O & {} & {} \\
 q^{-\lambda } B_0 & \cdots & q^{-\lambda }B_i + \bigl(1-q^{-\lambda}\bigr)I_m & \cdots & q^{-\lambda } B_N \\
 {} & {} & O & {} & {}
 \end{matrix}
\right)
 {\scriptstyle (i+1)}, \nonumber\\
 \qquad \qquad \qquad \vdots\nonumber \\
 G_N = \left(
 \begin{matrix}
 \vdots & {} & {} & \vdots & \vdots \\
 0 & \cdots & \cdots & 0 & 0 \\
 q^{-\lambda } B_0 & \cdots & \cdots & q^{-\lambda }B_{N-1} & q^{-\lambda } B_N + \bigl(1-q^{-\lambda}\bigr)I_m
 \end{matrix}
\right).  \label{eq:bG}
\end{gather}
The correspondence $c ^q _\lambda\colon( B_{0}, B_{1} ,\dots ,B_N ) \mapsto ( G_0, G_1, \dots , G_N )$ is called the $q$-convolution.
\end{Definition}
The $q$-convolution in Definition \ref{def:tqc} is related with a $q$-integral transformation on the linear $q$-difference equations in the form of equation \eqref{eq:YqxBYx2}.
Let $\xi \in \mathbb{C} \setminus \{ 0 \}$.
The Jackson integral on the open interval $(0,\infty )$ is defined as the infinite sum
\begin{equation*}
 \int^{\xi \infty }_{0}f(s) {\rm d}_{q}s = (1-q)\sum^{\infty}_{n=-\infty}q^{n} \xi f(q^n \xi ).
\end{equation*}
Let $K_{\lambda}(x, s) $ be a function which satisfies
\begin{equation}
q^{\lambda } K_{\lambda}(qx, s) = K_{\lambda}(x, s/q) = \frac{ x - q^{\lambda } s}{x - s} K_{\lambda}(x, s).  \label{eq:tPlambda}
\end{equation}
Then, the functions
\begin{equation}
K^{(1)}_{\lambda}(x, s) = x^{-\lambda } \frac{\bigl(q^{\lambda+1}s/x;q\bigr)_{\infty}}{(qs/x;q)_{\infty}} , \qquad K^{(2)}_{\lambda}(x, s) = s^{-\lambda } \frac{(x/s;q)_{\infty}}{\bigl(q^{-\lambda }x/s;q\bigr)_{\infty}}
 \label{eq:tP1tP2}
\end{equation}
satisfy equation \eqref{eq:tPlambda}.
The functions
\begin{equation*}
 \frac{\bigl(q^{\lambda+1}s/x , b x , q/(b x) ;q\bigr)_{\infty}}{\bigl(qs/x, q^{-\lambda } b x, q^{\lambda +1} /(bx) ;q\bigr)_{\infty}} , \qquad \frac{(x/s , b s , q/(b s) ;q)_{\infty}}{\bigl(q^{-\lambda }x/s , q^{-\lambda } b s, q^{\lambda +1} /(bs) ;q\bigr)_{\infty}}
\end{equation*}
also satisfy equation \eqref{eq:tPlambda} for any $b \in \mathbb{C} \setminus \{ 0 \} $.

The following theorem indicates a correspondence of the solutions to $q$-difference equations associated with the $q$-convolution.
Although this theorem is related to the results in \cite[Theorem~2.1]{SY} and \cite[Corollary 2.3]{AT}, it is not a mere reformulation of them.
In the present setting, the definition of the $q$-convolution is modified, and the variable $\xi$ is allowed to be either independent of $x$ or proportional to $x$, which is not covered in \cite{AT,SY}.
\begin{Theorem}[{cf.\ \cite[Theorem 2.1]{SY}, \cite[Corollary 2.3]{AT}}] \label{thm:qcint}

Set $b_0=0$.
Let $Y(x)$ be a solution~to
\begin{equation}
\frac{Y(q x) - Y(x )}{-x} = \Biggl[ \sum^{N}_{i = 0}\frac{B_{i}}{x -b_{i}} \Biggr] Y(x) ,
\label{eq:YqxBYx}
\end{equation}
and set
\begin{equation}
 \widetilde{Y}_{i}(x) = \int^{\xi \infty}_{0}\frac{K_{\lambda}(x, s)}{s-b_{i}}Y(s) {\rm d}_{q}s,\qquad i=0,\dots ,N, \qquad \widetilde{Y}(x) = \begin{pmatrix} \widetilde{Y}_{0}(x) \\ \vdots \\ \widetilde{Y}_{N}(x) \end{pmatrix}.  \label{eq:qcintadd}
\end{equation}
Assume that the variable $\xi $ is independent of the variable $x$ or it is proportional to $x$ {\rm(}i.e.,~$\xi =Ax$ where $A $ is independent of $x${\rm)}.
If every element of $\widetilde{Y}_{i} (x) $ converges for $i=0,1 , \dots ,N$ and
\begin{equation}
\lim _{K \to -\infty} K_{\lambda}\bigl(x, q^{K-1} \xi\bigr)Y\bigl(q^{K} \xi\bigr) =0,\qquad \lim _{L \to +\infty } K_{\lambda}\bigl(x, q^{L} \xi\bigr)Y\bigl(q^{L+1} \xi\bigr) =0,
\label{eq:tPKLlim0}
\end{equation}
then the function $\widetilde{Y}(x)$ satisfies
\begin{equation*}
\frac{\widetilde{Y}(q x) - \widetilde{Y}(x )}{-x} = \Biggl[ \sum^{N}_{i = 0}\frac{G_{i}}{x -b_{i}} \Biggr] \widetilde{Y}(x).
\end{equation*}
\end{Theorem}
To prove the theorem, we define the function $\widetilde{Y}_{i}^{[K,L]}(x)$ by
\begin{equation}
\widetilde{Y}_{i}^{[K,L]}(x) = (1-q) \sum _{n=K}^{L} \frac{K_{\lambda}(x, s)}{s-b_{i}} s Y(s) |_{s= q^n \xi }.
\label{eq:tYiKLdef}
\end{equation}
\begin{Proposition} \label{prop:convKL3}
Let $Y(x)$ be a solution to equation \eqref{eq:YqxBYx} and assume that the variable $\xi $ is independent of the variable $x$.
For $i=0,1,\dots ,N $, we have
\begin{gather}
 \frac{\widetilde{Y}_{i}^{[K,L]}(qx) - \widetilde{Y}_{i}^{[K,L]}(x)}{-x} \nonumber \\
\qquad = \frac{1 - q^{-\lambda } }{x - b_{i}} \widetilde{Y}_{i}^{[K,L]}(x) + \frac{1}{x - b_{i}}\sum^{N}_{j = 0} q^{-\lambda } B_{j} \widetilde{Y}_{j}^{[K,L]}(x) \nonumber \\
 \phantom{\qquad = }{} - \frac{(1-q) q^{-\lambda } }{x - b_{i}} \bigl\{ K_{\lambda}\bigl(x, q^{K-1} \xi\bigr)Y\bigl(q^{K} \xi\bigr) - K_{\lambda}\bigl(x, q^{L} \xi\bigr)Y\bigl(q^{L+1} \xi\bigr) \bigr\}.  \label{eq:wtYiqx3}
\end{gather}
\end{Proposition}
\begin{proof}
It follows from equation \eqref{eq:tPlambda} that
\begin{equation*}
 \frac{K_{\lambda}(qx, s) - K_{\lambda}(x, s)}{-x} = \frac{1- q^{-\lambda } }{x-s} K_{\lambda}(x, s) , \qquad \frac{K_{\lambda}(x, s/q) - K_{\lambda}(x, s)}{s} = \frac{1-q^{\lambda }}{x-s} K_{\lambda}(x, s).
\end{equation*}
Since
$
 \frac{1}{(x-s)(s-b_i)} = \frac{1}{(s-b_i)(x-b_i)} + \frac{1}{(x-s)(x-b_i)}$,
we have
\begin{align*}
 \frac{K_{\lambda}(qx, s) - K_{\lambda}(x, s)}{-x(s-b_i)}
&= \frac{ 1- q^{-\lambda } }{(s-b_i)(x-b_i)} K_{\lambda}(x, s) + \frac{ 1- q^{-\lambda } }{(x-s)(x-b_i)} K_{\lambda}(x, s) \\
& = \frac{ 1- q^{-\lambda } }{(s-b_i)(x-b_i)} K_{\lambda}(x, s) - \frac{q^{-\lambda } }{x-b_i} \frac{K_{\lambda}(x, s/q) - K_{\lambda}(x, s)}{s}.
\end{align*}
By multiplying by $(1-q) s Y(s)$ and taking summation, we have
\begin{align}
& \frac{\widetilde{Y}_{i}^{[K,L]}(qx) - \widetilde{Y}_{i}^{[K,L]}(x)}{-x} \nonumber \\
&\qquad =  \frac{ 1- q^{-\lambda } }{x - b_{i}} \widetilde{Y}_{i}^{[K,L]}(x)
- \frac{ q^{-\lambda } (1-q) }{x - b_{i}}\sum _{n=K}^{L} (K_{\lambda}(x, s/q) - K_{\lambda}(x, s)) Y(s)|_{s= q^n \xi }.  \label{eq:YiKLqxx}
\end{align}
It follows from equation \eqref{eq:YqxBYx} that
\begin{gather*}
 \sum _{n=K}^{L} (K_{\lambda}(x, s/q) - K_{\lambda}(x, s)) Y(s)|_{s= q^n \xi } \\
\qquad= K_{\lambda}\bigl(x, q^{K-1} \xi\bigr)Y\bigl(q^{K} \xi\bigr) - K_{\lambda}\bigl(x, q^{L} \xi\bigr)Y\bigl(q^{L+1} \xi\bigr) + \sum _{n=K}^{L} K_{\lambda}(x, s)( Y(q s) -Y(s))|_{s= q^n \xi } \\
\qquad = K_{\lambda}\bigl(x, q^{K-1} \xi\bigr)Y\bigl(q^{K} \xi\bigr) -K_{\lambda}\bigl(x, q^{L} \xi\bigr)Y\bigl(q^{L+1} \xi\bigr) - \sum _{n=K}^{L} \sum^{N}_{j = 0}\frac{ s K_{\lambda}(x, s)}{s - b_{j}}B_{j} Y(s) |_{s= q^n \xi } \\
\qquad = K_{\lambda}\bigl(x, q^{K-1} \xi\bigr)Y\bigl(q^{K} \xi\bigr) - K_{\lambda}\bigl(x, q^{L} \xi\bigr)Y\bigl(q^{L+1} \xi\bigr) - \frac{1 }{1-q}\sum^{N}_{j = 0} B_{j} \widetilde{Y}_{j}^{[K,L]}(x).
\end{gather*}
Therefore, we obtain equation \eqref{eq:wtYiqx3}.
\end{proof}

The following corollary follows from the definition of the matrices $G_{i}$ $(i=0,1,\dots ,N)$.
\begin{Corollary} \label{cor:convKLti}
Let $Y(x)$ be a solution to equation \eqref{eq:YqxBYx} and assume that the variable $\xi $ is independent of the variable $x$.
Set
\begin{align*}
& g^{[K,L]}(x) = (q-1) q^{-\lambda } \bigl\{ K_{\lambda}\bigl(x, q^{K-1} \xi\bigr)Y\bigl(q^{K} \xi\bigr) - K_{\lambda}\bigl(x, q^{L} \xi\bigr)Y\bigl(q^{L+1} \xi\bigr) \bigr\} , \\
& \widetilde{Y}^{[K,L]} (x) = \begin{pmatrix} \widetilde{Y}^{[K,L]}_{0}(x) \\ \vdots \\ \widetilde{Y}^{[K,L]}_{N}(x) \end{pmatrix} ,\qquad
 \widetilde{g}^{[K,L]} (x) = \begin{pmatrix} (x - b_{0})^{-1} g^{[K,L]}(x) \\ \vdots \\ (x - b_{N})^{-1} g^{[K,L]}(x) \end{pmatrix}.
\end{align*}
Then, we have
\begin{equation*}
\frac{\widetilde{Y} ^{[K,L]} (q x) - \widetilde{Y} ^{[K,L]} (x )}{-x} = \Biggl[ \sum^{N}_{i = 0}\frac{G_{i}}{x -b_{i}} \Biggr] \widetilde{Y} ^{[K,L]} (x) + \widetilde{g}^{[K,L]} (x).
\end{equation*}
\end{Corollary}
In the case $\xi =Ax$, Proposition~\ref{prop:convKL3} is replaced as follows.

\begin{Proposition} \label{prop:convKL3xiAx}
Let $Y(x)$ be a solution to equation \eqref{eq:YqxBYx} and assume that $\xi =Ax$.
For $i=0,1,\dots ,N $, we have
\begin{gather*}
 \frac{\widetilde{Y}_{i}^{[K,L]}(qx) - \widetilde{Y}_{i}^{[K,L]}(x)}{-x}\\
 \qquad= \frac{1 - q^{-\lambda } }{x - b_{i}} \widetilde{Y}_{i}^{[K,L]}(x) + \frac{1}{x - b_{i}}\sum^{N}_{j = 0} q^{-\lambda } B_{j} \widetilde{Y}_{j}^{[K,L]}(x) \\
 \phantom{\qquad=}{} - \frac{(1-q) q^{-\lambda } }{x - b_{i}} \bigl\{ K_{\lambda}\bigl(x, q^{K-1} \xi\bigr)Y\bigl(q^{K} \xi\bigr) - K_{\lambda}\bigl(x, q^{L} \xi\bigr)Y\bigl(q^{L+1} \xi\bigr) \bigr\} \\
 \phantom{\qquad=}{} + \frac{1-q}{x} \left[ \frac{K_{\lambda}(qx, s) }{s-b_i} s Y(s) |_{s=A q^K x } - \frac{K_{\lambda}(qx, s) }{s-b_i} s Y(s) |_{s=A q^{L+1} x } \right].
\end{gather*}
\end{Proposition}
\begin{proof}
In the case $\xi =Ax$, we have
\begin{align*}
& \widetilde{Y}_{i}^{[K,L]}(x) = (1-q) \sum _{n=K}^{L} \frac{K_{\lambda}(x, A q^n x )}{A q^n x -b_{i}} A q^n x Y(A q^n x) , \\
& \widetilde{Y}_{i}^{[K,L]}(qx) = (1-q) \sum _{n=K}^{L} \frac{K_{\lambda}(qx, s)}{s-b_{i}} s Y(s) |_{s= A q^{n+1} x }.
\end{align*}
It follows that
\begin{gather*}
 (1-q) \sum _{n=K}^{L} \frac{K_{\lambda}(qx, s) }{s-b_i} s Y(s) |_{s=A q^n x }\\
 \qquad = (1-q) \sum _{n=K-1}^{L-1} \frac{K_{\lambda}(qx, s) }{s-b_i} s Y(s) |_{s=A q^{n+1} x } \\
\qquad = \widetilde{Y}_{i}^{[K,L]}(qx) + (1-q) \frac{K_{\lambda}(qx, s) }{s-b_i} s Y(s) |_{s=A q^K x } - (1-q) \frac{K_{\lambda}(qx, s) }{s-b_i} s Y(s) |_{s=A q^{L+1} x }.
\end{gather*}
Therefore, in the case $\xi =Ax$, equation \eqref{eq:YiKLqxx} is replaced to
\begin{gather*}
 \frac{\widetilde{Y}_{i}^{[K,L]}(qx) - \widetilde{Y}_{i}^{[K,L]}(x)}{-x} + \frac{1-q}{-x} \left[ \frac{K_{\lambda}(qx, s) }{s-b_i} s Y(s) |_{s=A q^K x } - \frac{K_{\lambda}(qx, s) }{s-b_i} s Y(s) |_{s=A q^{L+1} x } \right] \\
\qquad = \frac{ 1- q^{-\lambda } }{x - b_{i}} \widetilde{Y}_{i}^{[K,L]}(x) - \frac{ q^{-\lambda } (1-q) }{x - b_{i}}\sum _{n=K}^{L} (K_{\lambda}(x, s/q) - K_{\lambda}(x, s)) Y(s)|_{s= q^n \xi } ,
\end{gather*}
and we obtain Proposition~\ref{prop:convKL3xiAx} by repeating the proof of Proposition~\ref{prop:convKL3}.
\end{proof}

We now prove Theorem~\ref{thm:qcint}.
The convergence of $\widetilde{Y}_{i} (x) $ is described as the convergence of $\widetilde{Y}_{i}^{[K,L]}(x) $ in equation \eqref{eq:tYiKLdef} as $ K \to -\infty $ and $L \to +\infty $, and we have
\[
 \widetilde{Y}_{i}(x) = \lim _{\substack{K \to -\infty \\ L \to +\infty}} \widetilde{Y}_{i}^{[K,L]}(x).
 \]
The theorem for the case that the variable $\xi $ is independent of the variable $x$ follows from Corollary~\ref{cor:convKLti}, because the function $ \widetilde{g}^{[K,L]} (x)$ in Corollary \ref{cor:convKLti} converges to $0$ as $ K \to -\infty $ and~${L \to +\infty }$ by equation \eqref{eq:tPKLlim0}.
To obtain the theorem for the case $\xi =Ax$, we use the fact that the convergence of the summation $a_1 + a_2 + \cdots $ implies the condition $a_n \to 0 $ as $n \to \infty $.
Since the Jackson integral $\widetilde{Y}_{i} (x) $ converges, we have
\begin{equation*}
 \frac{K_{\lambda}(qx, s) }{s-b_i} s Y(s) |_{s=A q^n x } \to 0
\end{equation*}
as $n \to +\infty $ and $n \to -\infty $, the supplementary terms in Proposition~\ref{prop:convKL3xiAx} converge as $ K \to -\infty $ and $L \to +\infty $, and we obtain the theorem.

We give a sufficient condition that the assumption of the convergence holds, which was essentially established in \cite[Proposition 2.5]{AT}.
\begin{Proposition} \label{prop:qcintconv}
Let $[Y(s )]_k $ be the $k$-th component of $Y(s ) (\in \mathbb{C}^{m} )$.
Assume that the integrand of $\widetilde{Y}_{i} (x) $ in equation \eqref{eq:qcintadd} is well defined $($i.e., does not diverge$)$ about $s= q^n \xi $ for any $i \in \{ 0 , 1, \dots , N \}$ and $n \in \mathbb{Z} $.
\begin{itemize}\itemsep=0pt
\item[$(i)$] If there exist $\varepsilon _1 , C_1 \in \mathbb{R} _{>0} $ and $M_1 \in \mathbb{Z} $ such that $|[Y(s )]_k | \leq C_1 |s |^{\varepsilon _1} $ for any $s \in \{ q^{n}\xi \mid n \geq M_1 ,\, n \in \mathbb{Z} \} $ and $k=1,\dots ,m$, then every component of \smash{$\widetilde{Y}_{j}^{[K,L]}(x) $} converges absolutely as~${L \to +\infty }$ for each $K \in \mathbb{Z} $ and $j=0,1,\dots ,N$, and
\begin{equation}
\lim _{L \to +\infty } K_{\lambda}\bigl(x, q^{L} \xi\bigr)Y\bigl(q^{L+1} \xi\bigr) =0.
\label{eq:limL}
\end{equation}
\item[$(ii)$] If there exist $\varepsilon _2 , C_2 \in \mathbb{R} _{>0} $ and $M_2 \in \mathbb{Z} $ such that $|[Y(s )]_k | \leq C_2 |s ^{\lambda }| |s | ^{ - \varepsilon _2}$ for any $s \in \{ q^{n}\xi \mid n \leq - M_2 ,\, n \in \mathbb{Z} \} $ and $k=1,\dots ,m$, then every component of \smash{$\widetilde{Y}_{j}^{[K,L]}(x) $} converges absolutely as $K \to -\infty $ for each $L \in \mathbb{Z} $ and $j=0,1,\dots ,N$, and
\begin{equation}
\lim _{K \to -\infty} K_{\lambda}\bigl(x, q^{K-1} \xi\bigr)Y\bigl(q^{K} \xi\bigr) =0.
\label{eq:limK}
\end{equation}
\end{itemize}
\end{Proposition}
\begin{proof}
We show (i).
It follows from equation \eqref{eq:tPlambda} that
\begin{equation*}
\frac{K_{\lambda}\bigl(x, q^{n+1} \xi \bigr)}{K_{\lambda}(x, q^{n} \xi )} = \frac{1- q^{ n+1}\xi /x }{1- q^{\lambda +n+1}\xi /x } \to 1
\end{equation*}
as $n \to +\infty $ for each $\xi $ and $x$.
Since $|q| ^{-\varepsilon _1 /2} >1$, there exists an integer $N'_{1} $ such that $\big| K_{\lambda}\bigl(x, q^{n+1} \xi \bigr) / K_{\lambda}(x, q^{n} \xi )\big| \leq |q| ^{-\varepsilon _1 /2} $ for any integer $n$ no less than $N'_{1} $.
Hence, we have
\begin{equation}
| K_{\lambda}(x, q^{n } \xi ) | \leq B' _1 |q|^{-n \varepsilon _1 /2 }
\label{eq:tPB1}
\end{equation}
for any integer $n$ no less than $N'_{1} $ by setting $B' _1 = |q|^{N'_{1} \varepsilon _1 /2 } \big| K_{\lambda}\bigl(x, q^{N'_{1} } \xi \bigr) \big| $.
By combining with the assumption, there exist an integer $M' _1$ and a positive number $C' _1 $ such that
\begin{equation*}
\left|\left. \frac{ K_{\lambda}(x, s)}{s-b_{j}} s [Y(s)]_k \right| _{s= q^{n}\xi } \right|\leq C'_1 |q |^{ n \varepsilon _1 /2 }
\end{equation*}
for any integer $ n$ no less than $ M'_1 $, $k \in \{ 1, \dots , m \}$ and $j \in \{ 0, 1, \dots , N \}$.
Absolute convergence of the summation
\begin{equation*}
\sum _{n=M'_1}^{\infty } \frac{ K_{\lambda}(x, s)}{s-b_{i}} s[Y(s)]_k |_{s= q^{n}\xi }
\end{equation*}
follows from the convergence of the majorant series $C'_1 |q |^{ n \varepsilon _1 /2 }$, and every component of \smash{$\widetilde{Y}_{j}^{[K,L]}(x) $} converges absolutely as $L \to +\infty $ for each $K \in \mathbb{Z} $ and $j=0,1,\dots ,N$.
Equation \eqref{eq:limL} is also shown by using equation \eqref{eq:tPB1} and the assumption.

We show (ii).
It follows from equation \eqref{eq:tPlambda} that
\begin{equation*}
\frac{K_{\lambda}\bigl(x, q^{-n-1} \xi \bigr)}{K_{\lambda}(x, q^{-n} \xi )} = \frac{q^{\lambda } - q^n x/\xi }{1- q^n x/\xi } \to q^{\lambda }
\end{equation*}
as $n \to +\infty $ for each $\xi $ and $x$.
Since $|q| ^{-\varepsilon _2 /2} >1$, there exists an integer $N'_{2}$ and a positive number $B' _2 $ such that $\big| K_{\lambda}\bigl(x, q^{-n-1} \xi \bigr) / K_{\lambda}(x, q^{-n} \xi ) \big| \leq |q| ^{-\varepsilon _2 /2} \big|q^{\lambda } \big|$ and
\begin{equation}
| K_{\lambda}(x, q^{-n } \xi ) | \leq B' _2 |q|^{-n \varepsilon _2 /2 } \big|q^{\lambda } \big|^{n }
\label{eq:tPB2}
\end{equation}
for any integer $n$ no less than $N'_{2} $.
By the assumption, we have
\begin{equation*}
| [Y(q^{-n}\xi )]_k | \leq C_2 \big|q ^{\lambda }\big|^{-n} \big| \xi ^{\lambda - \varepsilon _2} \big| |q |^{ n \varepsilon _2 }
\end{equation*}
for any integer $n$ no less than $M_2 $.
Hence, there exist an integer $M' _2$ and a positive number $C' _2 $ such that
\begin{equation*}
\left|\left. \frac{ K_{\lambda}(x, s)}{s-b_{j}} s [Y(s)]_k \right| _{s= q^{-n}\xi } \right|\leq C'_2 |q |^{ n \varepsilon _2 /2 }
\end{equation*}
for any integer $ n$ no less than $ M'_2 $, $k \in \{ 1, \dots , m \}$ and $j \in \{ 0, 1, \dots , N \}$.
Thus, every component of \smash{$\widetilde{Y}_{j}^{[K,L]}(x) $} converges absolutely as $K \to -\infty $ for each $L \in \mathbb{Z} $ and $j=0,1,\dots ,N$.
Equation~\eqref{eq:limK} is also shown by using equation \eqref{eq:tPB2} and the assumption.
\end{proof}

\begin{Proposition} \label{prop:qcintconvYx}
Set $b_0=0$.
Let $Y(x)$ be a solution to
\begin{equation}
\frac{Y(q x) - Y(x )}{-x} = \Biggl[ \sum^{N}_{i = 0}\frac{B_{i}}{x -b_{i}} \Biggr] Y(x).
\label{eq:YqxBYx00}
\end{equation}
\begin{itemize}\itemsep=0pt
\item[$(i)$] If the absolute value of any eigenvalue of $I_m - B_0$ is strictly less than $1$, then there exist~${\varepsilon _1 , C_1 \in \mathbb{R}_{>0}} $ and $M_1 \in \mathbb{Z} $ such that $|[Y(s )]_k | \leq C_1 |s |^{\varepsilon _1} $ for any $s \in \{ q^{n}\xi \mid n \geq M_1 ,\, n \in \mathbb{Z} \} $ and $k=1,\dots ,m$.
\item[$(ii)$] If the absolute value of any eigenvalue of $I_m - B_0 -B_1 - \dots - B_N$ is strictly more than~$\big|q^{\lambda }\big|$, then there exist $\varepsilon _2 , C_2 \in \mathbb{R}_{>0} $ and $M_2 \in \mathbb{Z} $ such that $|[Y(s )]_k | \leq C_2 \big|s ^{\lambda }\big| |s | ^{ - \varepsilon _2}$ for any~${s \in \{ q^{n}\xi \mid n \leq -M_2 ,\, n \in \mathbb{Z} \} }$ and $k=1,\dots ,m$.
\end{itemize}
\end{Proposition}
\begin{proof}
Write equation \eqref{eq:YqxBYx00} as
$
Y(q x) = B(x) Y(x) $.
Then, we have $B(x) \to I_m -B_0$ as $x \to 0$ and $B(x) \to I_m -B_0 -B_1 - \dots - B_N$ as $x \to \infty$.

We show (i).
Let $e_1, e_2, \dots ,e_m$ be the eigenvalues of the matrix $I_m- B_0 $.
By the assumption, we have $\max \{ |e_1|, |e_2|, \dots , |e_m| \} <1 $.
There exists $\varepsilon _1 >0$ such that $|e_k | <|q|^{\varepsilon _1} $ for any $k \in \{1, \dots ,m \}$.
Set $ \varepsilon = |q|^{\varepsilon _1} - \max \{ |e_1|, |e_2|, \dots , |e_m| \} (>0)$.
There exist numbers $p_j \in \{ 0, \varepsilon /2 \}$ $(j=1, \dots , m-1)$ and an invertible matrix $Q_0$ such that
\begin{equation*}
 Q_0^{-1} ( I_m -B_0 ) Q_0 =
\begin{pmatrix}
e_1 & p_1 & \cdots & 0 & 0 \\
0 & e_2 & \ddots & 0 & 0 \\
\vdots & \ddots & \ddots & \ddots & \vdots \\
0 & 0 & \ddots & e_{m-1} & p_{m-1} \\
0 & 0 & \cdots & 0 & e_{m}
\end{pmatrix}.
\end{equation*}
Set $p_m=0$ for later use.
The numbers $p_1, \dots ,p_{m-1}$ are determined by the Jordan normal form of $ I_m -B_0 $.
If $I_m- B_0 $ is diagonalizable, then $p_1 = \dots = p_{m-1} =0$.
Write
\begin{equation} \label{eq:Q0-BxQ0}
 Q_0^{-1} B(x) Q_0 = Q_0^{-1} ( I_m -B_0 ) Q_0 + \begin{pmatrix}
b_{1,1} (x) & \cdots & b_{1,m} (x) \\
\vdots & & \vdots \\
 b_{m,1} (x) & \cdots & b_{m,m} (x)
\end{pmatrix}.
\end{equation}
Since $Q_0^{-1} B(x) Q_0 \to Q_0^{-1} ( I_m -B_0 ) Q_0 $ as $x \to 0$, there exists $\delta >0$ such that $|b_{i,j}(x)| < \varepsilon /(2m)$ for any $i,j \in \{1, \dots ,m \}$ and any $x$ $(|x| < \delta )$.
Let $Y (x)$ be a solution to $Y (qx) = B(x)Y (x)$.
Then, the function $Z(x) = Q_0^{-1} Y (x)$ satisfies $Z(qx) = Q_0^{-1} B(x) Q_0 Z(x)$.
Let $[Z(s )]_k $ be the $k$-th component of $Z(s ) (\in \mathbb{C} ^{m} )$.
We have
\begin{equation*}
[Z(qx )]_k = e_k [Z(x )]_k + p_k [Z(x )]_{k+1} + b_{k,1}(x) [Z(x )]_1 + \dots + b_{k,m}(x) [Z(x )]_m.
\end{equation*}
We assume that $|x| < \delta $ and $|[Z(x )]_k| < \widetilde{A}$ for any $k \in \{1, \dots ,m \}$.
Then,
\begin{gather}
 | [Z(qx )]_k | \leq |e_k| | [Z(x )]_k| + |p_k | | [Z(x )]_{k+1}| + |b_{k,1}(x)| | [Z(x )]_1| + \dots + |b_{k,m}(x) | |[Z(x )]_m | \nonumber \\
\qquad \leq \widetilde{A} |e_k| + \widetilde{A} \varepsilon /2 + \widetilde{A} m \varepsilon /(2m) \leq \widetilde{A} ( \max \{ |e_1|, |e_2|, \dots , |e_m| \} + \varepsilon ) = |q|^{\varepsilon _1} \widetilde{A} \label{eq:Zqxineq}
\end{gather}
for any $k \in \{1, \dots ,m \}$.
We also have
\begin{equation}
 | [Z(q^n x )]_k | \leq |q|^{n \varepsilon _1} \widetilde{A} \label{eq:ineq}
\end{equation}
for $n \in \mathbb{Z} _{\geq 0} $ and $k \in \{1, \dots ,m \}$.
Let $\xi \in \mathbb{C} \setminus \{0\} $.
There exists an integer $M_0$ such that $| q^n \xi | < \delta $ for any integer $n \geq M_0$.
We take a positive number $C' $ such that $\big| \big[Z\bigl( q^{M_0} \xi \bigr)\big]_k\big| \leq C' $ for any~${k \in \{1, \dots ,m \}}$.
It follows from the inequality \eqref{eq:ineq} that
\[
\big| \big[Z\bigl( q^{M_0 +n} \xi \bigr)\big]_k\big| \leq |q|^{n \varepsilon _1 } C' = \big| q^{M_0 +n} \xi \big|^{\varepsilon _1} \big| q^{M_0 } \xi \big|^{-\varepsilon _1} C'
\]
for any $k \in \{1, \dots ,m \}$ and $n \in \mathbb{Z} _{\geq 0}$.
Since $Y (x) = Q_0 Z(x)$, there exists a positive number $C$ such that $\big| \big[Y\bigl( q^{M_0 +n} \xi \bigr)\big]_k\big| \leq \big|q^{M_0 +n} \xi\big|^{\varepsilon _1} C $ for any $k \in \{1, \dots ,m \}$ and $n \in \mathbb{Z} _{\geq 0}$.
Therefore, we obtain (i).

We show (ii).
Let $e'_1, e'_2, \dots ,e'_m$ be the eigenvalues of the matrix $I_m- B_0 -B_1 - \dots - B_N$.
By the assumption, we have $\min \{ |e'_1|, |e'_2|, \dots , |e'_m| \} > \big|q^{\lambda }\big|$.
Hence, the matrix $I_m -B_0 -B_1 - \dots - B_N$ is invertible and there exists $\varepsilon _2 >0$ such that \smash{$|e'_k |^{-1} < |q | ^{\varepsilon _2} \big|q ^{\lambda }\big|^{-1} $} for any $k \in \{1, \dots ,m \}$.
Set $e_k =1/e'_k $ $(k=1,\dots ,m )$ and \smash{$\varepsilon = |q | ^{\varepsilon _2} \big|q ^{\lambda }\big|^{-1} - \max \{ |e_1|, |e_2|, \dots , |e_m| \} (>0)$}.
There exist numbers $p_j \in \{ 0, \varepsilon /2 \}$ $(j=1, \dots , m-1)$ and an invertible matrix $Q_{\infty }$ such that $Q_{\infty }^{-1} ( I_m -B_0 -B_1 - \dots - B_N ) ^{-1} Q_{\infty }$ is expressed as the right-hand side of equation \eqref{eq:Q0-BxQ0}.
Set $p_m=0$ for later use.
Since $B(x) \to I_m -B_0 -B_1 - \dots - B_N $ as $x \to \infty $, the matrix $B(x)$ is invertible if $|x|$ is sufficiently large.
Write
\begin{align*}
 Q_{\infty }^{-1} B\bigl(q^{-1} x\bigr)^{-1} Q_{\infty } = {}& Q_{\infty }^{-1} ( I_m -B_0 -B_1 - \dots - B_N ) ^{-1} Q_{\infty } \\
& + \begin{pmatrix}
b_{1,1} (x) & \cdots & b_{1,m} (x) \\
\vdots & & \vdots \\
 b_{m,1} (x) & \cdots & b_{m,m} (x)
\end{pmatrix}.
\end{align*}
Since $Q_{\infty }^{-1} B\bigl(q^{-1} x\bigr)^{-1} Q_{\infty } \to Q_{\infty }^{-1} ( I_m -B_0 -B_1 - \dots - B_N ) ^{-1} Q_{\infty } $ as $x \to \infty $, there exists $D >0$ such that $|b_{i,j}(x)| < \varepsilon /(2m)$ for any $i,j \in \{1, \dots ,m \}$ and any $x$ $(|x| > D )$.
Let $Y (x)$ be a solution to $Y (qx) = B(x)Y (x)$.
Then, the function $Z(x) = Q_{\infty }^{-1} Y (x)$ satisfies $Z\bigl(q^{-1}x\bigr) = \smash{Q_{\infty }^{-1} B\bigl(q^{-1}x\bigr) ^{-1} Q_{\infty } Z(x)}$.
We assume that $|x| > D $ and \smash{$|[Z(x )]_k| < \widetilde{A}$} for any $k \in \{1, \dots ,m \}$.
It is shown as equation \eqref{eq:Zqxineq} that \smash{$\big| \big[Z\bigl(q^{-1}x \bigr)\big]_k \big| \leq |q | ^{\varepsilon _2} \big|q ^{\lambda }\big|^{-1} \widetilde{A} $} and
\begin{equation}
 | [Z(q^{-n} x )]_k | \leq |q | ^{n \varepsilon _2} \big|q ^{\lambda }\big|^{-n} \widetilde{A} \label{eq:ineq2}
\end{equation}
for $n \in \mathbb{Z} _{\geq 0} $ and $k \in \{1, \dots ,m \}$.

Let $\xi \in \mathbb{C} \setminus \{0\} $.
There exists an integer $M_{\infty }$ such that $| q^n \xi | < D $ for any integer $n \leq -M_{\infty }$.
We take a positive number $C' $ such that $\big| \big[Z\bigl( q^{-M_{\infty }} \xi \bigr)\big]_k\big| \leq C' $ for any $k \in \{1, \dots ,m \}$.
It follows from the inequality \eqref{eq:ineq2} that
\begin{align*}
\big| \big[Z\bigl( q^{-M_{\infty } -n} \xi \bigr)\big]_k\big| &\leq |q | ^{n \varepsilon _2} \big|q ^{\lambda }\big|^{-n} C'\\
& = \big|\bigl(q^{-M_{\infty } -n } \xi \bigr)^{\lambda } \big| \big|q^{-M_{\infty } -n} \xi \big|^{ - \varepsilon _2 }
  \big|\bigl(q^{-M_{\infty } } \xi \bigr)^{- \lambda } \big| \big|q^{-M_{\infty } } \xi \big|^{ \varepsilon _2 } C' \end{align*}
 for any $k \in \{1, \dots ,m \}$ and $n \in \mathbb{Z} _{\geq 0}$.
Since $Y (x) = Q_{\infty } Z(x)$, we obtain that there exists a~positive number $C$ such that
\[
\big| \big[Y\bigl( q^{-M_{\infty } -n} \xi \bigr)\big]_k\big| \leq \big|\bigl(q^{-M_{\infty } -n } \xi \bigr)^{\lambda } \big|
\cdot \big|q^{-M_{\infty } -n} \xi \big|^{ - \varepsilon _2 } C
\]
 for any $k \in \{1, \dots ,m \}$ and $n \in \mathbb{Z} _{\geq 0}$.
Therefore, we obtain (ii).
\end{proof}

By combining Propositions \ref{prop:qcintconv} and \ref{prop:qcintconvYx} with Theorem~\ref{thm:qcint}, we have the following.

\begin{Theorem} \label{thm:suffcondconv}
Set $b_0=0$.
Let $Y(x)$ be a solution to
\begin{equation*}
\frac{Y(q x) - Y(x )}{-x} = \Biggl[ \sum^{N}_{i = 0}\frac{B_{i}}{x -b_{i}} \Biggr] Y(x).
\end{equation*}
Assume that the integrand of $\widetilde{Y}_{i} (x) $ in equation \eqref{eq:qcintadd} is well defined $($i.e., does not diverge$)$ about~${s= q^n \xi }$ for any $i \in \{ 0 , 1, \dots , N \}$ and $n \in \mathbb{Z} $.
If the absolute value of any eigenvalue of~${I_m - B_0}$ is strictly less than $1$ and the absolute value of any eigenvalue of $I_m - B_0 -B_1 - \dots - B_N$ is strictly more than $\big|q^{\lambda }\big|$, then the function $\widetilde{Y}(x)$ in equation \eqref{eq:qcintadd} converges and it satisfies
\begin{equation*}
\frac{\widetilde{Y}(q x) - \widetilde{Y}(x )}{-x} = \Biggl[ \sum^{N}_{i = 0}\frac{G_{i}}{x -b_{i}} \Biggr] \widetilde{Y}(x).
\end{equation*}
\end{Theorem}
\begin{Remark}
It is natural to ask whether these conditions are also necessary.
We expect that the absence of eigenvalues $\nu$ of $I_m - B_0$ satisfying $|\nu| > 1$ and eigenvalues $\nu'$ of $I_m - B_0 -B_1 - \dots - B_N$ satisfying $|\nu'| < \big|q^{\lambda }\big|$ should be closely related to the convergence of $\widetilde{Y}(x )$.
However, the borderline cases, namely when the eigenvalues satisfy $|\nu| = 1$ or $|\nu'| = \big|q^{\lambda }\big|$, require delicate analysis.
For this reason, we focus on sufficient conditions and leave a precise characterization of necessary conditions for future investigation.
\end{Remark}

The definition of the $q$-middle convolution based on the $q$-convolution is quite the same as the one by Sakai and Yamaguchi \cite{SY} (see Definition \ref{def:qmcSY} in the appendix), which is a $q$-analogue of the one by Dettweiler and Reiter \cite{DR1}.
The following proposition, which is proved straightforwardly, is used to define the $q$-middle convolution.
\begin{Proposition}
Define the subspaces $\mathcal{K}$ and $\mathcal{L} $ of $(\mathbb{C}^m )^{N+1}$ as follows:
\begin{equation}
 \mathcal{K} =
 \left(
 \begin{matrix}
 \operatorname{ker}B_0 \\
 \vdots \\
 \operatorname{ker}B_N
 \end{matrix}
 \right), \qquad
 \mathcal{L} = \operatorname{ker}(G_0 +G_1 + \dots + G_N ).
\label{eq:qKL}
\end{equation}
Then, the spaces $\mathcal{K}$ and $\mathcal{L} $ are invariant under the action of $G_k$ for $k=0 ,1, \dots ,N $.
\end{Proposition}
\begin{Definition}[$q$-middle convolution]\label{def:qmc}
We denote the matrix induced from the action of $G_k$ on the quotient space $(\mathbb{C}^m )^{N+1}/(\mathcal{K} + \mathcal{L})$ by $\overline{G}_k$ $(k=0 ,1, \dots ,N )$.
The $q$-middle convolution $mc ^q _\lambda$ is defined by the correspondence $ ( B_{0}, B_{1} ,\dots , B_N ) \mapsto \bigl( \overline{G}_{0}, \overline{G}_{1} ,\dots ,\overline{G}_N \bigr) $.
\end{Definition}

\section[q-hypergeometric equation and generalization]{$\boldsymbol{q}$-hypergeometric equation and generalization} \label{sec:qhypgeomgen}

\subsection[q-hypergeometric equation]{$\boldsymbol{q}$-hypergeometric equation} \label{sec:qconvex1}

The $q$-hypergeometric equation and its solution have been studied from various viewpoints (e.g., see \cite{GR}), and the $q$-convolution was used in \cite{AT}.
We review the approach by the reformulated $q$-convolution.
Let $\mu \in \mathbb{C}$, $\alpha, \beta \in \mathbb{C} \setminus \{0\}$ and set $y(x)=x^{\mu}(\alpha x;q)_{\infty}/(\beta x;q)_{\infty}$.
Then, it satisfies the first-order linear $q$-difference equation
\begin{equation} \label{eq:qeab}
\frac{y(qx) - y(x)}{-x} = \bigg[ \frac{B_0}{x} + \frac{B_1}{x-1/\alpha} \bigg] y(x) ,\qquad B_0= 1 - q^{\mu},\qquad B_1= q^{\mu} \left(1 - \frac{\beta}{\alpha} \right).
\end{equation}
The function $y(x)= x^{\tilde{\mu }} (q/( \beta x) ;q)_{\infty}/(q/(\alpha x) ;q)_{\infty} $ with the condition $q^{\tilde{\mu }} \alpha /\beta = q^{\mu} $ also satisfies equation \eqref{eq:qeab}.
We apply $q$-convolution $c^q_{\lambda}$ to the pair $(B_0, B_1)$ of $1 \times 1$ matrices and set~${c^q_{\lambda}(B_0, B_1) = (G_0, G_1)}$.
The pair $(G_0 , G_1)$ is the $2\times 2$ matrices written as
\begin{align*}
& G_0 = \begin{pmatrix}
 q^{-\lambda}B_0 + 1 -q^{-\lambda} & q^{-\lambda}B_1 \\
 0 & 0
 \end{pmatrix}
 = \begin{pmatrix}
 1 - q^{\mu - \lambda} & q^{\mu - \lambda}(1 - \beta/\alpha) \\
 0 & 0
	 \end{pmatrix}, \\
& G_1 = \begin{pmatrix}
 0 & 0 \\
 q^{-\lambda}B_0 & q^{-\lambda}B_1 + 1 -q^{-\lambda}
 \end{pmatrix}
 = \begin{pmatrix}
 0 & 0 \\
 q^{-\lambda} - q^{\mu - \lambda} & q^{\mu - \lambda}(1 - \beta/\alpha) + 1 -q^{-\lambda}
 \end{pmatrix}.
\end{align*}
Note that $\det ( G_0 +G_1 )= \bigl( 1 -q^{-\lambda}\bigr)\bigl(1- q^{\mu - \lambda} \beta/\alpha \bigr)$, and the $q$-middle convolution coincides with the $q$-convolution (i.e., $ \mathcal{K} + \mathcal{L} =\{ 0 \}$), if $\lambda \neq 0 \neq \mu $, $\alpha \neq \beta $ and $q^{\mu - \lambda} \beta/\alpha \neq 1$.
Write the corresponding $q$-difference equation as
\begin{equation} \label{eq:qhg-sta}
\frac{\widetilde{Y}(qx)-\widetilde{Y}(x)}{-x} = \Biggl[ \frac{G_0}{x} +\frac{G_1}{x-1/\alpha } \Biggr] \widetilde{Y}(x) ,
\qquad \widetilde{Y}(x)=\begin{pmatrix}
 \widetilde{y}_0(x)\\
 \widetilde{y}_1(x)
 \end{pmatrix}.
\end{equation}
Then, the function $\widetilde{y}_0(x)$ satisfies the single second-order equation
\begin{equation} \label{eq:ytil}
\bigl(q^{-\lambda}\beta x-q\bigr)\widetilde{y}_0(x/q) + q^{\lambda-\mu}(\alpha x-q )\widetilde{y}_0(qx) -\bigl\{ \bigl(q^{-\mu}\alpha+\beta\bigr)x- q - q^{\lambda-\mu+1}\bigr\} \widetilde{y}_0(x)=0,
\end{equation}
and equation \eqref{eq:ytil} corresponds to the standard form of $q$-hypergeometric equation \eqref{eq:qhg-sta-sngl} by setting the parameters appropriately.

By the $q$-integral transformation associated with the $q$-convolution, the $q$-integral representation of solutions to equation \eqref{eq:qhg-sta} is obtained as
\begin{equation} \label{sol:Ytil-sta}
\widetilde{Y}(x)
= \begin{pmatrix}
 \widetilde{y}_0(x) \\
 \widetilde{y}_1(x)
 \end{pmatrix}
=
\begin{pmatrix}
 \displaystyle \int_0^{\xi \infty} \frac{K_{\lambda}(x,t)}{t} y(t) {\rm d}_q t \vspace{1mm}\\
 \displaystyle \int_0^{\xi \infty} \frac{K_{\lambda}(x,t)}{t - 1/\alpha} y(t) {\rm d}_q t
\end{pmatrix} ,
\end{equation}
if the assumption of Theorem~\ref{thm:qcint} holds.
We look into a sufficient condition that the $q$-integral representation converges and satisfies equation \eqref{eq:qhg-sta} by applying Theorem~\ref{thm:suffcondconv}.
In this case, the matrices $B_0$ and $B_1$ in Theorem~\ref{thm:suffcondconv} are scalars in equation \eqref{eq:qeab}, and we have ${1 - B_0 = q^{\mu}}$, ${1 - B_0 - B_1 = q^{\mu} \beta /\alpha}$.
Therefore, under the condition that the functions $y(t) K_{\lambda}(x,t)/ t $ and ${ y(t) K_{\lambda}(x,t)/(t-1/\alpha ) }$ are well defined about $t = q^n \xi $ for any $n \in \mathbb{Z} $, a sufficient condition that the function in \eqref{sol:Ytil-sta} converges and satisfies equation \eqref{eq:qhg-sta} is written as~$ \mu >0$ and ${ \big|q^{\lambda - \mu} \alpha /\beta \big|<1} $.
Note that the condition $ \mu >0$ and $ \big|q^{\lambda - \mu} \alpha /\beta \big|<1 $ leads to the condition~${|q^{\mu } |<1}$ and $ {|q^{ \mu} \beta / \alpha | > |q^{\lambda } |}$.

In the case \smash{$K_{\lambda}(x, s) = K^{(1)}_{\lambda}(x, s) $} (see equation \eqref{eq:tP1tP2}) and $y(x)=x^{\mu}(\alpha x;q)_{\infty}/(\beta x;q)_{\infty} $, the condition of well-definedness is $\xi \neq q^n/\beta $ and $\xi \neq x q^n $ for any $n \in \mathbb{Z}$, and the function $\widetilde{y}_0(x)$ in equation \eqref{sol:Ytil-sta} is written as
\begin{equation} \label{eq:y0P1}
\widetilde{y}_0(x) = (1 - q) x^{-\lambda} \sum_{n = -\infty}^{\infty} (q^n \xi)^{\mu} \frac{\bigl(q^{\lambda + n + 1} \xi/x, q^n \xi \alpha ;q\bigr)_{\infty}}{\bigl(q^{n + 1} \xi/x, q^n \xi \beta ;q\bigr)_{\infty}}.
\end{equation}
By specializing the value $\xi$, it can further be expressed as a unilateral $q$-hypergeometric series~${_2 \phi_1}$.
As shown in \cite{AT}, if we set $\xi =1 /\alpha $ in equation \eqref{eq:y0P1}, then it follows from $(q^n;q)_{\infty}=0$ for~${n\in \mathbb{Z}_{\le 0}}$ that
\begin{gather*}
\widetilde{y}_0(x) |_{\xi=1/\alpha }\\
\qquad = (1-q) x^{-\lambda} \sum_{n=1}^{\infty}(q^n/\alpha)^{\mu} \frac{\bigl(q^{\lambda+n+1}/(\alpha x),q^n;q\bigr)_{\infty}}{\bigl(q^{n+1}/(\alpha x),q^n\beta/\alpha;q\bigr)_{\infty}} \\
\qquad= (1-q)\frac{q^{\mu}}{\alpha^{\mu}} \frac{\bigl(q^{\lambda+2}/(\alpha x\bigr), q;q)_{\infty}}{\bigl(q^2/(\alpha x), q\beta/\alpha ;q\bigr)_{\infty}} x^{-\lambda}
\left\{ 1+ \sum_{n=2}^{\infty} \frac{\bigl(q^{2}/(\alpha x);q\bigr) _{n-1} (q \beta/\alpha;q)_{n-1}}{\bigl(q^{\lambda+2}/(\alpha x);q\bigr) _{n-1} (q ;q)_{n-1}} (q^{\mu})^{n-1} \right\} \\
\qquad = (1-q)\frac{q^{\mu}}{\alpha^{\mu}}\frac{\bigl(q^{\lambda+2}/(\alpha x),q;q\bigr)_{\infty}}{\bigl(q^2/(\alpha x),q\beta/\alpha;q\bigr)_{\infty}} x^{-\lambda } {}_2\phi_1 \left( \begin{matrix} q^2/(\alpha x), q\beta/\alpha\\ q^{\lambda+2}/(\alpha x) \end{matrix} ;q,q^{\mu} \right).
\end{gather*}
By applying Heine's transformation formula in \cite[equation~(1.4.1)]{GR},
\[
{}_2\phi_1 (a,b;c;q,z)= \frac{(b,az;q)_{\infty}}{(c,z;q)_{\infty}}  \cdot {}_2\phi_1 (c/b,z;az; q,b), \qquad |z|<1, \quad |b|<1,
\]
 we can rewrite it as
\begin{equation} \label{eq:y0al}
\widetilde{y}_0(x) |_{\xi=1/\alpha }
= (1-q) q^{\mu} \alpha^{-\mu} x^{-\lambda} \frac{\bigl(q, q^{\mu+1}\beta/\alpha ;q\bigr)_{\infty}}{(q\beta/\alpha,q^{\mu} ;q)_{\infty}}
{}_2\phi_1 \left( \begin{matrix} q^{\lambda},q^{\mu}\\
q^{\mu+1}\beta/\alpha \end{matrix} ;q, q^2/(\alpha x) \right),
\end{equation}
which was essentially obtained in version~2 of~\cite{Ar}.
We also have
\begin{align}
 \widetilde{y}_0(x) |_{\xi=q^{-\lambda}x } & = (1-q)q^{-\lambda\mu}x^{\mu -\lambda } \frac{\bigl(q^{-\lambda}\alpha x,q;q\bigr)_{\infty}}{\bigl(q^{-\lambda}\beta x,q^{1-\lambda};q\bigr)_{\infty}} {}_2\phi_1 \left( \begin{matrix} q^{-\lambda}\beta x,q^{1-\lambda}\\ q^{-\lambda}\alpha x \end{matrix} ;q,q^{\mu} \right) \nonumber \\
& = (1-q) q^{-\lambda\mu} x^{\mu - \lambda} \frac{\bigl(q,q^{-\lambda+\mu+1};q\bigr)_{\infty}}{\bigl(q^{-\lambda+1},q^{\mu};q\bigr)_{\infty}}
{}_2\phi_1 \left( \begin{matrix} \alpha/\beta, q^{\mu}\\
q^{-\lambda+\mu+1} \end{matrix} ;q,q^{-\lambda}\beta x \right).  \label{eq:y0la}
\end{align}
The condition for the parameters is $ \mu >0$ in these two cases, and we do not need the condition~${\big|q^{\lambda - \mu} \alpha /\beta \big|<1}$ (see \cite{AT}).
The functions in equations \eqref{eq:y0al} and \eqref{eq:y0la} satisfy the $q$-hypergeo\-metric equation in \eqref{eq:ytil}.

In the case \smash{$K_{\lambda}(x, s) = K^{(2)}_{\lambda}(x, s) $} and $y(x)= x^{\tilde{\mu }} (q/( \beta x) ;q)_{\infty}/(q/(\alpha x) ;q)_{\infty} $, the condition of well-definedness is $\xi \neq q^n/\alpha $ and $\xi \neq x q^{-\lambda +n } $ for any $n \in \mathbb{Z}$, and the function $\widetilde{y}_0(x)$ is written~as
\begin{gather} \label{eq:y0P2}
\widetilde{y}_0(x) = (1-q) \sum_{n=-\infty}^{\infty}(q^n\xi)^{\tilde{\mu } -\lambda } \frac{\bigl(x q^{-n} / \xi , q^{1-n} / (\beta \xi ) ; q\bigr)_{\infty}}{\bigl(x q^{-\lambda -n} / \xi , q^{1-n} / (\alpha \xi ) ; q\bigr)_{\infty}}.
\end{gather}
We set $\xi=1/\beta $ or $\xi= x $ in equation \eqref{eq:y0P2}.
It was essentially established in version~2 of~\cite{Ar} and~\cite{AT} that
\begin{gather}
 \widetilde{y}_0(x) |_{\xi=1/\beta }
 = (1-q) \beta^{\lambda-\tilde{\mu }} \frac{(\beta x, q ;q)_{\infty}}{\bigl(q^{-\lambda}\beta x, q\beta/\alpha ;q\bigr)_{\infty}}
{}_2\phi_1 \left( \begin{matrix} q^{-\lambda}\beta x,q\beta/\alpha\\
\beta x \end{matrix} ;q,q^{\lambda-\mu}\frac{\alpha}{\beta} \right)\nonumber \\
\phantom{ \widetilde{y}_0(x) |_{\xi=1/\beta }
 }{} = (1-q) \beta^{\lambda-\tilde{\mu }} \frac{\bigl(q, q^{\lambda-\mu+1} ;q\bigr)_{\infty}}{\bigl(q\beta/\alpha, q^{\lambda-\mu}\alpha/\beta ;q\bigr)_{\infty}}
{}_2\phi_1 \left( \begin{matrix} q^{\lambda},q^{\lambda-\mu}\alpha/\beta\\
 q^{\lambda-\mu+1} \end{matrix} ;q,q^{-\lambda}\beta x \right), \nonumber\\
 \widetilde{y}_0(x) |_{\xi= x }
 = (1-q) q^{\lambda - \mu} \frac{\alpha}{\beta} x^{\tilde{\mu } - \lambda}\frac{\bigl(q^2/(\beta x), q ;q\bigr)_{\infty}}{\bigl(q^2/(\alpha x),q^{-\lambda+1} ;q\bigr)_{\infty}}
{}_2\phi_1 \left( \begin{matrix} q^{-\lambda+1},q^2/(\alpha x)\\
q^2/(\beta x) \end{matrix} ;q,q^{\lambda-\mu}\frac{\alpha}{\beta} \right)\nonumber \\
\phantom{ \widetilde{y}_0(x) |_{\xi= x }
 }{} = (1-q) q^{\lambda-\mu} \frac{\alpha}{\beta} x^{\tilde{\mu } - \lambda} \frac{\bigl(q, q^{-\mu+1}\alpha/\beta ;q\bigr)_{\infty}}{\bigl(q^{-\lambda+1}, q^{\lambda-\mu}\alpha/\beta ;q\bigr)_{\infty}}
{}_2\phi_1\nonumber \\
\phantom{ \widetilde{y}_0(x) |_{\xi= x }
= }{}
\times\left( \begin{matrix} \alpha/\beta, q^{\lambda-\mu}\alpha/\beta\\
q^{-\mu+1}\alpha/\beta \end{matrix} ;q, q^2/(\alpha x) \right). \label{eq:y0be}
\end{gather}
The condition for the parameters is $|q^{\lambda - \mu} \alpha /\beta |<1$ in these two cases, and we do not need the condition $ \mu >0$ (see \cite{AT}).
The functions in equation \eqref{eq:y0be} satisfy the $q$-hypergeometric equation in \eqref{eq:ytil}.

\subsection[Generalized q-hypergeometric equation of order 3]{Generalized $\boldsymbol{q}$-hypergeometric equation of order 3} \label{sec:qmcex}

To obtain the generalized $q$-hypergeometric equation of order $3$, we use an addition and a $q$-middle convolution.
We apply the addition ${\rm add}_{\mu' } $ to the pair $(G_0, G_1) $ in equation \eqref{eq:qhg-sta}.
Then, we have the pair $(B'_0 , B'_1 )$, where
\begin{align*}
& B'_0 = \begin{pmatrix}
 1 -q^{\mu' + \mu -\lambda } & q^{\mu' + \mu -\lambda }(1-\beta/\alpha) \\
 0 & 1 - q^{\mu' } \\
 \end{pmatrix}, \\
& B'_1 = \begin{pmatrix}
 0 & 0\\
 q^{\mu' -\lambda } ( 1 -q^{\mu } ) & q^{\mu' } \bigl( q^{\mu -\lambda }(1-\beta/\alpha)+ 1 - q^{-\lambda} \bigr)
\end{pmatrix}.
\end{align*}
The mapping $(G_0, G_1) \mapsto (B'_0 , B'_1 ) $ corresponds to a gauge transformation of the $q$-difference equations.
Let $\widetilde{Y}(x)$ be a solution to equation \eqref{eq:qhg-sta}, and set
\begin{equation} \label{eq:Ygxxmu'}
 Y_g(x) = x^{\mu'}\widetilde{Y}(x).
\end{equation}
The function $ Y_g(x) $ satisfies
\begin{equation} \label{eq:B'0B'1m}
\frac{Y_g(qx) - Y_g(x)}{-x} = \bigg[ \frac{B'_0}{x} + \frac{B'_1}{x-1/\alpha} \bigg] Y_g(x).
\end{equation}
It follow from equation \eqref{eq:Ygxxmu'} and the argument around equation \eqref{sol:Ytil-sta} that the function
 \begin{equation*}
Y_g(x) =
\begin{pmatrix}
 \displaystyle x^{\mu'} \int_0^{\xi \infty} \frac{K_{\lambda}(x,t)}{t} y(t) {\rm d}_q t \vspace{1mm}\\
 \displaystyle x^{\mu'} \int_0^{\xi \infty} \frac{K_{\lambda}(x,t)}{t-1/\alpha} y(t) {\rm d}_q t
\end{pmatrix},
\end{equation*}
satisfies equation \eqref{eq:B'0B'1m}, if the function $y(x)$ is a solution to equation \eqref{eq:qeab}, $ \mu >0$, $ \big|q^{\lambda - \mu} \alpha / \beta \big|<1 $ and the functions $y(t) K_{\lambda}(x,t)/ t $ and $ y(t) K_{\lambda}(x,t)/(t-1/\alpha ) $ are well defined about $t = q^n \xi $ for any $n \in \mathbb{Z} $.

Before considering the $q$-middle convolution, we apply the $q$-convolution $c^q_{\lambda'}$ to $(B_0', B_1')$.
We obtain $c^q_{\lambda'} (B_0', B_1')= (G'_0, G'_1)$, where $G'_0$ and $G'_1$ are the $4 \times 4$ matrices described as
\begin{gather*}
 G'_0
= \begin{pmatrix}
 q^{-\lambda'}B'_0 + \bigl(1 - q^{-\lambda'}\bigr)I_2 & q^{-\lambda'}B'_1 \\
 O & O
\end{pmatrix},\qquad
G'_1
= \begin{pmatrix}
 O & O \\
 q^{-\lambda'}B'_0 & q^{-\lambda'}B'_1 + \bigl(1-q^{-\lambda'}\bigr) I_2
\end{pmatrix}.
\end{gather*}
Note that \smash{$\det ( G'_0 + G'_1 ) = \bigl(1- q^{-\lambda'} \bigr)^2 \bigl(1-q^{\mu'-\lambda-\lambda'}\bigr) \bigl(1-q^{\mu + \mu' -\lambda-\lambda'} \beta /\alpha \bigr)$}.
We investigate a~sufficient condition for the parameters that the $q$-difference equation
\begin{equation} \label{eq:G'0G'1m}
\frac{\widetilde{Y_g}(qx) - \widetilde{Y_g}(x)}{-x} = \bigg[ \frac{G'_0}{x} + \frac{G'_1}{x-1/\alpha} \bigg] \widetilde{Y_g} (x)
\end{equation}
admits solutions in terms of the $q$-integral representation by applying Theorem~\ref{thm:suffcondconv} for $c^q_{\lambda'} (B_0', B_1')\allowbreak= (G'_0, G'_1)$.
In this case, the condition for the eigenvalues is applied for the matrices $I_2 - B'_0 $ and~${I_2 - B'_0 - B'_1 }$.
The eigenvalues of $I_2 - B'_0 $ are $ q^{\mu'}$, $ q^{\mu+\mu'-\lambda} $, and the eigenvalues of $I_2 - B'_0 - B'_1 $ are $ q^{\mu' - \lambda }$, $ q^{\mu+\mu'-\lambda} \beta/ \alpha $.
The formal $q$-integral representation of solutions to the equation obtained after applying $q$-convolution $c_{\lambda'}^q$ is given by
\begin{align}
\widetilde{Y_g}(x)
& = \begin{pmatrix}
 \displaystyle \int_0^{\xi' \infty} \frac{\widehat{K}_{\lambda'}(x,s)}{s} Y_g (s) {\rm d}_q s \vspace{1mm}\\
 \displaystyle \int_0^{\xi' \infty} \frac{\widehat{K}_{\lambda'}(x,s)}{s - 1/\alpha} Y_g (s) {\rm d}_q s
 \end{pmatrix} \nonumber \\
& = \begin{pmatrix}
 \displaystyle \int_0^{\xi' \infty} \frac{\widehat{K}_{\lambda'}(x,s)}{s} s^{\mu' } \int_0^{\xi \infty} \frac{K_{\lambda}(s,t)}{t} y(t) {\rm d}_q t {\rm d}_q s \vspace{1mm}\\
 \displaystyle \int_0^{\xi' \infty} \frac{\widehat{K}_{\lambda'}(x,s)}{s} s^{\mu' } \int_0^{\xi \infty} \frac{K_{\lambda}(s,t)}{t - 1/\alpha} y(t) {\rm d}_q t {\rm d}_q s \vspace{1mm}\\
 \displaystyle \int_0^{\xi' \infty} \frac{\widehat{K}_{\lambda'}(x,s)}{s - 1/\alpha} s^{\mu' } \int_0^{\xi \infty} \frac{K_{\lambda}(s,t)}{t} y(t) {\rm d}_q t {\rm d}_q s \vspace{1mm}\\
 \displaystyle \int_0^{\xi' \infty} \frac{\widehat{K}_{\lambda'}(x,s)}{s - 1/\alpha} s^{\mu' } \int_0^{\xi \infty} \frac{K_{\lambda}(s,t)}{t - 1/\alpha} y(t) {\rm d}_q t {\rm d}_q s
 \end{pmatrix} , \label{eq:tYg}
\end{align}
where $\widehat{K}_{\lambda }(x,s) $ and $K_{\lambda }(x,s)$ are functions which satisfy equation \eqref{eq:tPlambda}, and $y(x)$ is a function which satisfies equation~\eqref{eq:qeab}.
Then, we obtain the following proposition by Theorem~\ref{thm:suffcondconv}.
\begin{Proposition} \label{prop:ghge3conv}
If the parameters satisfy
\begin{gather}
 \mu >0, \qquad \mu ' >0, \qquad \mu ' +\mu -\lambda >0 , \qquad \lambda ' + \lambda - \mu' >0 ,\qquad \big|q^{\lambda - \mu} \alpha /\beta \big|<1 , \nonumber\\ \big|q^{\lambda '+ \lambda - \mu ' - \mu} \alpha /\beta \big|<1,\label{eq:prop41}
\end{gather}
and the integrands in equation \eqref{eq:tYg} are well defined about $(s,t) = (q^{n'} \xi ', q^n \xi )$ for any $n', n \in \mathbb{Z} $, then the function \smash{$\widetilde{Y_g}(x)$} in equation \eqref{eq:tYg} converges and satisfies equation \eqref{eq:G'0G'1m}.
\end{Proposition}
Note that the condition $\mu >0$, $\mu ' >0$, $\mu ' +\mu -\lambda >0$ and $\lambda ' + \lambda - \mu' >0$ leads to the condition $|q^{\mu } |<1$, $|q^{\mu '} |<1$, $\big|q^{\mu ' +\mu -\lambda } \big|<1$ and $ \big|q^{ \mu ' - \lambda } \big| > \big|q^{\lambda '} \big| $.

We look into the $q$-middle convolution $mc^q_{\lambda'}$ to $(B_0', B_1')$.
Recall that the $q$-middle convolution is obtained as the pair of the operators on the quotient space $\mathbb{C}^4/(\mathcal{K} + \mathcal{L}) $.
The vector ${}^t\bigl(q^{\mu}(1-\beta/\alpha)+q^{\lambda}-1, q^{\mu}-1\bigr)$ belongs to the space $\ker (B_1')$, and we have $\dim \mathcal{K}=1$ under the condition ${\mu \neq 0}$, $ \mu ' \neq 0$ and $\mu+\mu'-\lambda \neq 0 $.
We have $\dim \mathcal{L}=0$, if $\lambda ' \neq 0$, $\mu'-\lambda-\lambda' \neq 0$ and $q^{\mu + \mu' -\lambda-\lambda'} \beta /\alpha \neq 1 $.
We continue the argument under these inequalities.
We have $\dim \bigl(\mathbb{C}^4/(\mathcal{K} + \mathcal{L}) \bigr) =3 $.
To obtain matrix representations of $mc^q_{\lambda'} (B_0', B_1') $, we consider simultaneous transformation for the matrices $G'_0 $ and $G' _1$.
Set
\begin{equation*} 
P = \begin{pmatrix}
 1 & 0 & 0 & 0 \\
 0 & 1 & 0 & 0 \\
 0 & 0 & 1 & f_0 \\
 0 & 0 & 0 & q^{\mu}-1
\end{pmatrix},\qquad f_0 = q^{\mu}(1 - \beta/\alpha)+q^{\lambda}-1.
\end{equation*}
Then, $P$ is invertible and we have
\begin{align*}
& P^{-1}G'_0 P = \begin{pmatrix}
 1-q^{\mu+\mu'-\lambda-\lambda'} & q^{\mu+\mu'-\lambda-\lambda'}(1 - \beta/\alpha) & 0 & 0 \\
 0 & 1-q^{\mu'-\lambda'} & q^{\mu'-\lambda-\lambda'} (1 - q^{\mu })& 0 \\
 0 & 0 & 0 & 0 \\
 0 & 0 & 0 & 0
 \end{pmatrix}, \\
& P^{-1}G'_1 P
= \begin{pmatrix}
 0 & 0 & 0 & 0 \\
 0 & 0 & 0 & 0 \\
 f_1 & f_2 & f_3 & 0 \\
 0 & -q^{-\lambda '}\frac{q^{\mu '}-1}{q^{\mu } -1} & -q^{\mu ' - \lambda ' - \lambda } & 1-q^{-\lambda '}
 \end{pmatrix},
\end{align*}
where
\begin{gather*}
f_1 = q^{-\lambda '} - q^{\mu ' +\mu - \lambda '- \lambda}, \qquad f_2 = q^{\mu ' +\mu - \lambda '- \lambda} (1 - \beta/\alpha) + q^{-\lambda'} \frac{q^{\mu '} -1}{q^{\mu } -1 } f_0,\qquad\text{and} \\ f_3 = 1 - q^{-\lambda'} + q^{\mu'-\lambda-\lambda'} f_0.
\end{gather*}
The upper-left $3\times 3$ submatrices of $P^{-1}G'_0 P$ and $P^{-1}G'_1 P$ are matrix representations of $G'_0$ and $G'_1$ on the quotient space $\mathbb{C}^4/(\mathcal{K} + \mathcal{L}) $.
Thus, we can write $mc^q_{\lambda'}(B'_0, B'_1) = \bigl(\overline{G'}_0, \overline{G'}_1\bigr)$, where
\begin{align*}
&\overline{G'}_0
= \begin{pmatrix}
 1-q^{\mu+\mu'-\lambda-\lambda'} & q^{\mu+\mu'-\lambda-\lambda'}(1 - \beta/\alpha) & 0 \\
 0 & 1-q^{\mu'-\lambda'} & q^{\mu'-\lambda-\lambda'} (1 - q^{\mu }) \\
 0 & 0 & 0
 \end{pmatrix}, \\
&\overline{G'}_1
= \begin{pmatrix}
 0 & 0 & 0 \\
 0 & 0 & 0 \\
 f_1 & f_2 & f_3
 \end{pmatrix}.
\end{align*}
The equation obtained after applying $q$-middle convolution $mc^q_{\lambda'}$ is written as
\begin{equation} \label{eq:111-21-111}
\frac{1}{-x}
\begin{pmatrix}
 g_1(qx) - g_1(x) \\
 g_2(qx) - g_2(x) \\
 g_3(qx) - g_3(x)
\end{pmatrix}
= \left[ \frac{\overline{G'}_0}{x} + \frac{\overline{G'}_1}{x-1/\alpha} \right]
 \begin{pmatrix}
 g_1(x) \\
 g_2(x) \\
 g_3(x)
 \end{pmatrix}.
\end{equation}
From equation \eqref{eq:111-21-111}, the single third-order $q$-difference equation for $g_1(x)$ is derived as
\begin{gather}
(\alpha x-1) g_1\bigl(q^3 x\bigr) + q^{-\lambda-\lambda'}
\bigl\{ -\bigl(q^{\mu+\mu'}\beta + q^{\mu'}\alpha + q^{\lambda}\alpha\bigr)x + q^{\mu+\mu'} + q^{\lambda+\mu'} + q^{\lambda+\lambda'} \bigr\} g_1\bigl(q^2 x\bigr) \nonumber \\
\qquad {}+q^{\mu'-2\lambda-2\lambda'}
\bigl\{\bigl(q^{\lambda}\alpha + q^{\mu+\mu'}\beta + q^{\lambda+\mu}\beta\bigr)x -q^{\lambda}\bigl(q^{\lambda+\lambda'} + q^{\mu+\mu'} + q^{\mu+\lambda'}\bigr) \bigr\} g_1(qx) \nonumber \\
\qquad {}-q^{\mu+2\mu'-\lambda-2\lambda'}\bigl(q^{-\lambda-\lambda'}\beta x - 1\bigr)g_1(x)
=0. \label{eq:111-21-111-g1}
\end{gather}
Each coefficient of $g_1 \bigl(q^j x\bigr)$ $(j=0,1,2,3)$ is a linear polynomial in $x$.
This is the generalized $q$-hypergeometric equation of order $3$.
The functions $g_2(x)$ and $g_3(x)$ also satisfy the generalized $q$-hypergeometric equation of order $3$, whose coefficients are of degree $1$ and some terms are different from equation \eqref{eq:111-21-111-g1}.

We investigate integral representations of solutions to equation \eqref{eq:111-21-111-g1}.
From now on, we assume the assumption of Proposition~\ref{prop:ghge3conv}.
If \smash{$\widetilde{Y_g}(x) $} is a solution to equation \eqref{eq:G'0G'1m} and we write
\begin{equation*}
\begin{pmatrix}
 g_1(x) \\
 g_2(x) \\
 g_3(x) \\
 g_4(x)
\end{pmatrix}
= P^{-1} \widetilde{Y_g}(x)
=
\begin{pmatrix}
 1 & 0 & 0 & 0 \\
 0 & 1 & 0 & 0 \\
 0 & 0 & 1 & -f_0/(q^{\mu} -1) \\
 0 & 0 & 0 & 1/(q^{\mu} -1)
 \end{pmatrix} \widetilde{Y_g}(x) ,
\end{equation*}
then the function ${}^t(g_1 (x), g_2 (x), g_3 (x) ) $ satisfies equation \eqref{eq:111-21-111}.
In particular, we obtain integral representations of solutions to equation \eqref{eq:111-21-111} by setting $\widetilde{Y_g}(x) $ as equation \eqref{eq:tYg}.

If we specialize to the case
\begin{equation}
\label{eq:KlpKly}
\widehat{K}_{\lambda'}(x,s) = K^{(1)}_{\lambda '}(x, s) ,\qquad K_{\lambda}(s, t) = K^{(1)}_{\lambda}(s, t) ,\qquad y(t)=t^{\mu}\frac{(\alpha t;q)_{\infty}}{(\beta t;q)_{\infty}} ,
\end{equation}
then we obtain
\begin{align}
g_1(x)
={}& x^{-\lambda'} \int_0^{\xi' \infty} \frac{\bigl(q^{\lambda '+1}s/x;q\bigr)_{\infty}}{(qs/x;q)_{\infty}} s^{\mu' - \lambda -1} \int_0^{\xi \infty} t^{\mu -1} \frac{\bigl(q^{\lambda+1}t/s, \alpha t ;q\bigr)_{\infty}}{(qt/s , \beta t ;q)_{\infty}} {\rm d}_q t {\rm d}_q s \nonumber \\
={}& (1 - q)^2 \xi^{\mu} \xi'^{\mu' - \lambda} \nonumber \\
& \times x^{-\lambda'} \sum_{n = -\infty}^{\infty} \sum_{m = -\infty}^{\infty} \bigl(q^{\mu' - \lambda}\bigr)^n (q^{\mu})^m \frac{\bigl(q^{\lambda' + n + 1}\xi'/x, q^{\lambda + m - n + 1}\xi/\xi', q^m \xi \alpha; q\bigr)_{\infty}}{\bigl(q^{n + 1}\xi'/x, q^{m - n + 1}\xi/\xi', q^m\xi \beta; q\bigr)_{\infty}}.  \label{sol:g1-11121111}
\end{align}
This function is a solution to equation \eqref{eq:111-21-111-g1}.
By specializing the values $\xi$ and $\xi'$ and using the $q$-binomial theorem, the solution \eqref{sol:g1-11121111} can be transformed to the $q$-hypergeometric series ${}_3\phi_2$.
\begin{Proposition}
If the parameters satisfy equation \eqref{eq:prop41}, $\lambda \not \in {\mathbb Z}$, $\lambda ' \not \in {\mathbb Z}$ and $\big|q^{ -\lambda ' -\lambda } \beta x\big| <1 $, then
\begin{align*}
g_1(x) | _{\xi = q^{ -\lambda ' -\lambda } x , \xi' = q ^{-\lambda ' }x }
={}& (1 - q)^2 q^{ -\lambda ' \mu -\lambda \mu -\lambda ' \mu' + \lambda ' \lambda} \frac{\bigl(q, q, q^{\mu -\lambda + 1} , q^{\mu' +\mu - \lambda -\lambda' + 1 } ; q\bigr)_{\infty}}{\bigl(q^{1-\lambda' } , q^{1- \lambda } , q^{\mu } , q^{\mu' +\mu - \lambda } ; q\bigr)_{\infty}} \\
&\times x^{\mu' + \mu -\lambda' -\lambda } {}_3\phi_2 \left( { \alpha /\beta , q^{\mu }, q^{\mu' +\mu - \lambda } \atop{q^{\mu -\lambda + 1} , q^{\mu' +\mu - \lambda -\lambda' + 1 }}} ; q, q^{ -\lambda ' -\lambda } \beta x \right).
\end{align*}
\end{Proposition}
\begin{proof}
Set $\xi ' = q ^{-\lambda ' }x $ and $\xi = q^{ -\lambda ' -\lambda } x$ in equation \eqref{sol:g1-11121111}.
We change the variables $(m,n)$ in equation \eqref{sol:g1-11121111} to the variables $(n,k)$ by setting $m-n = k$.
It follows from $\bigl(q^{\ell +1};q\bigr)_{\infty} = 0$ for~${\ell \in \mathbb{Z}_{<0}}$ that
\begin{gather}
 \sum_{n = -\infty}^{\infty} \sum_{m = -\infty}^{\infty} \bigl(q^{\mu' - \lambda}\bigr)^n (q^{\mu})^m \frac{\bigl(q^{\lambda' + n + 1}\xi'/x, q^{\lambda + m - n + 1}\xi/\xi', q^m \xi \alpha; q\bigr)_{\infty}}{\bigl(q^{n + 1}\xi'/x, q^{m - n + 1}\xi/\xi', q^m\xi \beta; q\bigr)_{\infty}} \nonumber \\
\qquad = \sum_{n = 0}^{\infty} \sum_{k = 0}^{\infty} \bigl(q^{\mu' - \lambda}\bigr)^n (q^{\mu})^{n+ k} \frac{\bigl(q^{n + 1}, q^{k + 1} , q^{ -\lambda ' -\lambda +n+k} \alpha x ; q\bigr)_{\infty}}{\bigl(q^{-\lambda' + n + 1} , q^{- \lambda + k + 1}, q^{ -\lambda ' -\lambda +n+k} \beta x ; q\bigr)_{\infty}}.  \label{eq:prop1}
\end{gather}
Note that the summand in equation \eqref{eq:prop1} is well defined, if $\lambda \not \in {\mathbb Z}$ and $\lambda ' \not \in {\mathbb Z}$.
By the $q$-binomial theorem, we have
\begin{equation*}
\frac{\bigl(q^{ -\lambda ' -\lambda +n+k} \alpha x; q\bigr)_{\infty}}{\bigl(q^{ -\lambda ' -\lambda +n+k} \beta x; q\bigr)_{\infty}}
= \sum_{N=0}^{\infty} \frac{(\alpha /\beta ; q)_N}{(q; q)_N} \bigl(q^{ -\lambda ' -\lambda +n+k} \beta x\bigr)^N ,
\end{equation*}
if $\big|q^{ -\lambda ' -\lambda } \beta x \big| <1$.
Hence, equation \eqref{eq:prop1} is equal to
\begin{gather}
\sum_{n = 0}^{\infty} \sum_{k = 0}^{\infty} \sum_{N=0}^{\infty} \bigl(q^{\mu' +\mu - \lambda +N }\bigr)^n \bigl(q^{\mu +N} \bigr)^ {k} \frac{\bigl(q^{n + 1}, q^{k + 1} ; q\bigr)_{\infty}}{\bigl(q^{-\lambda' + n + 1} , q^{- \lambda + k + 1} ; q\bigr)_{\infty}} \frac{(\alpha /\beta ; q)_N}{(q; q)_N} \bigl(q^{ -\lambda ' -\lambda } \beta x\bigr)^N \nonumber\\
\qquad= \frac{(q, q; q)_{\infty}}{\bigl(q^{-\lambda' + 1} , q^{- \lambda + 1} ; q\bigr)_{\infty}} \sum_{n = 0}^{\infty} \sum_{k = 0}^{\infty} \sum_{N=0}^{\infty} \bigl(q^{\mu' +\mu - \lambda +N }\bigr)^n \frac{\bigl(q^{-\lambda' + 1}; q\bigr)_n}{(q; q)_n} \bigl(q^{\mu +N} \bigr)^ {k} \frac{\bigl(q^{-\lambda + 1} ; q\bigr)_k}{(q; q)_k} \nonumber \\
\phantom{\qquad=}{} \times \frac{(\alpha /\beta ; q)_N}{(q; q)_N} \bigl(q^{ -\lambda ' -\lambda } \beta x\bigr)^N.  \label{eq:prop2}
\end{gather}
It follow from the conditions $ \mu' +\mu - \lambda >0$ and $\mu >0 $ and the $q$-binomial theorem that equation~\eqref{eq:prop2} is equal to
\begin{gather*}
 \frac{(q, q; q)_{\infty}}{\bigl(q^{-\lambda' + 1} , q^{- \lambda + 1} ; q\bigr)_{\infty}} \sum_{N=0}^{\infty}
\frac{ \bigl(q^{\mu' +\mu - \lambda -\lambda' + 1 +N } ; q\bigr)_{\infty} }{\bigl(q^{\mu' +\mu - \lambda +N } ; q\bigr)_{\infty}} \frac{ \bigl(q^{\mu -\lambda + 1 + N} ; q\bigr)_{\infty} }{\bigl(q^{\mu +N} ; q\bigr)_{\infty}} \frac{(\alpha /\beta ; q)_N}{(q; q)_N} \bigl(q^{ -\lambda ' -\lambda } \beta x\bigr)^N \\
\qquad = \frac{\bigl(q, q, q^{\mu -\lambda + 1} , q^{\mu' +\mu - \lambda -\lambda' + 1 } ; q\bigr)_{\infty}}{\bigl(q^{-\lambda' + 1} , q^{- \lambda + 1} , q^{\mu } , q^{\mu' +\mu - \lambda } ; q\bigr)_{\infty}} \sum_{N=0}^{\infty} \frac{\bigl(\alpha /\beta , q^{\mu }, q^{\mu' +\mu - \lambda } ; q\bigr)_N}{\bigl(q, q^{\mu -\lambda + 1} , q^{\mu' +\mu - \lambda -\lambda' + 1 } ; q\bigr)_N} \bigl(q^{ -\lambda ' -\lambda } \beta x\bigr)^N.
\end{gather*}
By combining with equation \eqref{sol:g1-11121111}, we obtain the proposition.
\end{proof}

If we choose $\widehat{K}_{\lambda'}(x, s)$, $K_{\lambda}(s, t)$ and $y(t)$ differently from equation \eqref{eq:KlpKly}, then we obtain a~solution different from equation \eqref{sol:g1-11121111}.
Furthermore, by suitably specializing the values $\xi$ and~$\xi'$, it can be transformed into the $q$-hypergeometric series ${}_3\phi_2$.

\subsection[Another generalized q-hypergeometric equation]{Another generalized $\boldsymbol{q}$-hypergeometric equation}

To obtain the $q$-difference equation which is different from the generalized $q$-hypergeometric equation in the previous subsection, we apply another operation to equation \eqref{eq:qhg-sta}.
Let $\widetilde{Y}(x)$ be a solution to equation \eqref{eq:qhg-sta}, and set
\begin{equation}
 Y_g(x) = \frac{(\gamma x;q)_{\infty}}{(\alpha x;q)_{\infty}}\widetilde{Y}(x).
\label{eq:Yggaal}
\end{equation}
As we discussed in Section~\ref{sec:add}, the $q$-difference equation which the function $Y_g(x) $ satisfies is written as
\begin{align}
& \frac{Y_g(qx) - Y_g(x)}{-x} = \bigg[ \frac{B'_0}{x} + \frac{B'_1}{x-1/\gamma} \bigg] Y_g(x) , \label{eq:B'0B'1ag} \\
&  B'_0
= q^{-\lambda}
 \begin{pmatrix}
 q^{\lambda} - q^{\mu} & q^{\mu} (1-\beta /\alpha ) \\
 0 & 0
 \end{pmatrix}
=
 \begin{pmatrix}
 1 - q^{\mu -\lambda } & q^{\mu -\lambda } (1- \beta/ \alpha ) \\
 0 & 0
 \end{pmatrix}, \nonumber\\
&  B'_1 = \begin{pmatrix}
 q^{\mu-\lambda}(1-\alpha /\gamma ) & -q^{\mu-\lambda} (1-\beta / \alpha )(1-\alpha /\gamma ) \\
 q^{-\lambda} \bigl(1-q^{\mu}\bigr)\alpha /\gamma & q^{-\lambda}\bigl\{q^{\mu}(\alpha -\beta)/\gamma +q^{\lambda}-\alpha / \gamma \bigr\}
\end{pmatrix}. \nonumber
\end{align}
Note that $\det B'_1 = q^{\mu-\lambda} (1-\alpha/ \gamma ) \bigl(1- q^{-\lambda} \beta /\gamma \bigr) $.
We apply the $q$-convolution $c^q_{\lambda'}$ to $(B_0', B_1')$.
Then, we obtain $c^q_{\lambda'} (B_0', B_1')= (G'_0, G'_1)$, where $G'_0$ and $G'_1$ are the $4 \times 4$ matrices described as
\begin{gather}
G'_0
= \begin{pmatrix}
 q^{-\lambda'}B'_0 + \bigl(1 - q^{-\lambda'}\bigr)I_2 & q^{-\lambda'}B'_1 \\
 O & O
\end{pmatrix} ,\nonumber\\
G'_1
 = \begin{pmatrix}
 O & O \\
 q^{-\lambda'}B'_0 & q^{-\lambda'}B'_1 + \bigl(1-q^{-\lambda'}\bigr) I_2
\end{pmatrix}. \label{eq:G'0G'1matman}
\end{gather}
Note that $\det ( G'_0 + G'_1 ) = \bigl(1- q^{-\lambda'}\bigr)^2 \bigl(1- q^{-\lambda'-\lambda} \alpha/\gamma\bigr) \bigl(1- q^{-\lambda'-\lambda +\mu } \beta /\gamma \bigr)$.

We now look into the $q$-middle convolution $mc^q_{\lambda'}$ to $(B_0', B_1')$.
The vector ${}^t\bigl(q^{\mu}(1-\beta / \alpha ), \allowbreak q^{\mu}-q^{\lambda}\bigr)$ belongs to the space $\ker (B_0 ')$, and we have $\dim \mathcal{K}=1$ under the condition $ \beta \neq \alpha \neq \gamma $ and~${\beta \neq q^{\lambda} \gamma }$.
We have $\dim \mathcal{L}=0$, if $\lambda ' \neq 0$, $q^{-\lambda'-\lambda} \alpha/\gamma \neq 1$ and $q^{-\lambda'-\lambda +\mu } \beta /\gamma \neq 1 $.
We continue the argument under these inequalities.
We have $\dim \bigl(\mathbb{C}^4/(\mathcal{K} + \mathcal{L}) \bigr) =3 $.
By the simultaneous transformation for $(G'_0, G'_1)$ by the matrix
\begin{equation*}
P = \begin{pmatrix}
 0 & 0 & 0 & q^{\mu}(1-\beta /\alpha ) \\
 0 & 1 & 0 & q^{\mu}-q^{\lambda} \\
 0 & 0 & 1 & 0 \\
 1 & 0 & 0 & 0
\end{pmatrix},
\end{equation*}
we can extract the $3 \times 3$ matrices $\overline{G'}_0$ and $\overline{G'}_1$, where $mc^q_{\lambda'}(B'_0, B'_1) = \bigl(\overline{G'}_0, \overline{G'}_1\bigr)$.
Namely, the upper-left $3\times 3$ submatrices of $P^{-1}G'_0 P$ and $P^{-1}G'_1 P$ are matrix representations of $G'_0$ and $G'_1$ on the quotient space $\mathbb{C}^4/(\mathcal{K} + \mathcal{L}) $, and they are expressed as
\begin{equation*}
\overline{G'}_0
= \begin{pmatrix}
 0 & 0 & 0 \\
 h_1 & h_2 & h_3 \\
 0 & 0 & 0
 \end{pmatrix}, \qquad
\overline{G'}_1
= \begin{pmatrix}
 j_1 & 0 & j_2 \\
 0 & 0 & 0 \\
 j_3 & j_4 & j_5
 \end{pmatrix},
\end{equation*}
where
\begin{align*}
& h_1 = q^{-\lambda-\lambda'} \bigl\{ q^{\mu} (1- \beta /\gamma ) + \bigl(q^{\lambda} -1\bigr) \alpha /\gamma \bigr\} , \qquad h_2 = 1-q^{\mu-\lambda-\lambda'}, \\
& h_3 = \frac{q^{\mu}-q^{\lambda}+\bigl\{ \bigl(q^{\lambda}-1\bigr)\alpha +\bigl(1-q^{\mu}\bigr)\beta \bigr\} /\gamma }{q^{\lambda+\lambda'}(\beta /\alpha -1)}, \qquad j_1 = 1 + q^{-\lambda-\lambda'} \{q^{\mu}( \alpha -\beta) - \alpha \} /\gamma , \\
& j_2 = q^{-\lambda-\lambda'}\bigl(1-q^{\mu}\bigr)\alpha /\gamma ,\qquad j_3 = -q^{\mu-\lambda-\lambda'} (1- \beta /\alpha )(1-\alpha / \gamma ), \\
& j_4 = q^{\mu-\lambda-\lambda'}(1- \beta /\alpha ), \qquad j_5 = 1 + q^{-\lambda-\lambda'}\bigl\{ q^{\mu}(1-\alpha/ \gamma )-q^{\lambda} \bigr\}.
\end{align*}
The equation obtained after applying $q$-middle convolution $mc^q_{\lambda'}$ is written as
\begin{equation} \label{eq:21-111-111}
\frac{1}{-x}
\begin{pmatrix}
 g_1(qx) - g_1(x) \\
 g_2(qx) - g_2(x) \\
 g_3(qx) - g_3(x)
\end{pmatrix}
= \bigg[ \frac{\overline{G'}_0}{x} + \frac{\overline{G'}_1}{x-1/\gamma} \bigg]
 \begin{pmatrix}
 g_1(x) \\
 g_2(x) \\
 g_3(x)
 \end{pmatrix}.
\end{equation}
The single third-order $q$-difference equation for $g_1(x)$ is derived as
\begin{gather}
(q \gamma x-1)\bigl(q^2 \gamma x-1\bigr) g_1\bigl(q^3 x\bigr) \nonumber \\
\qquad{}- \frac{q^{ - \lambda - \lambda'}(1- q^{-\mu})\alpha^2}{(\beta-\alpha)(\alpha-\gamma)} \bigl\{q^2\bigl(q^{\mu}\beta + \alpha +q^{\lambda}\gamma\bigr)x - q^{\lambda+\lambda'} (q+1) - q^{\mu+2} \bigr\} (q \gamma x-1) g_1\bigl(q^2 x\bigr) \nonumber \\
\qquad{}- \frac{q^{- 2\lambda - 2\lambda' + 1}(1 - q^{-\mu})\alpha^2}{(\beta-\alpha)(\alpha-\gamma)}
\bigl\{ -q^2\bigl(q^{\mu}\alpha\beta + q^{\lambda+\mu}\beta\gamma + q^{\lambda}\gamma\alpha\bigr)x^2+ q\bigl(q^{\lambda+\lambda'}\alpha + q^{\lambda+\mu+1}\alpha \nonumber \\
\qquad{}+ q^{\mu+1}\beta + q^{\lambda+\lambda'+\mu}\beta + q^{2\lambda+\lambda'}\gamma + q^{\lambda+\lambda'+\mu+1}\gamma \bigr)x -q^{\lambda+\lambda'}\bigl(q^{\lambda+\lambda'} + q^{\mu+1} + q^{\mu+2}\bigr) \bigr\} g_1(qx)\! \nonumber \\
\qquad {}+ \frac{q^{-\lambda - \lambda' + 3}(q^{\mu} -1)\alpha^2}{(\beta-\alpha)\gamma}
\bigl(q^{-\lambda'}\alpha x - 1\bigr)\bigl(q^{-\lambda - \lambda'}\beta x - 1\bigr) g_1(x)
=0. \label{eq:21-111-111-g1}
\end{gather}
Each coefficient of $g_1 \bigl(q^j x\bigr)$ $(j=0,1,2,3)$ is a quadratic polynomial in $x$.
This is a generalization of the $q$-hypergeometric equation, although it is different from the generalized $q$-hypergeometric equation of order $3$ discussed in Section~\ref{sec:qmcex}.
Equation \eqref{eq:21-111-111-g1} would be novel to the best of our knowledge.
The function $g_2(x)$ satisfies a single third-order $q$-difference equation whose coefficients are quadratic polynomials in $x$, but the function $g_3(x)$ satisfies a single third-order $q$-difference equation whose coefficients are cubic polynomials in $x$.

We discuss $q$-integral representations of solutions.
On the $q$-integral transformation with respect to the $q$-convolution $c^q_{\lambda'} (B_0', B_1')= (G'_0, G'_1)$, it is not possible to apply Theorem~\ref{thm:suffcondconv}, because the matrix $I_2 -B'_0$ has the eigenvalue $1$.
We use Proposition~\ref{prop:qcintconv} instead.
It follows from equations \eqref{sol:Ytil-sta} and \eqref{eq:Yggaal} that the function
\begin{equation} \label{eq:YgxB0'B1'}
Y_g(x)
=\begin{pmatrix}
 \widetilde{y_g}_{0}(x) \\
 \widetilde{y_g}_{1}(x) \\
 \end{pmatrix}
=
\begin{pmatrix}
 \displaystyle \frac{(\gamma x; q)_{\infty}}{(\alpha x; q)_{\infty}} \int_0^{\xi \infty} \frac{K_{\lambda}(x,t)}{t} y(t) {\rm d}_q t \vspace{1mm}\\
 \displaystyle \frac{(\gamma x; q)_{\infty}}{(\alpha x; q)_{\infty}} \int_0^{\xi \infty} \frac{K_{\lambda}(x,t)}{t-1/\alpha} y(t) {\rm d}_q t
\end{pmatrix}
\end{equation}
satisfies equation \eqref{eq:B'0B'1ag} under a suitable condition.
We specialize the value $\xi $ in equation \eqref{eq:YgxB0'B1'} to apply Proposition~\ref{prop:qcintconv} for $q$-integral representation of solutions to the $q$-difference equation related to $G'_0$ and $G'_1 $ in equation \eqref{eq:G'0G'1matman}.
We set
\begin{equation} \label{eq:condPyxi}
K_{\lambda}(x, s) = K^{(1)}_{\lambda}(x, s), \qquad y(s)=s^{\mu}(\alpha s;q)_{\infty}/(\beta s;q)_{\infty} , \qquad \xi = q^{-\lambda } x
\end{equation}
in equation \eqref{eq:YgxB0'B1'}.
It follows from equation \eqref{eq:y0la} that
\begin{equation*}
\widetilde{y_g}_{0}(x) |_{\xi=q^{-\lambda}x } = (1-q) q^{-\lambda\mu} x^{\mu - \lambda} \frac{\bigl(q,q^{-\lambda+\mu+1};q\bigr)_{\infty}}{\bigl(q^{-\lambda+1},q^{\mu};q\bigr)_{\infty}} \frac{(\gamma x; q)_{\infty}}{(\alpha x; q)_{\infty}}
{}_2\phi_1 \left( \begin{matrix} \alpha/\beta, q^{\mu}\\
q^{-\lambda+\mu+1} \end{matrix} ;q,q^{-\lambda}\beta x \right) ,
\end{equation*}
under the condition $\mu >0 $.
It follows from the $q$-difference equation that $\widetilde{y_g}_{1}(x)$ is expressed by a~combination of $\widetilde{y_g}_{0}(q x) $ and $\widetilde{y_g}_{0}(x) $.
Hence, if $\mu - \lambda >0 $, then the assumption of Proposition~\ref{prop:qcintconv}\,(i) to the function $Y_g(x) $ in equation \eqref{eq:YgxB0'B1'} holds because of the term $x^{\mu - \lambda} $.
We are going to apply Theorem~\ref{thm:qcint} in the situation that equations \eqref{eq:YqxBYx} and~\eqref{eq:qcintadd} correspond to equation~\eqref{eq:B'0B'1ag} and
\begin{equation} \label{eq:Ygan31}
\widetilde{Y_g}(x)
= \begin{pmatrix}
 \displaystyle \int_0^{\xi' \infty} \frac{K_{\lambda'}(x,s)}{s} Y_g (s) {\rm d}_q s \vspace{1mm}\\
 \displaystyle \int_0^{\xi' \infty} \frac{K_{\lambda'}(x,s)}{s - 1/\gamma} Y_g (s) {\rm d}_q s
 \end{pmatrix}.
\end{equation}
Since the eigenvalues of $I_2- B'_0 -B'_1$ are $q^{\mu-\lambda} \beta / \gamma $ and $q^{-\lambda} \alpha /\gamma $, the assumption of Proposition~\ref{prop:qcintconvYx}\,(ii) is written as $\big|q^{\lambda ' + \lambda -\mu } \gamma / \beta \big| <1 $ and $ \big|q^{\lambda ' +\lambda } \gamma / \alpha \big| <1 $, and Proposition~\ref{prop:qcintconv}\,(ii) holds under this condition.
By Proposition~\ref{prop:qcintconv}, the assumption of Theorem~\ref{thm:qcint} holds if $\mu >0$, $\mu - \lambda >0 $, $\big|q^{\lambda ' + \lambda -\mu } \gamma / \beta \big| <1 $, $ \big|q^{\lambda ' +\lambda } \gamma / \alpha \big| <1 $ and the condition in \eqref{eq:condPyxi} is satisfied.
Then, the function~\smash{$\widetilde{Y_g}(x) $} in equation \eqref{eq:Ygan31} is written as
\begin{align}
 \widetilde{Y_g}(x)
&=\begin{pmatrix}
 \widetilde{y_g}_{00}(x) \\
 \widetilde{y_g}_{01}(x) \\
 \widetilde{y_g}_{10}(x) \\
 \widetilde{y_g}_{11}(x)
 \end{pmatrix}
\nonumber \\
&= \begin{pmatrix}
 \displaystyle \int_0^{\xi' \infty} K_{\lambda'}(x,s) s^{-\lambda -1} \frac{(\gamma s; q)_{\infty}}{(\alpha s; q)_{\infty}} \int_0^{q^{-\lambda } s \infty} t^{\mu -1} \frac{\bigl(q^{\lambda+1}t/s , \alpha t ;q\bigr)_{\infty}}{(qt/s , \beta t ;q)_{\infty}} {\rm d}_q t {\rm d}_q s \vspace{1mm}\\
 \displaystyle - \alpha \int_0^{\xi' \infty} K_{\lambda'}(x,s) s^{-\lambda -1} \frac{(\gamma s; q)_{\infty}}{(\alpha s; q)_{\infty}} \int_0^{q^{-\lambda } s \infty} t^{\mu} \frac{\bigl(q^{\lambda+1}t/s , q \alpha t ;q\bigr)_{\infty}}{(qt/s , \beta t ;q)_{\infty}} {\rm d}_q t {\rm d}_q s \vspace{1mm}\\
 \displaystyle -\gamma \int_0^{\xi' \infty} K_{\lambda'}(x,s) s^{-\lambda } \frac{(q \gamma s; q)_{\infty}}{(\alpha s; q)_{\infty}} \int_0^{q^{-\lambda } s \infty} t^{\mu -1} \frac{\bigl(q^{\lambda+1}t/s , \alpha t ;q\bigr)_{\infty}}{(qt/s , \beta t ;q)_{\infty}} {\rm d}_q t {\rm d}_q s \vspace{1mm}\\
 \displaystyle \alpha \gamma \int_0^{\xi' \infty} K_{\lambda'}(x,s) s^{-\lambda } \frac{(q \gamma s; q)_{\infty}}{(\alpha s; q)_{\infty}} \int_0^{q^{-\lambda } s \infty} t^{\mu} \frac{\bigl(q^{\lambda+1}t/s , q \alpha t ;q\bigr)_{\infty}}{(qt/s , \beta t ;q)_{\infty}} {\rm d}_q t {\rm d}_q s
 \end{pmatrix}.  \label{eq:Ygan32}
\end{align}
Consequently, we obtain the following proposition.
\begin{Proposition}
Let $K_{\lambda }(x,s)$ be a function which satisfies equation \eqref{eq:tPlambda}.
If $\mu >0$, $\mu - \lambda >0 $, $\big|q^{\lambda ' + \lambda -\mu } \gamma / \beta \big| <1 $ and $ \big|q^{\lambda ' +\lambda } \gamma / \alpha \big| <1 $ and the integrands in equation \eqref{eq:Ygan32} are well defined about \smash{$(s,t) = \bigl(q^{n'} \xi ', q^{n +n' -\lambda } \xi ' \bigr)$} for any $n', n \in \mathbb{Z} $, then the function \smash{$\widetilde{Y_g}(x)$} in equation \eqref{eq:Ygan32} converges and satisfies
\begin{equation} \label{eq:G'0G'1man}
\frac{\widetilde{Y_g}(qx) - \widetilde{Y_g}(x)}{-x} = \bigg[ \frac{G'_0}{x} + \frac{G'_1}{x-1/\alpha} \bigg] \widetilde{Y_g}(x),
\end{equation}
where $G'_0$ and $G'_1$ are given in equation \eqref{eq:G'0G'1matman}.
\end{Proposition}

If $\widetilde{Y_g}(x) $ is a solution to equation \eqref{eq:G'0G'1man} and set
\begin{equation*}
\begin{pmatrix}
 g_1(x) \\
 g_2(x) \\
 g_3(x) \\
 g_4(x)
\end{pmatrix}
= P^{-1} \widetilde{Y_g}(x) = \begin{pmatrix}
 0 & 0 & 0 & 1 \\
 \displaystyle \frac{q^{\mu} - q^{\lambda}}{q^{\mu}(\beta /\alpha - 1)} & 1 & 0 & 0 \\
 0 & 0 & 1 & 0 \\
 \displaystyle \frac{-1}{q^{\mu}(\beta /\alpha - 1)} & 0 & 0 & 0
\end{pmatrix} \begin{pmatrix}
 \widetilde{y_g}_{00}(x) \\
 \widetilde{y_g}_{01}(x) \\
 \widetilde{y_g}_{10}(x) \\
 \widetilde{y_g}_{11}(x)
 \end{pmatrix},
\end{equation*}
then the function ${}^t(g_1 (x), g_2 (x), g_3 (x) ) $ satisfies equation \eqref{eq:21-111-111}.
In the case $K_{\lambda'}(x,s)\! =\! K^{(1)}_{\lambda '}(x, s) $, we have
\begin{align*}
g_1(x)={}& \alpha \gamma x^{-\lambda '} \int_0^{\xi' \infty} s^{-\lambda } \frac{\bigl(q^{\lambda '+1}s/x, q \gamma s ;q\bigr)_{\infty}}{(qs/x, \alpha s ;q)_{\infty} } \int_0^{q^{-\lambda } s \infty} t^{\mu} \frac{\bigl(q^{\lambda+1}t/s , q \alpha t ;q\bigr)_{\infty}}{(qt/s , \beta t ;q)_{\infty}} {\rm d}_q t {\rm d}_q s \\
={}& (1 - q)^2 \alpha \gamma \xi'^{\mu -\lambda + 2} q^{-\lambda (\mu +1 ) } x^{-\lambda'} \sum_{n = -\infty}^{\infty} \bigl(q^{\mu -\lambda + 2}\bigr)^n \frac{\bigl(q^{\lambda' + n + 1}\xi'/x, q^{n + 1}\xi' \gamma; q\bigr)_{\infty}}{\bigl(q^{n + 1}\xi'/x, q^n\xi' \alpha; q\bigr)_{\infty}} \\
& \times \sum_{m = 0}^{\infty} \bigl(q^{\mu + 1}\bigr)^m \frac{\bigl(q^{m + 1} , q^{-\lambda + n + m + 1} \xi' \alpha; q\bigr)_{\infty}}{\bigl(q^{-\lambda + m + 1} , q^{-\lambda +n+m} \xi' \beta; q\bigr)_{\infty}} ,
\end{align*}
and this is a solution to equation \eqref{eq:21-111-111-g1}.

\subsection[q-convolution and q-Jordan Pochhammer equation]{$\boldsymbol{q}$-convolution and $\boldsymbol{q}$-Jordan Pochhammer equation} \label{sec:JP}

We review the $q$-Jordan Pochhammer equation in terms of the $q$-convolution, which is a slight modification of \cite[Section~3.2]{AT}.
Let $N$ be a positive integer, and set
\begin{equation*}
y(x)=x^{\mu} \prod _{j=1}^N \frac{(\alpha _j x;q)_{\infty}}{(\beta _j x;q)_{\infty}}.
\end{equation*}
It satisfies the single $q$-difference equation
\begin{align}
& \frac{y(qx)-y(x)}{-x} = \Biggl[ \frac{B_0}{x} + \sum _{k=1}^N \frac{B_k}{x-1/\alpha _k} \Biggr] y(x) ,\nonumber \\
& B_{0}=1- q^{\mu}, \qquad B_k =q^{\mu} \frac{\alpha _k -\beta _k}{\alpha _k} \prod _{j=1, j\neq k}^N \frac{\alpha _k -\beta _j}{\alpha _k - \alpha _j} , \qquad k=1, \dots ,N. \label{eq:B0B1BN}
\end{align}
We apply the convolution $c^q_{\lambda }$.
The matrices $ c^q_{\lambda } (B_0,\dots, B_N) = (G_0,\dots,G_N)$ are written as
\begin{gather}
 G_i = \left(
 \begin{matrix}
 {} & {} & {} & O & {} & {} \\
 q^{-\lambda } B_0& q^{-\lambda } B_1 & \cdots & q^{-\lambda } B_i + 1 - q^{-\lambda } & \cdots & q^{-\lambda } B_N \\
 {} & {} & {} & O & {} & {}
 \end{matrix}
 \right)
 {\scriptstyle (i+1)},\! \qquad 0\leq i \leq N,\!\!\! \label{eq:GiJP}
\end{gather}
By applying Proposition~\ref{prop:qcintconv} and Theorem~\ref{thm:qcint}, we obtain the following theorem.
\begin{Theorem}[{cf.\ \cite[Theorem 3.5]{AT}}] \label{thm:qcintqJH}
Let $\xi \in \mathbb{C} \setminus \{ 0 \} $ and let $K_{\lambda}(x, s)$ be a function which satisfies equation \eqref{eq:tPlambda}.
If $\mu >0 $ and $\big| q ^{\lambda -\mu } \alpha _1 \cdots \alpha _N / (\beta _1 \cdots \beta _N)\big| < 1$, then the function $\widetilde{Y}(x) $ defined by
\begin{gather*}
 \widetilde{Y}(x) = \begin{pmatrix} \widetilde{y}_{0}(x) \\ \widetilde{y}_{1}(x) \\ \vdots \\ \widetilde{y}_{N}(x) \end{pmatrix} ,\qquad
 \widetilde{y}_{0}(x) = \int^{\xi \infty}_{0}\frac{K_{\lambda}(x, s)}{s} s^{\mu} \prod _{j=1}^N \frac{(\alpha _j s;q)_{\infty}}{(\beta _j s;q)_{\infty}} {\rm d}_{q}s,\\
 \widetilde{y}_{k}(x) = \int^{\xi \infty}_{0}\frac{K_{\lambda}(x, s)}{s-1/\alpha _k}s^{\mu}\prod _{j=1}^N \frac{(\alpha _j s;q)_{\infty}}{(\beta _j s;q)_{\infty}} {\rm d}_{q}s ,\qquad
 k=1,\dots ,N,
\end{gather*}
is convergent and it satisfies the equation
\begin{equation}
 \frac{\widetilde{Y}(qx)-\widetilde{Y}(x)}{-x} = \Biggl[ \frac{G_0}{x} + \sum _{k=1}^N \frac{G_k}{x-1/\alpha _k} \Biggr] \widetilde{Y}(x) , \label{eq:G0G1GN}
\end{equation}
where $G_0, G_1, \dots , G_N $ are determined in equation \eqref{eq:GiJP}.
\end{Theorem}
Note that the tuple $(G_0, G_1, \dots , G_N )$ is irreducible for generic parameters.
If the parameters are special, then the tuple $(G_0, G_1, \dots , G_N )$ is reducible and the system of the $q$-difference equations~\eqref{eq:G0G1GN} may have a subsystem of equations.
If $\mu =0$, then the space $\mathcal{K} $ in equation \eqref{eq:qKL}, which was used to define the $q$-middle convolution, is non-zero, i.e., $\dim \mathcal{K} \geq 1$.
If~${q ^{\lambda -\mu } \alpha _1 \cdots \alpha _N / (\beta _1 \cdots \beta _N) =1}$, then the space $\mathcal{L} $ in equation \eqref{eq:qKL} satisfies $\dim \mathcal{L} \geq 1$.

\subsection[q-middle convolution and variants of q-hypergeometric equation]{$\boldsymbol{q}$-middle convolution and variants of $\boldsymbol{q}$-hypergeometric equation} \label{sec:qconvex}

We investigate the case ($N=2 $ and $\mu =0 $) and the case ($N=3 $, $\mu =0 $ and $q^{\lambda}=\beta_1\beta_2\beta_3/\allowbreak(\alpha_1\alpha_2\alpha_3) $).
In the case $N=2 $ and $\mu =0 $, the single $q$-difference equation in \eqref{eq:B0B1BN} is written as
\begin{align}
&\frac{y(qx) - y(x)}{-x} = \bigg[ \frac{B_0}{x} + \frac{B_1}{x-1/\alpha_1} + \frac{B_2}{x-1/\alpha_2} \bigg] y(x), \nonumber\\
& B_0= 0, \qquad
B_1= \frac{(\alpha_1-\beta_1)(\alpha_1-\beta_2)}{\alpha_1(\alpha_1-\alpha_2)},\qquad
B_2= \frac{(\alpha_2-\beta_1)(\alpha_2-\beta_2)}{\alpha_2(\alpha_2-\alpha_1)}. \label{eq:qeqvar2}
\end{align}
The function $y(x)=(\alpha_1 x, \alpha_2 x;q)_{\infty}/(\beta_1 x, \beta_2 x;q)_{\infty}$ satisfies this equation.
By applying the $q$-convolution $c^q_{\lambda}$ to $(B_0, B_1, B_2)$, we obtain the $3\times 3$ matrices $G_0$, $G_1$, $G_2$ in equation \eqref{eq:GiJP}.
The $q$-middle convolution $mc^q_{\lambda}$ is formulated on the space $\mathbb{C} ^3/ ( \mathcal{K} + \mathcal{L} )$, and we have $\mathcal{L} =\{ 0\} $ and a~basis of $\mathcal{K} $ is ${}^t (1,0,0) $ for generic $\alpha_1, \alpha_2, \beta_1, \beta_2 \in \mathbb{C} \setminus \{0\}$.
Write $mc^q_{\lambda} (B_0, B_1, B_2) = \bigl( \overline{G}_0, \overline{G}_1, \overline{G}_2\bigr)$.
The matrix representations of $\overline{G}_0, \overline{G}_1, \overline{G}_2 \in \mbox{End} \bigl(\mathbb{C} ^3/ ( \mathcal{K} + \mathcal{L} )\bigr) $ are realized as the $2\times 2$ lower-right submatrices (see \cite{AT}).
Namely,
\begin{gather}
\overline{G}_0=O, \qquad
 \overline{G}_1= \begin{pmatrix}
 q^{-\lambda}B_1 + 1 -q^{-\lambda} & q^{-\lambda}B_2 \\
 0 & 0
 \end{pmatrix},\nonumber \\
 \overline{G}_2= \begin{pmatrix}
 0 & 0 \\
 q^{-\lambda}B_1 & q^{-\lambda}B_2 + 1 -q^{-\lambda}
 \end{pmatrix}, \label{eq:oG1oG2}
\end{gather}
and the $q$-difference equation corresponding to $mc^q_{\lambda} (B_0, B_1, B_2)$ is written as
\begin{equation} \label{eq:v-qhg-deg2}
 \frac{\overline{Y}(qx) - \overline{Y}(x)}{-x}
 = \Biggl[ \frac{\overline{G}_1}{x-1/\alpha _1} +\frac{\overline{G}_2}{x-1/\alpha _2} \Biggr] \overline{Y}(x).
%
\end{equation}
The single second-order $q$-difference equations derived from equation \eqref{eq:v-qhg-deg2} are the variant of $q$-hypergeometric equation of degree~2~\cite{HMST}.
Formal $q$-integral representations of solutions are written as
\begin{equation} \label{eq:qintvarqhg2}
\overline{Y}(x)
= \begin{pmatrix}
 \displaystyle \int_0^{\xi \infty} \frac{K_{\lambda}(x,t)}{t-1/\alpha_1} y(t) {\rm d}_q t \vspace{1mm}\\
 \displaystyle \int_0^{\xi \infty} \frac{K_{\lambda}(x,t)}{t-1/\alpha_2} y(t) {\rm d}_q t
 \end{pmatrix},
\end{equation}
where $y(t)$ is a solution to equation \eqref{eq:qeqvar2} and $K_{\lambda }(x,t)$ is a function which satisfies equation~\eqref{eq:tPlambda}.
However, they are not generally actual solutions, because the assumption of Theorem~\ref{thm:suffcondconv} does not hold.
As we discussed in \cite{AT}, the $q$-integral in equation \eqref{eq:qintvarqhg2} converges and satisfies an non-homogeneous version of equation \eqref{eq:v-qhg-deg2}.
If the functions $K_{\lambda}(x,t)$ and~$y(t)$ and the value~$\xi $ are chosen appropriately (e.g., $\xi =1/\beta_1, 1/\beta _2 $ or $x$), the $q$-integral in equation \eqref{eq:qintvarqhg2} satisfies equation \eqref{eq:v-qhg-deg2}.
For details see \cite[Section~4.1.1]{AT}.

We discuss the case $N=3 $ and $\mu =0 $.
In this case, the single $q$-difference equation in \eqref{eq:B0B1BN} is written as
\begin{gather}
\frac{y(qx) - y(x)}{-x} = \bigg[ \frac{B_0}{x} + \frac{B_1}{x-1/\alpha_1} + \frac{B_2}{x-1/\alpha_2} + \frac{B_3}{x-1/\alpha_3} \bigg] y(x), \nonumber\\
 B_0= 0, \qquad B_1= \frac{(\alpha_1-\beta_1)(\alpha_1-\beta_2)(\alpha_1-\beta_3)}{\alpha_1(\alpha_1-\alpha_2)(\alpha_1-\alpha_3)}, \nonumber \\
 B_2= \frac{(\alpha_2-\beta_1)(\alpha_2-\beta_2)(\alpha_2-\beta_3)}{\alpha_2(\alpha_2-\alpha_3)(\alpha_2-\alpha_1)}, \qquad B_3= \frac{(\alpha_3-\beta_1)(\alpha_3-\beta_2)(\alpha_3-\beta_3)}{\alpha_3(\alpha_3-\alpha_2)(\alpha_3-\alpha_1)}. \label{eq:B1B2B3}
\end{gather}
By applying the $q$-convolution $c^q_{\lambda}$ to $(B_0, B_1, B_2 ,B_3)$, we obtain the $4\times 4$ matrices $G_0$, $G_1$, $G_2$, $G_3$ in equation \eqref{eq:GiJP}.
We impose the condition
\begin{equation}
\mu=0, \qquad q^{\lambda}=\beta_1\beta_2\beta_3/(\alpha_1\alpha_2\alpha_3).
\label{eq:qlb1b2b3a1a2a3}
\end{equation}
The $q$-middle convolution $mc^q_{\lambda}$ is formulated on the space $\mathbb{C} ^4/ ( \mathcal{K} + \mathcal{L} )$, and it is shown that $\dim \mathcal{K} =1= \dim \mathcal{L} $ and $\dim ( \mathcal{K} + \mathcal{L} ) =2$ for generic $\alpha_1, \alpha_2, \alpha_3, \beta_1, \beta_2, \beta_3 \in \mathbb{C} \setminus \{0\}$ which satisfy equation \eqref{eq:qlb1b2b3a1a2a3} (see \cite{AT}).
By choosing an appropriate basis of $\mathbb{C} ^4/ ( \mathcal{K} + \mathcal{L} ) \simeq \mathbb{C} ^2$, we obtain the tuple of $2\times 2$ matrices $\bigl(\overline{G}_0, \overline{G}_1, \overline{G}_2 , \overline{G}_3 \bigr)$ and the $q$-difference equation described as
\begin{align}
&\frac{\overline{Y}(qx) - \overline{Y}(x)}{-x} = \bigg[ \frac{\overline{G}_1}{x-1/\alpha_1} + \frac{\overline{G}_2}{x-1/\alpha_2} + \frac{\overline{G}_3}{x-1/\alpha_3} \bigg] \overline{Y}(x), \nonumber \\
& \overline{G}_0=O, \qquad \overline{G}_1= \begin{pmatrix}
 q^{-\lambda}B_1 + 1 -q^{-\lambda} & q^{-\lambda}B_2 \\
 0 & 0
 \end{pmatrix}, \nonumber \\
& \overline{G}_2= \begin{pmatrix}
 0 & 0 \\
 q^{-\lambda}B_1 & q^{-\lambda}B_2 + 1 -q^{-\lambda}
 \end{pmatrix}, \qquad \overline{G}_3= \begin{pmatrix}
 -q^{-\lambda}B_1 & -q^{-\lambda}B_2 \\
 -q^{-\lambda}B_1 & -q^{-\lambda}B_2
 \end{pmatrix}.\label{eq:v-qhg-deg3}
\end{align}
The single second-order $q$-difference equations derived from equation \eqref{eq:v-qhg-deg3} correspond to the variant of $q$-hypergeometric equation of degree~$3$~\cite{HMST}.
What we have mentioned above is a~recalculation of what is written in \cite{AT} by using the reformulated $q$-convolution.
Note that $q$-integral representations of solutions to the variant of $q$-hypergeometric equation of degree $3$ were discussed in \cite{AT} from the framework of the $q$-middle convolution, and Fujii and Nobukawa obtained significant results on solutions to the variant of $q$-hypergeometric equation of degree $3$ in \cite{FN}.

In the following subsections, we apply gauge-transformations and the $q$-middle convolution to the tuple $\bigl(O, \overline{G}_1, \overline{G}_2 \bigr) $ in equation \eqref{eq:oG1oG2} or the tuple $\bigl(O, \overline{G}_1, \overline{G}_2, \overline{G}_3 \bigr) $ in equation \eqref{eq:v-qhg-deg3}.

\subsection[Third-order extension of variants of q-hypergeometric equation]{Third-order extension of variants of $\boldsymbol{q}$-hypergeometric equation} \label{subsec:3exte}

\subsubsection[Extension of variant of q-hypergeometric equation of degree 2]{Extension of variant of $\boldsymbol{q}$-hypergeometric equation of degree 2} \label{subsec:3extedeg2}

In this subsection, we apply the addition with respect to gauge-transformation
\begin{equation}
 Y_g(x) = \frac{(\gamma_1 x;q)_{\infty}}{(\alpha_1 x;q)_{\infty}}\overline{Y}(x)
\label{eq:Ygdeg2-1}
\end{equation}
to the tuple $\bigl(O, \overline{G}_1, \overline{G}_2 \bigr) $ in equation \eqref{eq:oG1oG2}.
Then, we obtain the tuple $(B'_0, B'_1, B'_2 ) $ and the $q$-difference equation
\begin{equation} \label{eq:B'1B'2gm1}
\frac{Y_g(qx) - Y_g(x)}{-x} = \bigg[ \frac{B'_1}{x-1/\gamma_1} + \frac{B'_2}{x-1/\alpha_2} \bigg] Y_g(x) ,
\end{equation}
where
\begin{align*}
& B'_0 = O, B'_1
= \begin{pmatrix}
 b^1_{11} & b^1_{12} \\
 b^1_{21} & b^1_{22}
 \end{pmatrix},\qquad
 B'_2
= \begin{pmatrix}
 0 & 0 \\
 b^2_{21} & b^2_{22}
\end{pmatrix}
, \qquad b^1_{11} = 1 + q^{- \lambda} (B_1 - 1) \frac{\alpha_1 }{\gamma_1 },
\\
& b^1_{12} = q^{-\lambda} B_2 \frac{\alpha_1 }{\gamma_1}, \qquad b^1_{21} = q^{-\lambda} B_1 \frac{\alpha_2(\alpha_1 - \gamma_1)}{\gamma_1(\alpha_2 - \gamma_1)} ,\qquad b^1_{22} = \left\{1 + q^{-\lambda} (B_2 - 1) \frac{\alpha_2 }{\gamma_1} \right\} \frac{\alpha_1 - \gamma_1}{\alpha_2 - \gamma_1} , \\
& b^2_{21} = q^{- \lambda} B_1 \frac{\alpha_1-\alpha_2}{\gamma_1-\alpha_2},\qquad b^2_{22} = \bigl\{ 1 + q^{-\lambda} (B_2- 1) \bigr\} \frac{\alpha_1-\alpha_2}{\gamma_1-\alpha_2} , \\
& B_1=\frac{(\alpha_1-\beta_1)(\alpha_1-\beta_2)}{\alpha_1(\alpha_1-\alpha_2)},\qquad B_2=\frac{(\alpha_2-\beta_1)(\alpha_2-\beta_2)}{\alpha_2(\alpha_2-\alpha_1)}.
\end{align*}
We apply the $q$-convolution.
Then, we have $c^q_{\lambda'} (O, B'_1, B'_2) = (G'_0, G'_1, G'_2)$, where
\begin{gather}
G'_0
= \begin{pmatrix}
 \bigl(1 - q^{-\lambda'}\bigr)I_2 & q^{-\lambda '} B'_1 & q^{-\lambda '} B'_2\\
 O & O & O \\
 O & O & O
\end{pmatrix},\nonumber \\
G'_1
= \begin{pmatrix}
 O & O & O \\
 O & q^{-\lambda'} B'_1 + \bigl(1-q^{-\lambda'}\bigr) I_2 & q^{-\lambda'} B'_2 \\
 O & O & O
\end{pmatrix},\nonumber \\
G'_2
= \begin{pmatrix}
 O & O & O \\
 O & O & O \\
 O & q^{-\lambda'} B'_1 & q^{-\lambda'} B'_2 + \bigl(1-q^{-\lambda'}\bigr) I_2
\end{pmatrix}.\label{eq:G'012}
\end{gather}
The vectors ${}^t(1, 0)$ and ${}^t(0, 1)$ are a basis for the space $\ker (B'_0)$, and the vector ${}^t\bigl(q^{\lambda} + B_2 - 1, -B_1\bigr)$ belongs to the space $\ker (B'_2)$.
We have $\dim \mathcal{K}=3$, $\dim \mathcal{L}=0$ and $\dim \bigl(\mathbb{C}^6/(\mathcal{K} + \mathcal{L}) \bigr) =3 $ for generic parameters.
Write $mc^q_{\lambda'} (O, B'_1, B'_2) = \bigl(\overline{G'}_0, \overline{G'}_1, \overline{G'}_2\bigr)$.
By the simultaneous transformation for $(G'_0, G'_1, G'_2)$ by the matrix
\begin{equation*}
P = \begin{pmatrix}
 0 & 0 & 0 & 1 & 0 & 0 \\
 0 & 0 & 0 & 0 & 1 & 0 \\
 1 & 0 & 0 & 0 & 0 & 0 \\
 0 & 1 & 0 & 0 & 0 & 0 \\
 0 & 0 & 1 & 0 & 0 & q^{\lambda} + B_2 - 1 \\
 0 & 0 & 0 & 0 & 0 & -B_1
\end{pmatrix},
\end{equation*}
we obtain matrix representations $\overline{G'}_0$, $\overline{G'}_1$ and $\overline{G'}_2$ on the space~${\mathbb{C}^6/(\mathcal{K} + \mathcal{L}) }$ by the upper-left~${3\times 3}$ submatrices of $P^{-1}G'_0 P$, $P^{-1}G'_1 P$ and $P^{-1}G'_2 P$.
They are expressed as
\begin{gather*}
\overline{G'}_0 = O, \overline{G'}_1
= \begin{pmatrix}
 1 + q^{-\lambda'}(b^1_{11} - 1) & q^{-\lambda'}b^1_{12} & 0 \\
 q^{-\lambda'}b^1_{21} & 1 + q^{-\lambda'}\bigl(b^1_{22} - 1\bigr) & q^{-\lambda'}b^2_{21} \\
 0 & 0 & 0
 \end{pmatrix}, \\
 \overline{G'}_2
= \begin{pmatrix}
 0 & 0 & 0 \\
 0 & 0 & 0 \\
 l_1 & l_2 & l_3
 \end{pmatrix},
\qquad
l_1 = q^{-\lambda'} \bigl\{ b^1_{11} + \bigl(B_2 + q^{\lambda}-1\bigr) b^1_{21} / B_1 \bigr\} , \\
l_2 = q^{-\lambda'} \bigl\{ b^1_{12} + \bigl(B_2 + q^{\lambda}-1\bigr) b^1_{22}/B_1 \bigr\},\qquad
l_3 = 1 - q^{-\lambda'} + q^{-\lambda'} \bigl\{ \bigl(B_2 + q^{\lambda}-1\bigr) b^2_{21} /B_1 \bigr\}.
\end{gather*}
The equation obtained after applying $q$-middle convolution $mc^q_{\lambda'}$ is written as
\begin{equation} \label{eq:3-111-21-111}
\frac{1}{-x}
\begin{pmatrix}
 g_1(qx) - g_1(x) \\
 g_2(qx) - g_2(x) \\
 g_3(qx) - g_3(x)
\end{pmatrix}
= \bigg[ \frac{\overline{G'}_1}{x-1/\gamma_1} + \frac{\overline{G'}_2}{x-1/\alpha_2} \bigg]
 \begin{pmatrix}
 g_1(x) \\
 g_2(x) \\
 g_3(x)
 \end{pmatrix}.
\end{equation}
We derive single third-order $q$-difference equations from equation \eqref{eq:3-111-21-111}.
The functions $g_1(x)$ and~$g_3(x)$ satisfy single third-order $q$-difference equations whose coefficients are cubic polynomials in $x$, and $g_2(x)$ satisfies one whose coefficients are quartic polynomials in $x$.

The single third-order $q$-difference equation satisfied by $g_1(x)$ and $q$-integral representations of solutions are described in Appendix \ref{appsec:B}.

Next, we apply the addition with respect to the gauge-transformation
\[
 Y_g(x) = \smash{\frac{(\gamma_1 x, \gamma_2 x;q)_{\infty}}{(\alpha_1 x, \alpha_2 x;q)_{\infty}}} \overline{Y}(x)
\]
to the tuple $\bigl(O, \overline{G}_1, \overline{G}_2 \bigr) $ in equation \eqref{eq:oG1oG2}.
Then, we obtain the tuple $(B'_0, B'_1, B'_2 ) $ and the $q$-difference equation
\begin{equation*} 
\frac{Y_g(qx) - Y_g(x)}{-x} = \bigg[ \frac{B'_1}{x-1/\gamma_1} + \frac{B'_2}{x-1/\gamma_2} \bigg] Y_g(x) ,
\end{equation*}
where
\begin{gather*}
 B'_0 = O, B'_1 = \begin{pmatrix}
 b^1_{11} & b^1_{12} \\
 b^1_{21} & b^1_{22}
 \end{pmatrix},\qquad
 B'_2
= \begin{pmatrix}
 b^2_{11} & b^2_{12} \\
 b^2_{21} & b^2_{22}
\end{pmatrix},\\
 b^1_{11} = \left\{ 1 + q^{-\lambda} (B_1 - 1) \frac{\alpha_1 }{\gamma_1} \right\} \frac{\alpha_2 - \gamma_1 }{\gamma_2 - \gamma_1} , \\
 b^1_{12} = q^{-\lambda} B_2 \frac{\alpha_1(\alpha_2 - \gamma_1)}{\gamma_1(\gamma_2 - \gamma_1)}, \qquad b^1_{21} = q^{-\lambda}B_1 \frac{\alpha_2(\alpha_1 - \gamma_1)}{\gamma_1(\gamma_2 - \gamma_1)} , \\
 b^1_{22} = \left\{ 1 + q^{-\lambda} (B_2 - 1) \frac{\alpha_2 }{\gamma_1} \right\} \frac{\alpha_1 - \gamma_1 }{\gamma_2 - \gamma_1 } , \\
 b^2_{11} = \left\{ 1 + q^{-\lambda} (B_1 - 1) \frac{\alpha_1 }{\gamma_2} \right\} \frac{\alpha_2 - \gamma_2 }{\gamma_1-\gamma_2},\qquad
 b^2_{12} = q^{-\lambda} B_2 \frac{\alpha_1(\alpha_2 - \gamma_2)}{\gamma_2(\gamma_1 - \gamma_2)} , \\
 b^2_{21} = q^{-\lambda} B_1 \frac{\alpha_2(\alpha_1 - \gamma_2)}{\gamma_2(\gamma_1 - \gamma_2)} , \qquad b^2_{22}= \Big\{ 1 + q^{-\lambda} (B_2 - 1) \frac{\alpha_2 }{\gamma_2} \Big\} \frac{\alpha_1 - \gamma_2 }{\gamma_1-\gamma_2} , \\
 B_1= \frac{(\alpha_1-\beta_1)(\alpha_1-\beta_2)}{\alpha_1(\alpha_1-\alpha_2)},\qquad B_2=\frac{(\alpha_2-\beta_1)(\alpha_2-\beta_2)}{\alpha_2(\alpha_2-\alpha_1)}.
\end{gather*}
Then, we have $c^q_{\lambda'} (O, B'_1, B'_2) = (G'_0, G'_1, G'_2)$, where the tuple $(G'_0, G'_1, G'_2)$ is the $6 \times 6$ matrices given as equation \eqref{eq:G'012}.
The vectors ${}^t(1, 0)$ and ${}^t(0, 1)$ are a basis for the space $\ker (B'_0)$.
Hence, we have $\dim \mathcal{K} \geq 2 $.

We look for a description of the condition $\dim \mathcal{L}\geq 1$.
If $q^{\lambda+\lambda'}=\alpha_1\alpha_2/(\gamma_1\gamma_2)$ holds, then the vector ${}^t(B_2,-B_1,B_2,-B_1,B_2,-B_1)$ belongs to the space $\mathcal{L}$.
We have $\dim \mathcal{K}=2$ and ${\dim \mathcal{L}=1}$ for the generic parameters which satisfy $q^{\lambda+\lambda'}=\alpha_1\alpha_2/(\gamma_1\gamma_2)$.
Write
\[
mc^q_{\lambda'} (O, B'_1, B'_2) = \bigl(\overline{G'}_0, \overline{G'}_1, \overline{G'}_2\bigr).
\]
By the simultaneous transformation for $(G'_0, G'_1, G'_2)$ by the matrix
\begin{equation*}
P = \begin{pmatrix}
 0 & 0 & 0 & 1 & 0 & B_2 \\
 0 & 0 & 0 & 0 & 1 & -B_1 \\
 1 & 0 & 0 & 0 & 0 & B_2 \\
 0 & 1 & 0 & 0 & 0 & -B_1 \\
 0 & 0 & 1 & 0 & 0 & B_2 \\
 0 & 0 & 0 & 0 & 0 & -B_1
\end{pmatrix},
\end{equation*}
we obtain matrix representations $\overline{G'}_0$, $\overline{G'}_1$ and $\overline{G'}_2$ on the space $\mathbb{C}^6/(\mathcal{K} + \mathcal{L}) $ by the upper-left~${3\times 3}$ submatrices of $P^{-1}G'_0 P$, $P^{-1}G'_1 P$ and $P^{-1}G'_2 P$.
They are expressed as
\begin{align*}
&\overline{G'}_0 = O , \qquad \overline{G'}_1
= \begin{pmatrix}
 1 - q^{-\lambda'} + q^{-\lambda'} b^1_{11} & q^{-\lambda'}b^1_{12} & q^{-\lambda'}b^2_{11} \\
 q^{-\lambda'}b^1_{21} & 1 - q^{-\lambda'} + q^{-\lambda'} b^1_{22} & q^{-\lambda'}b^2_{21} \\
 0 & 0 & 0
 \end{pmatrix}, \\
&\overline{G'}_2
= q^{-\lambda'} \begin{pmatrix}
 \displaystyle b^1_{21} B_2 /B_1
 & \displaystyle b^1_{22} B_2 /B_1
 & \displaystyle b^2_{21} B_2 /B_1 \\
 - b^1_{21}
 & - b^1_{22}
 & - b^2_{21} \\
 \displaystyle b^1_{11} + b^1_{21} B_2 /B_1
 & \displaystyle b^1_{12} + b^1_{22} B_2 /B_1
 & q^{\lambda'} -1 + \displaystyle b^2_{11}+b^2_{21} B_2 /B_1
 \end{pmatrix}.
\end{align*}
The equation obtained after applying $q$-middle convolution $mc^q_{\lambda'}$ is written as
\begin{equation} \label{eq:3-111-111-21}
\frac{1}{-x}
\begin{pmatrix}
 g_1(qx) - g_1(x) \\
 g_2(qx) - g_2(x) \\
 g_3(qx) - g_3(x)
\end{pmatrix}
= \bigg[ \frac{\overline{G'}_1}{x-1/\gamma_1} + \frac{\overline{G'}_2}{x-1/\gamma_2} \bigg]
 \begin{pmatrix}
 g_1(x) \\
 g_2(x) \\
 g_3(x)
 \end{pmatrix}.
\end{equation}
We derive single third-order $q$-difference equations from equation \eqref{eq:3-111-111-21}.
The functions $g_1(x)$, $g_2(x)$, and $g_3(x)$ satisfy single third-order $q$-difference equations whose coefficients are polynomials in $x$ of degree $6$, $5$, and $4$, respectively.

The single third-order $q$-difference equation satisfied by $g_3(x)$ and $q$-integral representations of solutions are described in Appendix~\ref{appsec:B}.

\subsubsection[Extension of variant of q-hypergeometric equation of degree 3]{Extension of variant of $\boldsymbol{q}$-hypergeometric equation of degree 3} \label{subsec:3extedeg3}

In this subsection, we apply the addition with respect to the gauge-transformation
\begin{gather*}
 Y_g(x) = \frac{(\gamma_1 x, \gamma_2 x;q)_{\infty}}{(\alpha_1 x, \alpha_2 x;q)_{\infty}}\overline{Y}(x)
\end{gather*}
to the tuple $\bigl(O, \overline{G}_1, \overline{G}_2 , \overline{G}_3 \bigr) $ in equation \eqref{eq:v-qhg-deg3}.
Then, we obtain the tuple $(B'_0, B'_1, B'_2 , B'_3 )$ and the $q$-difference equation
\begin{gather*} 
\frac{Y_g(qx) - Y_g(x)}{-x} = \bigg[ \frac{B'_0}{x} + \frac{B'_1}{x-1/\gamma_1} + \frac{B'_2}{x-1/\gamma_2} + \frac{B'_3}{x-1/\alpha_3} \bigg] Y_g(x) ,
\end{gather*}
where
\begin{gather*}
B'_0 = O,
B'_1
= \begin{pmatrix}
 b^1_{11} & b^1_{12} \\
 b^1_{21} & b^1_{22}
 \end{pmatrix},\qquad
B'_2
= \begin{pmatrix}
 b^2_{11} & b^2_{12} \\
 b^2_{21} & b^2_{22}
 \end{pmatrix},\qquad
B'_3
= \begin{pmatrix}
 b^3_{11} & b^3_{12} \\
 b^3_{21} & b^3_{22}
 \end{pmatrix} , \\
 b^1_{11} = \frac{1}{(\gamma_1-\gamma_2)(\alpha_3-\gamma_1)}
\bigl\{ -(\alpha_1\alpha_2+\alpha_2\alpha_3+\alpha_3\alpha_1) + \gamma_1(\alpha_1+\alpha_2+\alpha_3-\gamma_1) \\
\phantom{b^1_{11} =}{} + \frac{\alpha_1\alpha_2\alpha_3}{q^{\lambda}\gamma_1}
+ (\alpha_1-\alpha_3)(\alpha_2-\gamma_1) q^{-\lambda} B_1 + \alpha_1(\alpha_2+\alpha_3-\gamma_1)\bigl(1-q^{-\lambda}\bigr) \bigr\}, \\
 b^1_{12} = \frac{(\alpha_1-\alpha_3)(\alpha_2-\gamma_1)}{(\gamma_1-\gamma_2)(\alpha_3-\gamma_1)} q^{-\lambda} B_2, \qquad b^1_{21} = \frac{(\alpha_2-\alpha_3)(\alpha_1-\gamma_1)}{(\gamma_1-\gamma_2)(\alpha_3-\gamma_1)} q^{-\lambda} B_1, \\
 b^1_{22} = \frac{1}{(\gamma_1-\gamma_2)(\alpha_3-\gamma_1)}
\bigl\{ -(\alpha_1\alpha_2+\alpha_2\alpha_3+\alpha_3\alpha_1) + \gamma_1(\alpha_1+\alpha_2+\alpha_3-\gamma_1) \\
\phantom{b^1_{22} =}{} + \frac{\alpha_1\alpha_2\alpha_3}{q^{\lambda}\gamma_1}
+ (\alpha_2-\alpha_3)(\alpha_1-\gamma_1) q^{-\lambda} B_2 + \alpha_2(\alpha_3+\alpha_1-\gamma_1)\bigl(1-q^{-\lambda}\bigr) \bigr\}, \\
 b^2_{11} = \frac{1}{(\gamma_1-\gamma_2)(\gamma_2-\alpha_3)}
\bigl\{ -(\alpha_1\alpha_2+\alpha_2\alpha_3+\alpha_3\alpha_1) + \gamma_2(\alpha_1+\alpha_2+\alpha_3-\gamma_2) \\
\phantom{b^1_{22} =}{}+ \frac{\alpha_1\alpha_2\alpha_3}{q^{\lambda}\gamma_2}
+ (\alpha_1-\alpha_3)(\alpha_2-\gamma_2) q^{-\lambda} B_1 + \alpha_1(\alpha_2+\alpha_3-\gamma_2)\bigl(1-q^{-\lambda}\bigr) \bigr\}, \\
 b^2_{12} = \frac{(\alpha_1-\alpha_3)(\alpha_2-\gamma_2)}{(\gamma_1-\gamma_2)(\gamma_2-\alpha_3)} q^{-\lambda} B_2, \qquad b^2_{21} = \frac{(\alpha_2-\alpha_3)(\alpha_1-\gamma_2)}{(\gamma_1-\gamma_2)(\gamma_2-\alpha_3)} q^{-\lambda} B_1, \\
 b^2_{22} = \frac{1}{(\gamma_1-\gamma_2)(\gamma_2-\alpha_3)}
\bigl\{ -(\alpha_1\alpha_2+\alpha_2\alpha_3+\alpha_3\alpha_1) + \gamma_2(\alpha_1+\alpha_2+\alpha_3-\gamma_2) \\
\phantom{b^1_{22} =}{} + \frac{\alpha_1\alpha_2\alpha_3}{q^{\lambda}\gamma_2}
+ (\alpha_2-\alpha_3)(\alpha_1-\gamma_2) q^{-\lambda} B_2 + \alpha_2(\alpha_3+\alpha_1-\gamma_2)\bigl(1-q^{-\lambda}\bigr) \bigr\}, \\
 b^3_{11} = b^3_{21} = -\frac{(\alpha_2-\alpha_3)(\alpha_3-\alpha_1)}{(\gamma_2-\alpha_3)(\alpha_3-\gamma_1)} q^{-\lambda} B_1, \\ b^3_{12} = b^3_{22} = -\frac{(\alpha_2-\alpha_3)(\alpha_3-\alpha_1)}{(\gamma_2-\alpha_3)(\alpha_3-\gamma_1)} q^{-\lambda} B_2,
\end{gather*}
and $B_1$, $B_2$, $B_3$ are defined in equation \eqref{eq:B1B2B3}.
Then, we have
\[
c^q_{\lambda'} (O, B'_1, B'_2, B'_3 ) = (G'_0, G'_1, G'_2, G'_3),
\]
 where the tuple $(G'_0, G'_1, G'_2, G'_3)$ is the $ 8\times 8$ matrices given by the $q$-convolution $c^q_{\lambda'} $.
The vectors~${}^t(1, 0)$ and ${}^t(0, 1)$ are a basis for the space $\ker (B'_0)$, and the vector ${}^t(-B_2, B_1)$ belongs to the space $\ker (B'_3)$.
Hence, we have $\dim \mathcal{K} \geq 3 $.
If $q^{\lambda+\lambda'}=\alpha_1\alpha_2/(\gamma_1\gamma_2)$,
then the vectors ${}^t(1,0,1,0,1,0,1,0)$ and ${}^t(0,1,0,1,0,1,0,1)$ can be taken as a basis of the vector space $\mathcal{L}=\ker(G'_0 + G'_1 + G'_2 + G'_3)$.
Therefore, we have $\dim \mathcal{K}=3$ and $\dim \mathcal{L}=2$.
Write
\[
mc^q_{\lambda'} (O, B'_1, B'_2, B'_3) = \bigl(\overline{G'}_0, \overline{G'}_1, \overline{G'}_2, \overline{G'}_3\bigr).
\]
By the simultaneous transformation for $(G'_0, G'_1,\allowbreak G'_2, G'_3)$ by the matrix
\begin{equation*}
P = \begin{pmatrix}
 0 & 0 & 0 & 1 & 0 & 0 & 1 & 0 \\
 0 & 0 & 0 & 0 & 1 & 0 & 0 & 1 \\
 1 & 0 & 0 & 0 & 0 & 0 & 1 & 0 \\
 0 & 1 & 0 & 0 & 0 & 0 & 0 & 1 \\
 0 & 0 & 0 & 0 & 0 & 0 & 1 & 0 \\
 0 & 0 & 0 & 0 & 0 & 0 & 0 & 1 \\
 0 & 0 & 1 & 0 & 0 & -B_2 & 1 & 0 \\
 0 & 0 & 0 & 0 & 0 & B_1 & 0 & 1
\end{pmatrix},
\end{equation*}
we obtain the expression
\begin{align*}
&\overline{G'}_0 = O , \qquad \overline{G'}_1
= \begin{pmatrix}
 1 + q^{-\lambda'}\bigl(b^1_{11}-1\bigr) & q^{-\lambda'}b^1_{12} & q^{-\lambda'}b^3_{11} \\
 q^{-\lambda'}b^1_{21} & 1 + q^{-\lambda'}\bigl(b^1_{22}-1\bigr) & q^{-\lambda'}b^3_{21} \\
 0 & 0 & 0
 \end{pmatrix}, \\
&\overline{G'}_2
= \begin{pmatrix}
 -q^{-\lambda'}b^1_{21} & -q^{-\lambda'}b^1_{12} & -q^{-\lambda'}b^3_{11} \\
 -q^{-\lambda'}b^1_{21} & -q^{-\lambda'}b^1_{22} & -q^{-\lambda'}b^3_{21} \\
 -\displaystyle \frac{b^1_{11}}{q^{\lambda'}}-\frac{B_2 b^1_{21}}{q^{\lambda'}B_1}
 & -\displaystyle \frac{b^1_{12}}{q^{\lambda'}} -\frac{B_2 b^1_{22}}{q^{\lambda'}B_1}
 & -\displaystyle \frac{b^3_{11}}{q^{\lambda'}}-\frac{B_2 b^3_{21}}{q^{\lambda'}B_1}
 \end{pmatrix}, \\
&\overline{G'}_3
= \begin{pmatrix}
 0 & 0 & 0 \\
 0 & 0 & 0 \\
 \displaystyle \frac{b^1_{11}}{q^{\lambda'}}+\frac{B_2 b^1_{21}}{q^{\lambda'}B_1}
 & \displaystyle \frac{b^1_{12}}{q^{\lambda'}}+\frac{B_2 b^1_{22}}{q^{\lambda'}B_1}
 & 1+\displaystyle \frac{b^3_{11}-1}{q^{\lambda'}}+\frac{B_2 b^3_{21}}{q^{\lambda'}B_1}
 \end{pmatrix}.
\end{align*}
The equation obtained after applying $q$-middle convolution $mc^q_{\lambda'}$ is written as
\begin{equation} \label{eq:3-111-111-21-3}
\frac{1}{-x}
\begin{pmatrix}
 g_1(qx) - g_1(x) \\
 g_2(qx) - g_2(x) \\
 g_3(qx) - g_3(x)
\end{pmatrix}
= \bigg[ \frac{\overline{G'}_1}{x-1/\gamma_1} + \frac{\overline{G'}_2}{x-1/\gamma_2} + \frac{\overline{G'}_3}{x-1/\alpha_3} \bigg]
 \begin{pmatrix}
 g_1(x) \\
 g_2(x) \\
 g_3(x)
 \end{pmatrix}.
\end{equation}
We derive single third-order $q$-difference equations from equation \eqref{eq:3-111-111-21-3}.
The functions $g_1(x)$ and $g_2(x)$ satisfy single third-order $q$-difference equations whose coefficients are polynomials in~$x$ of degree $6$, and $g_3(x)$ satisfies one whose coefficients are polynomials in~$x$ of degree $5$.

The single third-order $q$-difference equation satisfied by $g_3(x)$ and $q$-integral representations of solutions are described in Appendix \ref{appsec:B}.

\subsection{Spectral type of equation}\label{subsec:specttype}

Sakai and Yamaguchi obtained other important results in the paper \cite{SY}, where the $q$-middle convolution was launched.
Namely, they introduced the spectral type and the index of rigidity for some class of linear $q$-difference equations, and they established that the index of rigidity is preserved by the $q$-middle convolution.
On the system of the $q$-difference equations
\[
 Y(qx)= B(x) Y(x), \qquad B(x)= B_{\infty}+\sum_{i=1}^N \frac{B_i}{1-x/b_i} ,
 \]
the spectral type $( S_0; S_{\infty }; S_{\rm div} )$ is determined by the Jordan normal forms of $B_0 = I_m - B_\infty - B_{1} - \dots -B_N$ and $B_{\infty } $, and the Smith normal form of the matrix $A(x)= B(x) \prod _{i=1} ^N (1-x/b_i)$ whose elements are polynomials in $x$.
If $B_0$ (resp.~$B_{\infty}$) is diagonalizable, then $S_0$ (resp.~$S_{\infty}$) is the partition whose element is the multiplicity of each eigenvalue of $B_0$ (resp.~$B_{\infty}$).
If the Smith normal form of the matrix $A(x)$ is written as
\[
\operatorname{diag}\left(1, \dots , 1, \prod _{i=1}^{k_3} (x-a_i), \prod _{i=1}^{k_2} (x-a_i), \prod _{i=1}^{k_1} (x-a_i)\right),
\] where $k_3 \leq k_2 \leq k_1$, then $S_{\rm div}$ is the partition $3\cdots 32\cdots 21\cdots 1$ where the total number of $3$ (resp.~$2$, $1$) is $k_3$ (resp.~$k_2-k_3$, $k_1-k_2$).
In particular, if the roots of $\det A(x) =0$ are different mutually, then $S_{\rm div} = 11\cdots 1$.
For the precise definition of the spectral type and the index of rigidity, see \cite{SY}.
The spectral type would be a nice labelling to the system of the linear $q$-difference equations.
In Figure~\ref{fg:spectraltype}, the spectral types of the system of the linear q-difference equations which appeared in this paper are described.

We explain the spectral type of equation \eqref{eq:3-111-21-111} as an example.
The matrix $\overline{G'}_0 $ in equation~\eqref{eq:3-111-21-111} is the zero matrix, which is a diagonal matrix with the eigenvalue $0$ of multiplicity~$3$.
The matrix $I_3 - \overline{G'}_0 -\overline{G'}_1 -\overline{G'}_2 $ is diagonalizable for generic parameters and each eigenvalue is of multiplicity $1$.
Hence, the part $3;111$ appears in the spectral type.
To interpret the part~$21111$, we write equation \eqref{eq:3-111-21-111} as $Y(qx)=B(x) Y(x)$ and set $C(x)= (x-1/\gamma _1) (x-1/\alpha _2 ) B(x)$.
Then, each coefficient of $C(x)$ is a polynomial of degree $2$, and the matrix $C(x)$ is equivalent to
\begin{equation*}
\begin{pmatrix}
 1 & 0 & 0 \\
 0 & x-a_1 & 0 \\
 0 & 0 & (x-a_1) (x-a_2)(x-a_3)(x-a_4)(x-a_5)
 \end{pmatrix}
\end{equation*}
as the Smith normal form of the polynomial matrix.
Thus, the part $21111$ appears in the spectral type as the multiplicities of the Smith normal form.
Note that the index of rigidity for the cases in Figure \ref{fg:spectraltype} is equal to $2$ (see \cite{SY} for the definition of the index of rigidity), and we can regard these equations rigid.
\begin{figure}[t]
\centering
\begin{tikzpicture}
\node (a1) at (0,0) { $\underset{\eqref{eq:qeab}}{1;1;1}$ };
\node (a2) at (4,0) { $\underset{\eqref{eq:qhg-sta}}{11;11;11}$};
\node (a3) at (8,0) { $\underset{\eqref{eq:111-21-111}}{111;111;21}$};
\node (a4) at (8,-2) { $\underset{\eqref{eq:21-111-111}}{21;111;111}$};
\draw[->] (a1) -- (a2);
\draw[->] (a2) -- (a3);
\draw[->] (a2) -- (a4);
\node (b1) at (0,-4) { $\underset{\eqref{eq:qeqvar2}}{1;1;11}$};
\node (b2) at (4,-4) { $\underset{\eqref{eq:v-qhg-deg2}}{2;11;1111}$};
\node (b3) at (8,-4) { $\underset{\eqref{eq:3-111-21-111}}{3;111;21111}$};
\node (b4) at (8,-6) { $\underset{\eqref{eq:3-111-111-21}}{3;21;111111}$};
\node (b5) at (4,-7) { $\underset{\eqref{eq:G0G1GN},\ N=2}{21;21;2211}$};
\draw[->] (b1) -- (b2);
\draw[->] (b2) -- (b3);
\draw[->] (b2) -- (b4);
\draw[->] (b1) -- (b5);
\node (c1) at (0,-9) { $\underset{\eqref{eq:B1B2B3}}{1;1;111}$};
\node (c2) at (4,-9) { $\underset{\eqref{eq:v-qhg-deg3}}{2;2;111111}$};
\node (c3) at (8,-9) { $\underset{\eqref{eq:3-111-111-21-3}}{3;3;21111111}$};
\node (c5) at (4,-11) { $\underset{\eqref{eq:G0G1GN},\ N=3}{31;31;333111}$};
\draw[->] (c1) -- (c2);
\draw[->] (c2) -- (c3);
\draw[->] (c1) -- (c5);
\end{tikzpicture}
\caption{Change of the spectral type.}\label{fg:spectraltype}
\end{figure}

On the other hand, the spectral type of the system of the $q$-difference equations written as
\begin{equation*}
\frac{Y(q x) - Y(x )}{-x} = \Biggl[ \frac{B_{0}}{x} + \frac{B_{1}}{x -b_{1}} + \frac{B_{2}}{x -b_{2}} \Biggr] Y(x) 
\end{equation*}
for generic $2\times 2$ matrices $B_0$, $B_1$, $B_2 $ is $11;11;1111$.
This equation is related to the sixth $q$-difference Painlev\'e equation in \cite{JS} and the $q$-Heun equation in \cite{STT,TakR,TakqH}.
On this equation, the index of rigidity is equal to $0$, and it indicates that a nature of this equation would be different from that of the equations appeared in this paper.

\section[Composition of q-middle convolution]{Composition of $\boldsymbol{q}$-middle convolution} \label{sec:compqmc}

In this section, we establish the additivity of the $q$-middle convolution in equation \eqref{eq:mcaddi0} by refocusing the arguments in Dettweiler and Reiter \cite{DR1, DR2}.
A key formula for the additivity is equation \eqref{eq:Gjql1l2GDR}, which is valid because the definition of the $q$-middle convolution is reformulated.

We change notations around the middle convolution to have more precise description.
Let~$V$ be a finite-dimensional vector space and $\boldsymbol{B} = ( B_0, B_1, \dots , B_N)$ be a tuple of endomorphisms (linear transformations) of $V$.
By fixing a basis of the vector space $V$, we obtain the tuple of square matrices of size $\dim V$ from the tuple of endomorphisms of $V$.

In \cite{DR1,DR2}, Dettweiler and Reiter discussed several properties of the middle convolution, which are related to the following definition.
\begin{Definition} \label{def:condqmc}
Let $V$ and $W$ be finite-dimensional vector spaces and let $\boldsymbol{B} = ( B_0, B_1, \dots , B_N)$ (resp.~$\boldsymbol{C} = ( C_0, C_1, \dots , C_N)$) be a tuple of endomorphisms of~$V$ (resp.~$W$).
\begin{itemize}\itemsep=0pt
\item[(i)] $(V, \boldsymbol{B} )$ is isomorphic to $(W, \boldsymbol{C})$, if there exists an isomorphism $\phi \colon V \to W$ of the vector spaces such that $\phi \circ B_j = C_j \circ \phi $ for $j=0,1,\dots ,N$.
\item[(ii)] $(V,\boldsymbol{B} )$ is irreducible, if the subspace $W$ of $V$ such that $B_j W \subset W$ for all $j=0,1,\dots ,N$ is only $W=V$ or $W=\{ 0 \} $.
\item[(iii)] $(V,\boldsymbol{B} )$ satisfies $(*)$, if
\begin{equation}
 \forall i=0,1,\dots ,N,\quad \forall \tau \in \mathbb{C} , \qquad \Biggl( \bigcap _{i' =0 \atop{i'\neq i} }^N \operatorname{ker}(B_{i'} ) \Biggr) \bigcap \operatorname{ker}(B_{i} +\tau ) =\{ 0\}.
\label{eq:def*}
\end{equation}
\item[(iv)] $(V,\boldsymbol{B} )$ satisfies $(**)$, if
\begin{equation}
 \forall i=0,1,\dots ,N, \quad \forall \tau \in \mathbb{C} , \qquad \Biggl( \sum_{i' =0 \atop{i'\neq i} }^N \operatorname{im}(B_{i'} ) \Biggr) + \operatorname{im}(B_{i} +\tau ) = V.
\label{eq:def**}
\end{equation}
\end{itemize}
\end{Definition}

Set $V_{j} =V$ $(j=0,\dots ,N)$ and $V' = V_0 \oplus V_{1} \oplus \dots \oplus V_{N} \bigl(= V^{\oplus N+1}\bigr)$.
Write $u \in V'$ as
\begin{equation*}
 u = \begin{pmatrix} u_{0} \\ \vdots \\ u_{N} \end{pmatrix}, \qquad u_j \in V_j , \quad j=0,\dots , N.
\end{equation*}
Let $\lambda \in \mathbb{C} $, $j \in \{ 0,1, \dots ,N \}$ and $G_j^{q} (\lambda )$ be an endomorphism of $V'$ which corresponds to the $q$-convolution given in equation \eqref{eq:bG}.
Then, it admits the expression
\begin{gather}
G_j^{q} (\lambda ) \begin{pmatrix} u_{0} \\ \vdots \\ u_{N} \end{pmatrix} =\begin{pmatrix} w_{0} \\ \vdots \\ w_{N} \end{pmatrix}, \qquad w_j = \bigl(1- q^{-\lambda } \bigr) u_j + q^{-\lambda } \sum _{i=0}^{N} B_i u_i, \nonumber\\ w _k =0, \qquad k\neq j.
\label{eq:Gjq}
\end{gather}
Write the $q$-convolution as
$
c_\lambda ^{q} (B_0, \dots , B_N ) = \bigl( G_0^{q} (\lambda ), \dots , G_N^{q} (\lambda )\bigr) $.
Let $ \mathcal{K}^{q} $ and $\mathcal{L} ^{q} ( \lambda ) $ be the subspaces of $V'$ defined in equation~\eqref{eq:qKL}.
They are expressed as
\begin{align*}
& \mathcal{K}^{q} = \operatorname{ker}B_0 \oplus \operatorname{ker}B_1 \oplus \cdots \oplus \operatorname{ker}B_N, \\
& \mathcal{L}^{q} ( \lambda ) = \operatorname{ker}\bigl(G_0^{q} (\lambda ) +G_1^{q} (\lambda ) + \cdots + G_N^{q} (\lambda ) \bigr).
\end{align*}
We denote the action of $G_j^{q} (\lambda ) $ on the quotient space $ {\mathcal M } = V^{\oplus N+1}/( \mathcal{K}^{q} + \mathcal{L}^{q} ( \lambda ) )$ by $\overline{G}_j^{q} (\lambda )$ $(j=0 ,1, \dots ,N )$.
Then, the $q$-middle convolution $mc_\lambda ^{q}$ is defined by the correspondence of the tuples of endomorphisms $ (B_{0}, B_{1} ,\dots ,B_N ) \mapsto \bigl( \overline{G}_0^{q} (\lambda ) , \overline{G}_1^{q} (\lambda ) ,\dots ,\overline{G}_N^{q} (\lambda )\bigr) $.
Write $mc_\lambda ^{q} (V) = V^{\oplus N+1}/( \mathcal{K}^{q} + \mathcal{L}^{q} ( \lambda ) ) $.

We express the convolution and the middle convolution of Dettweiler and Reiter in our notation.
Let $\lambda \in \mathbb{C} $ and $j \in \{ 0,1, \dots ,N \}$.
Let $G_j (\lambda )$ be an endomorphism of $V'$ which corresponds to the convolution given in equation \eqref{eq:Fi}.
It is expressed as
\begin{equation}
G_j (\lambda ) \begin{pmatrix} u_{0} \\ \vdots \\ u_{N} \end{pmatrix} =\begin{pmatrix} w_{0} \\ \vdots \\ w_{N} \end{pmatrix}, \qquad w_j = \lambda u_j + \sum _{i=0}^{N} B_i u_i , \qquad w _k =0,\quad k\neq j.
\label{eq:GjDR}
\end{equation}
Write the convolution as $c_\lambda (B_0, \dots , B_N ) = ( G_0 (\lambda ), \dots , G_N (\lambda )) $.
Let $ \mathcal{K} $ and $\mathcal{L} ( \lambda ) $ be the subspaces of $V'$ such that
\begin{align*}
\begin{split}
& \mathcal{K} = \operatorname{ker}B_0 \oplus \operatorname{ker}B_1 \oplus \cdots \oplus \operatorname{ker}B_N, \\
& \mathcal{L} ( \lambda ) = \operatorname{ker}(G_0 (\lambda ) +G_1 (\lambda ) + \cdots + G_N (\lambda ) ).
\end{split}
\end{align*}
We denote the action of $G_j (\lambda ) $ on the quotient space $ {\mathcal M } = V^{\oplus N+1}/( \mathcal{K} + \mathcal{L} ( \lambda ) )$ by $\overline{G}_j (\lambda )$ $(j=0 ,1, \dots ,N )$.
Then, the middle convolution $mc_\lambda $ is defined by the correspondence of the tuples of endomorphisms $ (B_{0}, B_{1} ,\dots ,B_N ) \mapsto \bigl( \overline{G}_0 (\lambda ) , \overline{G}_1 (\lambda ) ,\dots ,\overline{G}_N (\lambda )\bigr) $.
Write $mc_\lambda (V) = V^{\oplus N+1}/( \mathcal{K} + \mathcal{L} ( \lambda ) ) $.
We describe some results obtained by Dettweiler and Reiter \cite{DR1} in our notation.
\begin{Proposition}[{\cite[Proposition 3.4]{DR1}}] \label{prop:*,**}
If $(V ,\boldsymbol{B}) $ satisfies $(*)$ {\rm(}resp.~$(**))$, then $(mc _{\mu }(V),\allowbreak mc _{\mu } (\boldsymbol{B} ))$ also satisfies $(*)$ {\rm(}resp.~$(**))$.
\end{Proposition}
\begin{Proposition}[{\cite[Proposition 3.2]{DR1}}] \label{prop:mc0VV}
If $(V ,\boldsymbol{B}) $ satisfies $(**)$, then $(mc _{0} (V), mc _{0}(\boldsymbol{B})) $ is isomorphic to $(V ,\boldsymbol{B}) $.
\end{Proposition}
\begin{Theorem}[{\cite[Theorem 3.5]{DR1}}] \label{thm:mc12isom}
If $(V ,\boldsymbol{B}) $ satisfies $(*)$ and $(**)$, then
\[
(mc _{\mu _2} (mc _{\mu _1} (V)), mc _{\mu _2} (mc _{\mu _1} (\boldsymbol{B} )))
\]
is isomorphic to $(mc _{\mu _1 +\mu _2} (V), mc _{\mu _1 +\mu _2} (\boldsymbol{B} ))$ for any $\mu _1, \mu _2 \in \mathbb{C} $.
Moreover,
\[
(mc _{-\mu } (mc _{\mu } (V)), mc _{-\mu } (mc _{\mu } (\boldsymbol{B} )))
\]
 is isomorphic to $(V, \boldsymbol{B} )$ for any $\mu \in \mathbb{C} $.
\end{Theorem}
\begin{Theorem}[{\cite[Corollary 3.6]{DR1}}] \label{thm:irr}
If $(V,\boldsymbol{B})$ is irreducible, then $(mc _{\mu }(V), mc _{\mu } (\boldsymbol{B} ))$ is irreducible or $V=\{ 0 \}$.
\end{Theorem}

Some properties of the $q$-middle convolution are obtained by comparing the description of the $q$-deformed ones with the original ones.
It follows from equations \eqref{eq:Gjq} and \eqref{eq:GjDR} that
\begin{equation}
 G_j^{q} (\lambda ) = q^{-\lambda } G_j \bigl(q^{\lambda } -1 \bigr),\qquad j=0,\dots ,N,\qquad \mathcal{K}^{q} = \mathcal{K} , \qquad \mathcal{L}^{q} ( \lambda ) = \mathcal{L} \bigl(q^{\lambda } -1\bigr).  \label{eq:GjqlGjDR}
\end{equation}
The space \smash{$mc_\lambda ^{q} (V) = V^{\oplus N+1}/( \mathcal{K}^{q} + \mathcal{L}^{q} ( \lambda ) ) $} coincides with the space $\smash{mc_{q^{\lambda } -1} (V) = V^{\oplus N+1}}/\bigl( \mathcal{K} + \mathcal{L} \bigl(q^{\lambda } -1\bigr) \bigr) $, and we have
\begin{equation}
 \overline{G}_j^{q} (\lambda ) = q^{-\lambda } \overline{G}_j \bigl(q^{\lambda } -1 \bigr), \qquad j=0,\dots ,N. \label{eq:GjqlGjDRb}
\end{equation}
The following propositions are obtained by Propositions \ref{prop:*,**} and \ref{prop:mc0VV}, and the correspondence above.
\begin{Proposition}
The conditions $(*)$ and $(**)$ are preserved by the $q$-middle convolution $mc_\lambda ^{q} $.
\end{Proposition}
\begin{Proposition}[{cf.~\cite[Proposition 3.2]{DR1}}] \label{prop:qmc0VV}
If $(V ,\boldsymbol{B}) $ satisfies $(**)$, then $\bigl(mc ^{q}_{0} (V), mc ^{q}_{0}(\boldsymbol{B})\bigr) $ is isomorphic to $(V ,\boldsymbol{B}) $.
\end{Proposition}

We are going to consider composition of the middle convolutions, which is related to Theorem~\ref{thm:mc12isom}.
We apply the $q$-convolution with the parameter $\lambda _2 $ to the tuples of endomorphisms ${ \bigl( G_0^{q} (\lambda _1) , G_1^{q} (\lambda _1) ,\dots , G_N^{q} (\lambda _1 )\bigr) \bigl(= c _{\lambda _1}^q(B_0, \dots ,B_N)\bigr)}$ on the space $V' =V^{\oplus N+1}$.
Write
\begin{equation*}
c _{\lambda _2}^q\bigl(c _{\lambda _1}^q(B_0, \dots ,B_N)\bigr) = \bigl( G_0^{q} (\lambda _1, \lambda _2 ) , \dots , G_N^{q} (\lambda _1, \lambda _2 ) \bigr).
\end{equation*}
We denote an element $\mathbf{u} $ in $ (V') ^{\oplus N+1} \bigl(= \bigl(V^{\oplus N+1}\bigr) ^{\oplus N+1}\bigr)$ by
\begin{equation*}
 \mathbf{u} = \begin{pmatrix} \mathbf{u}_{0} \\ \vdots \\ \mathbf{u}_{N} \end{pmatrix} \in (V') ^{\oplus N+1} , \qquad \mathbf{u}_{j} = \begin{pmatrix} \mathbf{u}_{j,0} \\ \vdots \\ \mathbf{u}_{j,N} \end{pmatrix} \in V ^{\oplus N+1}, \qquad j=0,\dots , N.
\end{equation*}
Let $j \in \{ 0,\dots ,N \}$.
The action of $G_j^{q} (\lambda _1, \lambda _2 ) $ on the space $(V') ^{\oplus N+1} $ is written as
\begin{gather*}
G_j^{q} (\lambda _1, \lambda _2 ) \begin{pmatrix} \mathbf{u}_{0} \\ \vdots \\ \mathbf{u}_{N} \end{pmatrix} =\begin{pmatrix} \mathbf{w}_{0} \\ \vdots \\ \mathbf{w}_{N} \end{pmatrix}, \qquad \mathbf{w}_j =(1- q^{-\lambda _2} ) \mathbf{u}_j + q^{-\lambda _2} \sum _{i=0}^{N} G_i^{q} (\lambda _1) \mathbf{u}_i , \\ \mathbf{w}_k =0 , \qquad k\neq j,
\end{gather*}
and $\mathbf{w}_j $ is written as
\begin{gather}
\mathbf{w}_{j} = \begin{pmatrix} \mathbf{w}_{j,0} \\ \vdots \\ \mathbf{w}_{j,N} \end{pmatrix} ,\nonumber\\ \mathbf{w}_{j,k} = \bigl(1- q^{-\lambda _2} \bigr) \mathbf{u}_{j,k} + q^{-\lambda _2} \bigl(1- q^{-\lambda _1} \bigr) \mathbf{u}_{k,k} + q^{-\lambda _2} q^{-\lambda _1} \sum _{l=0}^{N} B_{l} \mathbf{u}_{k,l} \label{eq:wjq}
.
\end{gather}
We investigate composition of two convolutions.
Let $\lambda _1 , \lambda _2 \in \mathbb{C} $ and write the composition of two convolutions as
\[
c _{\lambda _2} (c _{\lambda _1} (B_0, \dots ,B_N)) = ( G_0 (\lambda _1, \lambda _2 ) , \dots , G_N (\lambda _1, \lambda _2 ) ).
\]
Then, $G_j (\lambda _1, \lambda _2 )$ is an endomorphism of \smash{$(V' )^{\oplus N+1}$}, and it admits the expression
\begin{equation*}
G_j (\lambda _1, \lambda _2 ) \begin{pmatrix} \mathbf{u}_{0} \\ \vdots \\ \mathbf{u}_{N} \end{pmatrix} =\begin{pmatrix} \mathbf{w}_{0} \\ \vdots \\ \mathbf{w}_{N} \end{pmatrix}, \qquad \mathbf{w}_j = \lambda _2 \mathbf{u}_j + \sum _{i=0}^{N} G_i (\lambda _1) \mathbf{u}_i ,\qquad \mathbf{w}_k =0, \quad k\neq j,
\end{equation*}
and $\mathbf{w}_j $ is written as
\begin{equation}
\mathbf{w}_{j} = \begin{pmatrix} \mathbf{w}_{j,0} \\ \vdots \\ \mathbf{w}_{j,N} \end{pmatrix} ,\qquad \mathbf{w}_{j,k} = \lambda _2 \mathbf{u}_{j,k} + \lambda _1 \mathbf{u}_{k,k} + \sum _{l=0}^{N} B_{l} \mathbf{u}_{k,l}. \label{eq:wjDR}
\end{equation}
By combining equation \eqref{eq:wjDR} with equation \eqref{eq:wjq}, we have
\begin{equation}
G_j^{q} (\lambda _1, \lambda _2 ) = q^{-\lambda _1 -\lambda _2} G_j \bigl(q^{\lambda _1} -1 , q^{\lambda _1 + \lambda _2} - q^{\lambda _1 } \bigr) ,\qquad j=0,\dots ,N.
\label{eq:Gjql1l2GDR}
\end{equation}

Let us describe the vector space $mc ^{q} _{\lambda _2} (mc ^{q}_{\lambda _1} (V))$, which is obtained by composition of two $q$-middle convolutions.
Set
\[
 {\mathcal M } = mc ^{q} _{\lambda _1} (V) = V^{\oplus N+1}/( \mathcal{K}^{q} + \mathcal{L}^{q} (\lambda _1 ) ) .
 \]
Write a representative element \smash{$\mathbf{v} \bigl(\in \bigl(V^{\oplus N+1}\bigr) ^{\oplus N+1}\bigr)$} in the space ${\mathcal M } ^{\oplus N+1} $ as
\begin{equation*}
 \mathbf{v} = \begin{pmatrix} \mathbf{v}_{0} \\ \vdots \\ \mathbf{v}_{N} \end{pmatrix}.
\end{equation*}
The condition that a representative element $\mathbf{v} $ belongs to the space
\[
{\mathcal K }_{\mathcal M }^{q} = \operatorname{ker} \overline{G}_0^{q} (\lambda _1 ) \oplus \dots \oplus \operatorname{ker} \overline{G}_N^{q} (\lambda _1 )
\]
 is
\begin{equation}
 G_j^{q} (\lambda _1 ) \mathbf{v}_j \in \mathcal{K}^{q} + \mathcal{L}^{q} (\lambda _1 ) , \qquad j=0,\dots ,N , \label{eq:deftK}
\end{equation}
and the condition that a representative element $\mathbf{v}$ belongs to the space
\[
{\mathcal L }_{\mathcal M } ^{q} (\lambda _2 ) = \operatorname{ker} \bigl(\bigl(G_0^{q} (\lambda _1,\allowbreak \lambda _2 ) + \dots + G_N^{q} (\lambda _1, \lambda _2 )\bigr)|_{{\mathcal M } ^{\oplus N+1} }\bigr)
\]
 is
\begin{equation}
 \bigl( G_0^{q} (\lambda _1, \lambda _2 ) + \dots + G_N^{q} (\lambda _1, \lambda _2 )\bigr) \mathbf{v} \in \bigl(\mathcal{K}^{q} + \mathcal{L}^{q} (\lambda _1 )\bigr)^{\oplus N+1}.  \label{eq:deftL}
\end{equation}

Let $ \mathcal{K}^{q} (\lambda _1) $ and $\mathcal{L} ^{q} ( \lambda _1, \lambda _2 ) $ be the subspaces of \smash{$\bigl(V^{\oplus N+1}\bigr) ^{\oplus N+1}$} determined by equations \eqref{eq:deftK} and \eqref{eq:deftL} respectively.
Then, the space
\[
mc _{\lambda _2} ^{q} (mc _{\lambda _1} ^{q} (V)) = {\mathcal M }^{\oplus N+1} /({\mathcal K }_{\mathcal M } + {\mathcal L_{\mathcal M }} (\lambda _2 ) )
\]
 is isomorphic to \smash{$\bigl(V^{\oplus N+1}\bigr)^{\oplus N+1} /( \mathcal{K}^{q} (\lambda _1) + \mathcal{L} ^{q} ( \lambda _1, \lambda _2 ) ) $} as the vector space, because
\begin{gather*}
 mc _{\lambda _2} ^{q} (mc _{\lambda _1} ^{q} (V)) \\
 \qquad= mc _{\lambda _2} ({\mathcal M }) = {\mathcal M }^{\oplus N+1} /({\mathcal K }_{\mathcal M } + {\mathcal L_{\mathcal M }} (\mu _2 ) ) \\
 \qquad= \bigl(\bigl(V^{\oplus N+1} \bigr)^{\oplus N+1} \!/ \bigl(\mathcal{K}^{q} + \mathcal{L}^{q} (\lambda _1 )\bigr) ^{\oplus N+1} \bigr) /\bigl(( \mathcal{K}^{q} (\lambda _1) + \mathcal{L} ^{q} ( \lambda _1, \lambda _2 ) )/(\mathcal{K}^{q} + \mathcal{L}^{q} (\lambda _1 )) ^{\oplus N+1} \bigr) \\
 \qquad\simeq \bigl(V^{\oplus N+1}\bigr)^{\oplus N+1} /( \mathcal{K}^{q} (\lambda _1) + \mathcal{L} ^{q} ( \lambda _1, \lambda _2 ) ).
\end{gather*}
The action of $mc ^{q}_{\lambda _2} (mc^{q}_{\lambda _1} (B_0, \dots ,B_N)) $ on the space $mc ^{q}_{\lambda _2} (mc ^{q}_{\lambda _1} (V)) $ is induced by the tuple $(G_0^{q} (\lambda _1, \lambda _2 ), \dots , G_N^{q} (\lambda _1, \lambda _2 ) )$ acting on \smash{$ \bigl(V^{\oplus N+1}\bigr)^{\oplus N+1} $}, and we write
\begin{equation} \label{eq:mcqmcq}
mc ^{q}_{\lambda _2} \bigl(mc^{q}_{\lambda _1} (B_0, \dots ,B_N)\bigr) = \bigl(\overline{G}_0^{q} (\lambda _1, \lambda _2 ), \dots , \overline{G}_N^{q} (\lambda _1, \lambda _2 )\bigr).
\end{equation}
The vector space $mc _{\lambda _2} (mc _{\lambda _1} (V))$, which is obtained by composition of two middle convolutions, is described as
\begin{equation*}
 mc _{\lambda _2} (mc _{\lambda _1} (V)) \simeq \bigl(V^{\oplus N+1}\bigr)^{\oplus N+1} /( \mathcal{K} (\lambda _1) + \mathcal{L} ( \lambda _1, \lambda _2 ) ) ,
\end{equation*}
where the space $\mathcal{K} (\lambda _1)$ (resp.~$\mathcal{L} ( \lambda _1, \lambda _2 )$) is characterized by the condition
\begin{gather}
 G_j (\lambda _1 ) \mathbf{v}_j \in \mathcal{K} + \mathcal{L} (\lambda _1 ) , \qquad j=0,\dots ,N , \nonumber\\
 \bigl(\mbox{resp. } ( G_0 (\lambda _1, \lambda _2 ) + \dots + G_N (\lambda _1, \lambda _2 )) \mathbf{v} \in (\mathcal{K} + \mathcal{L} (\lambda _1 ))^{\oplus N+1} \bigr). \label{eq:deftKLDR}
\end{gather}
The action of $mc _{\lambda _2} (mc_{\lambda _1} (B_0, \dots ,B_N)) $ on the space $mc _{\lambda _2} (mc _{\lambda _1} (V))$ is induced by the tuple $(G_0 (\lambda _1, \lambda _2 ), \dots , G_N (\lambda _1, \lambda _2 ) )$ acting on \smash{$ \bigl(V^{\oplus N+1}\bigr)^{\oplus N+1} $}, and we write
\begin{equation} \label{eq:mcmc}
mc _{\lambda _2} (mc_{\lambda _1} (B_0, \dots ,B_N)) = \bigl(\overline{G}_0 (\lambda _1, \lambda _2 ), \dots , \overline{G}_N (\lambda _1, \lambda _2 ) \bigr).
\end{equation}
It follows from equations \eqref{eq:GjqlGjDR}, \eqref{eq:Gjql1l2GDR}, \eqref{eq:deftK}, \eqref{eq:deftL} and \eqref{eq:deftKLDR} that
\begin{equation*}
 \mathcal{K}^{q} (\lambda _1) = \mathcal{K} \bigl(q^{\lambda _1} -1 \bigr), \qquad \mathcal{L}^{q} ( \lambda _1, \lambda _2 ) = \mathcal{L} \bigl(q^{\lambda _1} -1 , q^{\lambda _1 + \lambda _2} - q^{\lambda _1 } \bigr).
\end{equation*}
Therefore, the vector space $mc _{\lambda _2} ^{q} \bigl(mc _{\lambda _1} ^{q} (V)\bigr) $ is isomorphic to $mc _{q^{\lambda _1 + \lambda _2} - q^{\lambda _1 }} (mc _{q^{\lambda _1} -1 } (V)) $, and it follows from equation \eqref{eq:Gjql1l2GDR} that
\begin{equation}
\overline{G}_j^{q} (\lambda _1, \lambda _2 ) = q^{-\lambda _1 -\lambda _2} \overline{G}_j \bigl(q^{\lambda _1} -1 , q^{\lambda _1 + \lambda _2} - q^{\lambda _1 } \bigr) ,\qquad j=0,\dots ,N.
\label{eq:Gjql1l2GDRb}
\end{equation}
The following theorem is obtained by Theorem~\ref{thm:mc12isom} \cite[Theorem 3.5]{DR1}.
\begin{Theorem} \label{thm:qmc12isom}
Assume that $(V ,\boldsymbol{B}) $ satisfies $(*)$ and $(**)$.
\begin{itemize}\itemsep=0pt
\item[$(i)$] $\bigl(mc _{\lambda _2} ^{q} \bigl(mc _{\lambda _1} ^{q} (V)\bigr), mc _{\lambda _2} ^{q} \bigl(mc _{\lambda _1} ^{q} (\boldsymbol{B} )\bigr)\bigr)$ is isomorphic to $\bigl(mc _{\lambda _1 +\lambda _2} ^{q} (V), mc _{\lambda _1 +\lambda _2} ^{q} (\boldsymbol{B} )\bigr)$ for any $\lambda _1, \lambda _2 \in \mathbb{C} $.

\item[$(ii)$] $\bigl(mc _{-\lambda } ^{q} \bigl(mc _{\lambda } ^{q} (V)\bigr), mc _{-\lambda } ^{q} \bigl(mc _{\lambda } ^{q} (\boldsymbol{B} )\bigr)\bigr) $ is isomorphic to $(V, \boldsymbol{B} )$ for any $\lambda \in \mathbb{C} $.
\end{itemize}
\end{Theorem}
\begin{proof}
It follows from Theorem~\ref{thm:mc12isom} that
\[
\bigl(mc _{q^{\lambda _1 + \lambda _2} - q^{\lambda _1 }} (mc _{q^{\lambda _1} -1} (V)) , mc _{ q^{\lambda _1 + \lambda _2} - q^{\lambda _1 }} (mc _{q^{\lambda _1} -1} (\boldsymbol{B} ))\bigr)
\]
 is isomorphic to $(mc _{q^{\lambda _1 + \lambda _2} -1 } (V), mc _{q^{\lambda _1 + \lambda _2} -1 } (\boldsymbol{B} ))$.
By combining with equations \eqref{eq:mcmc}, \eqref{eq:mcqmcq}, \eqref{eq:Gjql1l2GDRb} and \eqref{eq:GjqlGjDRb} for $\lambda = \lambda _1 +\lambda _2$, we obtain (i).
We set $\lambda _2 =-\lambda _1$ in (i).
Then, we obtain that~$\bigl(mc _{-\lambda _1} ^{q} \bigl(mc _{\lambda _1} ^{q} (V)\bigr), mc _{-\lambda _1} ^{q} \bigl(mc _{\lambda _1} ^{q} (\boldsymbol{B} )\bigr)\bigr)$ is isomorphic to $\bigl(mc _{0} ^{q} (V), mc _{0} ^{q} (\boldsymbol{B} )\bigr) $, and it follows from Proposition~\ref{prop:qmc0VV} that it is isomorphic to $(V, \boldsymbol{B} )$.
\end{proof}

The following theorem is obtained by Theorem~\ref{thm:irr} \cite[Corollary~3.6]{DR1}.
\begin{Theorem}
If $(V,\boldsymbol{B})$ is irreducible, then $\bigl(mc ^{q}_{\lambda }(V), mc ^{q}_{\lambda } (\boldsymbol{B} )\bigr)$ is irreducible or $V=\{ 0 \}$.
\end{Theorem}
We describe an isomorphism between $\bigl(mc _{\lambda _2} ^{q} (mc _{\lambda _1} ^{q} (V)), mc _{\lambda _2} ^{q} \bigl(mc _{\lambda _1} ^{q} (\boldsymbol{B} )\bigr)\bigr)$ and $\bigl(mc _{\lambda _1 +\lambda _2} ^{q} (V),\allowbreak mc _{\lambda _1 +\lambda _2} ^{q} (\boldsymbol{B} )\bigr)$ in the case $\lambda _1 \neq 0 \neq \lambda _2$.
\begin{Proposition}
Define the map \smash{$\phi ^{q} \colon \bigl(V^{\oplus N+1}\bigr)^{\oplus N+1} \rightarrow V^{\oplus N+1} $} by
\begin{equation*}
\phi ^{q} (\mathbf{v}) = \sum _{j=0}^N G_j^{q} (\lambda _1) \mathbf{v}_{j}.  
\end{equation*}
This map induces an isomorphism $\overline{\phi }^{q} $ between
\[
\bigl(mc _{\lambda _2} ^{q} \bigl(mc _{\lambda _1} ^{q} (V)\bigr), mc _{\lambda _2} ^{q} \bigl(mc _{\lambda _1} ^{q} (\boldsymbol{B} )\bigr)\bigr)\qquad \text{and} \qquad\bigl(mc _{\lambda _1 +\lambda _2} ^{q} (V),
mc _{\lambda _1 +\lambda _2} ^{q} (\boldsymbol{B} )\bigr),
\]
 if $\lambda _1 \neq 0 \neq \lambda _2$.
\end{Proposition}
\begin{proof}
It was shown in \cite{DR1} (see also \cite[Section~7.5]{Har}) that the map \smash{$\phi \colon \bigl(V^{\oplus N+1}\bigr)^{\oplus N+1} \rightarrow V^{\oplus N+1} $} defined by
\begin{equation*}
\phi (\mathbf{v}) = \sum _{j=0}^N G_j (\mu _1) \mathbf{v}_{j} 
\end{equation*}
induces an isomorphism $\overline{\phi } $ between
\[
(mc _{\mu _2} (mc _{\mu _1} (V)), mc _{\mu _2} (mc _{\mu _1} (\boldsymbol{B} )))\qquad \text{and} \qquad (mc _{\mu _1 +\mu _2} (V), mc _{\mu _1 +\mu _2} (\boldsymbol{B} )),
\]
 if $\mu _1\! \neq \! 0 \! \neq \! \mu _2$.
Namely, the map $\overline{\phi } $ is an isomorphism of the vector spaces between $ mc _{\mu _2} (mc _{\mu _1} (V))$ and $mc _{\mu _1 +\mu _2} (V) $, and it satisfies
\begin{equation}
\overline{\phi } \bigl( \overline{G}_j (\mu _1 , \mu _2 ) \mathbf{v}\bigr) = \overline{G}_j (\mu _1 + \mu _2 ) \overline{\phi } ( \mathbf{v}) \label{eq:phiG}
\end{equation}
for $j=0,\dots ,N$ and $ \mathbf{v} \in mc _{\mu _2} (mc _{\mu _1} (V))$.
Set $\mu _1 = q^{\lambda _1} -1$ and $ \mu _2 = q^{\lambda _1 + \lambda _2} - q^{\lambda _1 }$.
Then, we have $\phi ^{q} = q^{\lambda _1} \phi $.
The map $\overline{\phi } ^{q} $ is an isomorphism of the vector spaces between $ mc _{\lambda _2} ^{q} \bigl(mc _{\lambda _1} ^{q} (V)\bigr) (= mc _{\mu _2} (mc _{\mu _1} (V)))$ and $mc _{\lambda _1 +\lambda _2} ^{q} (V) (= mc _{\mu _1 +\mu _2} (V)) $.
It follows from equations \eqref{eq:GjqlGjDRb}, \eqref{eq:Gjql1l2GDRb} and \eqref{eq:phiG} that
\begin{equation*}
\overline{\phi } ^{q} \bigl( \overline{G}_j^{q} (\lambda _1 , \lambda _2 ) \mathbf{v}\bigr) = \overline{G}_j^{q} (\lambda _1 + \lambda _2 ) \overline{\phi } ^{q} ( \mathbf{v}) ,
\end{equation*}
because the scalar $q^{\lambda _1 +\lambda _2}$ commutes with the linear map $\phi ^{q} $.
Hence, the map $\overline{\phi }^{q} $ is an isomorphism between $\bigl(mc _{\lambda _2} ^{q} \bigl(mc _{\lambda _1} ^{q} (V)\bigr), mc _{\lambda _2} ^{q} \bigl(mc _{\lambda _1} ^{q} (\boldsymbol{B} )\bigr)\bigr)$ and $\bigl(mc _{\lambda _1 +\lambda _2} ^{q} (V), mc _{\lambda _1 +\lambda _2} ^{q} (\boldsymbol{B} )\bigr)$.
\end{proof}

\section{Concluding remarks} \label{sec:rmk}

In this paper, we reformulated the $q$-middle convolution.
By applying the $q$-middle convolution and $q$-analogues of addition, we obtained several linear $q$-difference equations and their solutions.
We attached the spectral type, which was introduced by Sakai and Yamaguchi \cite{SY}, to these $q$-difference equations in Section~\ref{subsec:specttype}.
The spectral type would be useful for the study of undiscovered linear $q$-difference equations.

On the limit $q \to 1$, the $q$-difference equation
\begin{equation*}
\frac{Y(q x) - Y(x )}{-x} = \Biggl[ \frac{B_{0}}{x} + \sum^{N}_{i = 1}\frac{B_{i}}{x -b_{i}} \Biggr] Y(x)
\end{equation*}
might be related with the Fuchsian system of differential equation.

The problem on finding similarities or differences between the linear system of differential equation and the linear system of $q$-difference equation would not be studied very well.
In reference to this problem, we observe the limit $q \to 1$ of some $q$-difference equations appeared in this paper

Let us recall equation \eqref{eq:qhg-sta}, which is related to the $q$-hypergeometric equation.
Set $\beta/\alpha = q^{\nu}$, and assume that the parameter $\alpha $ does not depend on $q$.
Then, equation \eqref{eq:qhg-sta} is written as
\begin{align}
\frac{\widetilde{Y}(qx) - \widetilde{Y}(x)}{qx-x}
={}& \left[ \frac{1}{x} \begin{pmatrix}
 \displaystyle \frac{1 - q^{\mu - \lambda}}{1-q} & \displaystyle q^{\mu - \lambda}\frac{1 - q^{\nu}}{1-q} \\
 0 & 0
 \end{pmatrix} \nonumber\right. \\
& \left.+ \frac{1}{x-1/\alpha} \begin{pmatrix}
 0 & 0 \\
 \displaystyle q^{-\lambda} \frac{1 - q^{\mu}}{1-q} & \displaystyle q^{\mu - \lambda}\frac{1 - q^{\nu}}{1-q} + \frac{1 -q^{-\lambda}}{1-q}
 \end{pmatrix}
\right] \widetilde{Y}(x). \label{eq:qhg-sta-nu}
\end{align}
As $q\to 1$, equation \eqref{eq:qhg-sta-nu} tends to
\begin{equation} \label{eq:hg-sta-nu}
\frac{{\rm d}}{{\rm d}x} \widetilde{Y}(x)
= \bigg[ \frac{1}{x} \begin{pmatrix}
 \mu - \lambda & \nu \\
 0 & 0
 \end{pmatrix}
+ \frac{1}{x-1/\alpha} \begin{pmatrix}
 0 & 0 \\
 \mu & \nu - \lambda
 \end{pmatrix}
\bigg] \widetilde{Y}(x).
\end{equation}
Note that this equation is obtained by the convolution $c_{-\lambda }$ in equation \eqref{eq:Fi} from the single differential equation $y'= \{ \mu /x +\nu /(x-1/\alpha ) \} y $, which has a solution $y=x^{\mu } (x-1/\alpha )^{\nu }$.
Let~$\widetilde{y}_0(x) $ be the upper element of $\widetilde{Y}(x) $.
It follows from equation \eqref{eq:hg-sta-nu} that the function $\widetilde{y}_0(x) $ satisfies the hypergeometric differential equation by setting $\alpha =1$.
Conversely, equation \eqref{eq:qhg-sta} is a $q$-deformation of equation \eqref{eq:hg-sta-nu}.

We investigate the limit $q \to 1$ on the $q$-difference equation which is related to the generalized $q$-hypergeometric equation of order $3$.
Set $\beta/\alpha = q^{\nu}$ in equation \eqref{eq:111-21-111}, and assume that $\alpha $ does not depend on $q$.
By dividing by $1-q$ on equation \eqref{eq:111-21-111} and taking the limit as $q \to 1$, we have
\begin{align}
\frac{{\rm d}}{{\rm d}x} G (x)
={} &\left[ \frac{1}{x} \begin{pmatrix}
 \mu - \lambda + \mu' - \lambda' & \nu & 0 \\
 0 & \mu' - \lambda' & \mu \\
 0 & 0 & 0
 \end{pmatrix} \nonumber\right. \\
&\left. + \frac{1}{x-1/\alpha} \begin{pmatrix}
 0 & 0 & 0 \\
 0 & 0 & 0 \\
\mu ' + \mu - \lambda & \displaystyle \nu + \frac{(\nu - \lambda) \mu '}{\mu } & \nu - \lambda ' - \lambda
 \end{pmatrix}
\right] G (x). \label{eq:ghge}
\end{align}
The rank of the residue matrix about $x=1/\alpha $ is one, and we obtain the generalized hypergeometric equation of order $3$ from equation \eqref{eq:ghge} by setting $\alpha =1$ (see Oshima \cite[Example~5.3]{O}).

We discuss the limit $q \to 1$ on the equation \eqref{eq:21-111-111}.
Set $\beta/\alpha = q^{\nu}$ and $\alpha/\gamma = q^{\nu'}$ in equation~\eqref{eq:21-111-111}, and assume that $\alpha $ does not depend on $q$.
By dividing by $1-q$ on equation \eqref{eq:21-111-111} and taking the limit as $q \to 1$, we have
\begin{align}
\frac{{\rm d}}{{\rm d}x} G (x)
={}& \left[ \frac{1}{x} \begin{pmatrix}
 0 & 0 & 0 \\
 \nu - \lambda & \mu - \lambda - \lambda' & \displaystyle \frac{-\lambda\nu' + \mu\nu + \mu\nu'}{\nu} \\
 0 & 0 & 0
 \end{pmatrix} \nonumber\right. \\
&\left. + \frac{1}{x-1/\alpha } \begin{pmatrix}
 -\lambda - \lambda' + \nu + \nu' & 0 & \mu \\
 0 & 0 & 0 \\
 0 & \nu & -\lambda' + \nu'
 \end{pmatrix}
\right] G(x). \label{eq:ghgev}
\end{align}
The rank of the residue matrix about $x=0 $ is one, and we can obtain the generalized hypergeometric equation of order $3$ from equation \eqref{eq:ghgev} by setting $\alpha =1$ and replacing $x$ with~${1-x}$.

\appendix

\section[q-middle convolution by Sakai and Yamaguchi]{$\boldsymbol{q}$-middle convolution by Sakai and Yamaguchi} \label{sec:qmcSY}

Let $\boldsymbol{B}= ( B_{\infty};B_{1} ,\dots ,B_N ) $ be the tuple of the square matrices of the same size and $\mathbf{b}= (b_1, b_2, \dots ,b_N)$ be the tuple of the non-zero complex numbers which are different one another.
Sakai and Yamaguchi investigated the $q$-difference equation of the form
\begin{equation}
Y(q x) = B(x)Y(x), \qquad B(x) = B_{\infty} + \sum^{N}_{i = 1}\frac{B_{i}}{1 - x/b_{i}}.
\label{eq:YqxBYx1}
\end{equation}
\begin{Definition}[$q$-convolution \cite{SY}] \label{def:qc}
Let $ ( B_{\infty}; B_{1} ,\dots ,B_N ) $ be the tuple of $m\times m $ matrices and~${\lambda \in \mathbb{C}} $.
Set $B_0 = I_m - B_\infty - B_{1} - \dots -B_N$.
We define the $q$-convolution $c_\lambda^{\rm SY}\colon ( B_{\infty};B_{1} ,\dots ,\allowbreak B_N ) \mapsto ( F_{\infty};F_{1} ,\dots ,F_N )$ as follows:
\begin{gather}
 \boldsymbol{F} = ( F_\infty ; F_1, \dots , F_N) \mbox{ \rm is a tuple of $(N+1)m \times (N+1)m$ matrices,}\nonumber \\
 F_i = \left(
 \begin{matrix}
 {} & {} & O & {} & {} \\
 B_0 & \cdots & B_i - \bigl(1-q^\lambda\bigr)I_m & \cdots & B_N \\
 {} & {} & O & {} & {}
 \end{matrix}
 \right)
 {\scriptstyle (i+1)} , \qquad 1\leq i \leq N, \nonumber \\
 F_\infty = I_{(N+1)m} - \widehat{F}, \qquad
 \widehat{F}
 = \left(
 \begin{matrix}
 B_0 & \cdots & B_N \\
 \vdots & \cdots & \vdots \\
 B_0 & \cdots & B_N
 \end{matrix}
 \right).  \label{eq:bF}
\end{gather}
\end{Definition}
The $q$-convolution in Definition \ref{def:qc} induces the correspondence of the linear $q$-difference equations
\begin{align*}
& Y(q x) = B(x)Y(x) \mapsto \widehat{Y}(q x) = F(x)\widehat{Y}(x) , \\
&B(x) = B_{\infty} + \sum^{N}_{i = 1}\frac{B_{i}}{1 - x/b_{i}} ,\qquad F(x) = F_{\infty} + \sum^{N}_{i = 1}\frac{F_{i}}{1 - x/b_{i}} ,
\end{align*}
and it is related with a $q$-integral transformation, which was established by Sakai and Yamaguchi~\cite{SY}.
\begin{Theorem}[{cf.~\cite[Theorem 2.1]{SY}}] \label{thm:qcintSY}
Let $Y(x)$ be a solution to equation \eqref{eq:YqxBYx1}.
Set $b_0 =0$ and
\begin{equation}
 P_{\lambda}(x, s) = \frac{\bigl(q^{\lambda +1} s/x;q\bigr)_{\infty }}{(q s/x;q)_{\infty }}.  \label{eq:Plaxss/x}
\end{equation}
Define the function $\widehat{Y}(x) $ by
\begin{equation}
 \widehat{Y}_{i}(x) = \int^{\xi \infty}_{0}\frac{P_{\lambda}(x, s)}{s-b_{i}}Y(s) {\rm d}_{q}s,\qquad i=0,\dots ,N, \qquad \widehat{Y}(x) = \begin{pmatrix} \widehat{Y}_{0}(x) \\ \vdots \\ \widehat{Y}_{N}(x) \end{pmatrix}.  \label{eq:qcint0}
\end{equation}
Then, the function $\widehat{Y}(x)$ formally satisfies
\[
\widehat{Y}(q x) = \left( F_{\infty} + \sum^{N}_{i = 1}\frac{F_{i}}{1 - x/b_{i}} \right) \widehat{Y}(x).
\]
\end{Theorem}
Note that the original theorem by Sakai and Yamaguchi was restricted to the case $\xi =1$ in the Jackson integral, and convergence was discussed in \cite{AT}.

\begin{Definition}[$q$-middle convolution, \cite{SY}]\label{def:qmcSY}
We define the subspaces $\mathcal{K}$ and $\mathcal{L} $ of $(\mathbb{C}^m )^{N+1}$ as follows;
\begin{equation*}
 \mathcal{K} =
 \left(
 \begin{matrix}
 \operatorname{ker}B_0 \\
 \vdots \\
 \operatorname{ker}B_N
 \end{matrix}
 \right), \qquad
 \mathcal{L} =
 \operatorname{ker}\bigl(\widehat{F} - \bigl(1 - q^{\lambda}\bigr)I_{(N+1)m}\bigr).
\end{equation*}
Then, the spaces $\mathcal{K}$ and $\mathcal{L} $ are invariant under the action of $F_k$, and we denote the matrix induced from the action of $F_k$ on the quotient space $(\mathbb{C}^m )^{N+1}/(\mathcal{K} + \mathcal{L})$ by $\overline{F}_k$ $(k=0 ,1, \dots ,N, \infty )$.
The $q$-middle convolution $mc_\lambda ^{\rm SY}$ is defined by the correspondence
\[
( B_{\infty}; B_{1} ,\dots ,B_N ) \mapsto \bigl(\overline{F}_{\infty}; \overline{F}_1 \dots ,\allowbreak\overline{F}_N \bigr)
\]
 (or the correspondence of associated $q$-difference equations).
\end{Definition}

In \cite{SY}, the transformation $\psi _{\mu } $ $(\mu \in \mathbb{C} )$ on the tuple of square matrices of size $m'$ was defined~by
$
 \psi _{\mu }\colon ( F_{\infty } ; F_{1} ,\dots ,F_N ) \mapsto ( F_{\infty } +(1- q^{\mu }) I_{m'} ; F_{1} ,\dots ,F_N )$,
and the map $\overline{\Psi }_{\lambda } $ was defined by
\smash{$
 \overline{\Psi }_{\lambda } = \psi _{\lambda } \circ mc_{\lambda } ^{\rm SY} $}.
The relation which was mentioned in equation~\eqref{eq:Psimulam0} is described precisely as the following proposition.
\begin{Proposition}[{\cite[Proposition 4.12]{SY}}]
Under the conditions $(*)$, $(**)$ in \eqref{eq:def*} and \eqref{eq:def**}, we have
\smash{$
 \overline{\Psi }_{\mu } \circ \overline{\Psi }_{\lambda } \simeq \overline{\Psi }_{\log ( q^{\lambda } + q^{\mu } -1)/ \log q } $}.
\end{Proposition}
We compare the $q$-convolution and the $q$-middle convolution of Sakai and Yamaguchi with the ones in this manuscript.
The $q$-difference equation \eqref{eq:YqxBYx1} is equivalent to
\begin{equation*}
\frac{Y(q x) - Y(x )}{-x} = \Biggl[ \sum^{N}_{i = 0}\frac{B_{i}}{x -b_{i}} \Biggr] Y(x)
\end{equation*}
by setting $B_0 = I_m - B_\infty - B_{1} - \dots -B_N$ and $b_0=0$.
On the matrices $G_0 ,G_1 , \dots , G_N $ in equation \eqref{eq:bG} and the ones $F_1, \dots , F_N , F_\infty , \widehat{F}$ in equation \eqref{eq:bF}, we have
\begin{equation*}
 G_1 = q^{-\lambda } F_1 , \dots , G_N = q^{-\lambda } F_N , \qquad G_0 +G_1 + \dots + G_N = \bigl(1-q^{-\lambda }\bigr) I_{(N+1)m} + q^{-\lambda } \widehat{F}.
\end{equation*}
Set $G_{\infty} = I_{(N+1)m} - G_0 - G_1 - \dots - G_N $.
Then, we have $
G_{\infty} = q^{-\lambda } I_{(N+1)m}- q^{-\lambda } \widehat{F} = q^{-\lambda } F_{\infty} $.
On the function \smash{$K^{(1)}_{\lambda}(x, s)$} in equation \eqref{eq:tP1tP2} and the function $ P_{\lambda}(x, s) $ in equation~\eqref{eq:Plaxss/x}, there is a relation
\smash{$K^{(1)}_{\lambda}(x, s) = x^{-\lambda } P_{\lambda}(x, s) $}.
If we take \smash{$K_{\lambda}(x, s) = K^{(1)}_{\lambda}(x, s) $}, then we have \smash{$\widetilde{Y}(x) = x^{-\lambda } \widehat{Y} (x)$} on the functions $\widetilde{Y}(x) $ and $\widehat{Y} (x)$ in equations \eqref{eq:qcintadd} and \eqref{eq:qcint0}.
Hence, we may regard the $q$-convolution $c^{q} _{\lambda} $ as the composition of the $q$-convolution $c^{\rm SY} _{\lambda} $ of Sakai and Yamaguchi and the addition induced by the gauge transformation with respect to $x^{-\lambda }$.
Namely, we have~${c^{q} _{\lambda} = {\rm add}_{-\lambda } \circ c^{\rm SY} _{\lambda} }$.
If $\lambda $ is an integer, then we have
\begin{equation*}
K^{(1)}_{\lambda}(x, s) = x^{-\lambda } \frac{\bigl(q^{\lambda+1}s/x;q\bigr)_{\infty}}{(qs/x;q)_{\infty}} \to x^{-\lambda } (1-s/x) ^{-\lambda } = (x-s)^{-\lambda }
\end{equation*}
as $q \to 1$.
Thus, the $q$-integral transformation with respect to the kernel function $K^{(1)}_{\lambda}(x, s) $ is a direct $q$-deformation of the Euler's integral transformation, and the $q$-analogue by the $q$\nobreakdash-convolution $c^{q} _{\lambda} $ would be more suitable from the aspect of the integral transformation.

\section[Supplement to Section 4.6]{Supplement to Section~\ref{subsec:3exte}} \label{appsec:B}

We derive single third-order $q$-difference equations from the systems of three $q$-difference equations obtained in Section~\ref{subsec:3exte} and also discuss their $q$-integral representations of solutions.

\subsection[Supplement to Section 4.6.1]{Supplement to Section~\ref{subsec:3extedeg2}}

The single third-order $q$-difference equation for the function $g_1(x)$ from equation \eqref{eq:3-111-21-111} is derived~as
\begin{gather}
(\alpha_2 x-1)(q \gamma_1 x-1)\bigl(q^2 \gamma_1 x-1\bigr) g_1\bigl(q^3x\bigr)- q^{-\lambda-\lambda'} [ c_{2,2} x^2
- c_{2,1} x + q^{\lambda + \lambda'}\bigl(1+q+q^2\bigr) ] \nonumber \\
\qquad\times
 (q \gamma_1 x-1) g_1\bigl(q^2x\bigr)+ q^{-2\lambda-2\lambda'+1} \big[ c_{1,3} x^3 - c_{1,2} x^2 + c_{1,1} x - q^{2\lambda+2\lambda'}\bigl(1+q+q^2\bigr) \big] g_1(qx) \nonumber \\
\phantom{\qquad\times}{}- q^3\bigl(q^{-\lambda'}\alpha_1 x - 1\bigr)\bigl(q^{-\lambda-\lambda'}\beta_1 x - 1\bigr)\bigl(q^{-\lambda-\lambda'}\beta_2 x - 1\bigr) g_1(x)=0, \label{eq:3-111-21-111-g1}
\end{gather}
where
\begin{gather*}
 c_{2,2} = q^2\bigl(q^{\lambda}\alpha_2\gamma_1 + \alpha_1\alpha_2 + \beta_1\beta_2\bigr), \\
 c_{2,1} = q^{\lambda + \lambda' + 2}\gamma_1 + q^{\lambda + 2}\alpha_1 + q^2(\beta_1 + \beta_2) + q^{\lambda + \lambda'}(1 + q)\alpha_2 , \\
 c_{1,3} = q^2\bigl(q^{\lambda}\beta_1\beta_2\gamma_1 + q^{\lambda}\alpha_1\alpha_2\gamma_1 + \alpha_1\beta_1\beta_2\bigr) , \\
 c_{1,2} = q \bigl\{q^{\lambda+1}\bigl(q^{\lambda'}\gamma_1 + \alpha_1\bigr)(\beta_1 + \beta_2) + q\beta_1\beta_2 + q^{2\lambda+\lambda'}\gamma_1(q\alpha_1 + \alpha_2) + q^{\lambda+\lambda'}(\alpha_1\alpha_2 + \beta_1\beta_2) \bigr\} , \\
 c_{1,1} = q^{\lambda+\lambda'} \bigl\{q^{\lambda+\lambda'+1}(1 + q)\gamma_1 + q^{\lambda+\lambda'}\alpha_2 + q^{\lambda+1}(1 + q)\alpha_1 + q(1 + q)(\beta_1 + \beta_2) \bigr\}.
\end{gather*}
Each coefficient of $g_1 \bigl(q^j x\bigr)$ $(j=0,1,2,3)$ is a cubic polynomial in $x$.

We discuss $q$-integral representations of solutions.
We obtain $q$-integral representations of solutions $Y_g(x) $ to equation \eqref{eq:B'1B'2gm1} by using the $q$-integral in equation \eqref{eq:qintvarqhg2} and the transformation in equation \eqref{eq:Ygdeg2-1}.
By the $q$-integral transformation induced from the $q$-convolution $c_{\lambda'}^q$ and the transformation by the matrix $P$, we obtain formal $q$-integral representations of solutions to equation \eqref{eq:3-111-21-111}.
The formal $q$-integral representation for $g_1(x)$ is
expressed as
\begin{equation*}
 g_1(x) = \int_0^{\xi' \infty} \frac{\widehat{K}_{\lambda'}(x,s)}{s - 1/\gamma_1} \frac{(\gamma_1 s; q)_{\infty}}{(\alpha_1 s; q)_{\infty}} \int_0^{\xi \infty} \frac{K_{\lambda}(s,t)}{t - 1/\alpha_1} y(t) {\rm d}_q t {\rm d}_q s,
\end{equation*}
where $\widehat{K}_{\lambda' }(x,s) $ and $K_{\lambda }(s,t)$ are functions which satisfy equation \eqref{eq:tPlambda}, and $y(t)$ is a function which satisfies equation \eqref{eq:qeqvar2}.
By specializing to \smash{$\widehat{K}_{\lambda' }(x,s) = K_{\lambda' }^{(1)}(x,s)$}, \smash{$K_{\lambda }(s,t) = K_{\lambda }^{(1)}(s,t)$} and $y(t) = (\alpha_1 t , \alpha_2 t; q)_{\infty}/(\beta_1 t , \beta_2 t; q)_{\infty}$, the formal solution to equation \eqref{eq:3-111-21-111-g1} is written as
\begin{equation*}
 g_1(x) = \alpha_1 \gamma_1 x^{-\lambda '} \int_0^{\xi' \infty} \frac{\bigl(q^{\lambda '+1}s/x, q \gamma_1 s ;q\bigr)_{\infty}}{(qs/x, \alpha_1 s ;q)_{\infty}} s^{-\lambda } \int_0^{\xi \infty} \frac{\bigl(q^{\lambda+1}t/s, q \alpha_1 t, \alpha_2 t ;q\bigr)_{\infty}}{(qt/s, \beta_1 t, \beta_2 t ;q)_{\infty}} {\rm d}_q t {\rm d}_q s.
\end{equation*}
We need to examine the condition that a formal solution converges and it is an actual solution to equation \eqref{eq:3-111-21-111-g1}, and we do not discuss it here.

The single third-order $q$-difference equation for the function $g_3(x)$ from equation \eqref{eq:3-111-111-21} is derived as
\begin{gather}
(\gamma_1 x-1)(q \gamma_1 x-1)(q \gamma_2 x-1)\bigl(q^2 \gamma_2 x-1\bigr)g_3\bigl(q^3x\bigr) \nonumber \\
\qquad{}- q^{-\lambda'}\Big\{ q^2 \frac{\gamma_1\gamma_2}{\alpha_1\alpha_2}\bigl(\alpha_1\alpha_2 + q\alpha_1\alpha_2 + q^{\lambda'+1}\beta_1\beta_2\bigr)x^2 +r_1 x + q^{\lambda'}\bigl(1+q+q^2\bigr) \Big\}(\gamma_1 x-1) \nonumber \\
\qquad\phantom{-}{}\times (q \gamma_2 x-1)g_3\bigl(q^2x\bigr) + q^{-2\lambda'+1} \bigl\{ s_4 x^4 + s_3 x^3 + s_2 x^2 + s_1 x + q^{2\lambda'}\bigl(1+q+q^2\bigr) \bigr\} g_3(qx)\! \nonumber \\
\qquad{}- q^3 \bigl(q^{-\lambda'}\alpha_1 x - 1\bigr)\bigl(q^{-\lambda'}\alpha_2 x - 1\bigr)\bigg(\frac{\beta_1\gamma_1\gamma_2}{\alpha_1\alpha_2} x - 1 \bigg)\bigg(\frac{\beta_2\gamma_1\gamma_2}{\alpha_1\alpha_2} x - 1 \bigg) g_3(x)=0, \label{eq:3-111-111-21-g3}
\end{gather}
where
\begin{gather*}
 r_1 = - \frac{q}{\alpha_1\alpha_2}
\bigl\{ q^{\lambda' + 1}(\beta_1+\beta_2)\gamma_1\gamma_2 + q^{\lambda'}\alpha_1\alpha_2(\gamma_1 + q\gamma_2) + q\alpha_1\alpha_2(\alpha_1+ \alpha_2) \bigr\} , \\
s_4 = q^2\frac{\gamma_1^2\gamma_2^2}{\alpha_1\alpha_2} \bigl\{q^{\lambda'}\beta_1\beta_2(1+q) + \alpha_1\alpha_2 \bigr\}, \\
s_3 = -q\frac{\gamma_1\gamma_2}{\alpha_1^2\alpha_2^2}
\bigl\{ q^{\lambda'+1}\alpha_1\alpha_2(\alpha_1\alpha_2 + \gamma_1\gamma_2)(\beta_1 + \beta_2) + q\bigl(\alpha_1^2\alpha_2^2 +
q^{2\lambda'}\beta_1\beta_2\gamma_1\gamma_2\bigr)(\alpha_1 + \alpha_2) \\
\phantom{s_3 =}{} + q^{\lambda'}\alpha_1\alpha_2(\alpha_1\alpha_2 + q\beta_1\beta_2)(\gamma_1 + q\gamma_2) \bigr\}, \\
s_2 = \frac{q}{\alpha_1^2\alpha_2^2}
\big[ q^{\lambda'}\alpha_1^2\alpha_2^2(\gamma_1 + q\gamma_2)(\alpha_1 + \alpha_2) + q^{\lambda'}\alpha_1^2\alpha_2^2\gamma_1\gamma_2\bigl(q^{\lambda'} + 1 + q\bigr) \\
\phantom{s_2 =}{} + q^{\lambda'}\alpha_1\alpha_2\gamma_1\gamma_2 \bigl\{ q(\alpha_1 + \alpha_2) + q^{\lambda'}(\gamma_1 + q\gamma_2) \bigr\} (\beta_1 + \beta_2) \\
\phantom{s_2 =}{} + q^{2\lambda'+1}\beta_1\beta_2\gamma_1\gamma_2(\alpha_1\alpha_2 + \gamma_1\gamma_2) + q\alpha_1^3\alpha_2^3\big], \\
s_1 = -\frac{q^{\lambda'}(1+q)}{\alpha_1\alpha_2}
\bigl\{ q^{\lambda'+1}\gamma_1\gamma_2(\beta_1 + \beta_2) + q^{\lambda'}\alpha_1\alpha_2(\gamma_1 + q\gamma_2) + q\alpha_1\alpha_2(\alpha_1 + \alpha_2) \bigr\}.
\end{gather*}
Each coefficient of $g_3 \bigl(q^j x\bigr)$ $(j=0,1,2,3)$ is a quartic polynomial in $x$.
Remark that the coefficients in equation \eqref{eq:3-111-111-21-g3} are expected to be characterized by some underlying structural conditions of equation \eqref{eq:3-111-111-21}.

By the $q$-integral transformation induced from the $q$-convolution $c_{\lambda'}^q$ and the transformation by the matrix $P$, we obtain formal $q$-integral representations of solutions to equation \eqref{eq:3-111-111-21}, and they include the expression
\begin{align*}
g_3(x) ={}& \frac{-\alpha_1 (\alpha_1 - \alpha_2) \gamma_2}{(\alpha_1 - \beta_1)(\alpha_1 - \beta_2)} x^{-\lambda'} \int_0^{\xi' \infty} \frac{\bigl(q^{\lambda' + 1}s/x, \gamma_1 s, q\gamma_2 s; q\bigr)_{\infty}}{(qs/x, \alpha_1 s, \alpha_2 s; q)_{\infty}} s^{-\lambda} \\
&\times \int_0^{\xi \infty} \frac{\bigl(q^{\lambda + 1}t/s, q\alpha_1 t, q\alpha_2 t; q\bigr)_{\infty}}{(qt/s, \beta_1 t, \beta_2 t; q)_{\infty}} (-\alpha_1-\alpha_2 + \beta_1 + \beta_2 + (\alpha_1\alpha_2 - \beta_1\beta_2)t) {\rm d}_q t {\rm d}_q s
,
\end{align*}
where $q^{\lambda + \lambda'} = \alpha_1\alpha_2/(\gamma_1\gamma_2)$.
Since the assumption of Theorem~\ref{thm:qcint} is not satisfied in general for this case, the formal $q$-integral representations may not satisfy the designated $q$-difference equation, and we need to discuss more to obtain actual solutions.

\subsection[Supplement to Section 4.6.2]{Supplement to Section~\ref{subsec:3extedeg3}}

The single third-order $q$-difference equation for the function $g_3(x)$ from equation \eqref{eq:3-111-111-21-3} is derived~as
\begin{gather}
\bigl(q^2 \alpha_3 x-1\bigr)(\gamma_1 x-1)(q \gamma_1 x-1)(q \gamma_2 x-1)\bigl(q^2 \gamma_2 x-1\bigr)g_3\bigl(q^3x\bigr) \nonumber \\
\qquad{}+ \bigl\{ t_3 x^3 + t_2 x^2 + t_1 x +\bigl(1+q+q^2\bigr) \bigr\} (\gamma_1 x-1)(q \gamma_2 x-1) g_3\bigl(q^2x\bigr) \nonumber \\
\qquad{}+ \bigl\{ u_5 x^5 + u_4 x^4 + u_3 x^3 + u_2 x^2 + u_1 x - q\bigl(1+q+q^2\bigr) \bigr\} g_3(qx) \nonumber \\
\qquad{}- q^3 \Biggl(\frac{\beta_1\gamma_1\gamma_2}{\alpha_1\alpha_2} x - 1 \Biggr) \Biggl(\frac{\beta_2\gamma_1\gamma_2}{\alpha_1\alpha_2} x - 1 \Biggr) \Biggl(\frac{\beta_3\gamma_1\gamma_2}{\alpha_1\alpha_2} x - 1 \Biggr) \nonumber \\
\phantom{\qquad-}{}\times \Biggl( \frac{\beta_1\beta_2\beta_3\gamma_1\gamma_2}{\alpha_1^2\alpha_2\alpha_3} x - 1 \Biggr) \Biggl( \frac{\beta_1\beta_2\beta_3\gamma_1\gamma_2}{\alpha_1\alpha_2^2\alpha_3} x - 1 \Biggr) g_3(x)=0, \label{eq:3-111-111-21-3-g3}
\end{gather}
where
\begin{gather*}
t_3 = -q^3\bigl(1+q+q^2\bigr) \frac{\beta_1\beta_2\beta_3\gamma_1^2\gamma_2^2}{\alpha_1^2\alpha_2^2}, \\
t_2 = q^2\frac{\gamma_1\gamma_2}{\alpha_1^2\alpha_2^2\alpha_3}
\bigl\{ (1 + q)\beta_1\beta_2\beta_3\gamma_1\gamma_2 + q\alpha_3\beta_1\beta_2\beta_3(\gamma_1 + q\gamma_2) \\
\phantom{t_2 =}{} + q\alpha_1\alpha_2\alpha_3\beta_1(\beta_2 + \beta_3) + q\alpha_1\alpha_2\alpha_3^2(\alpha_1 + \alpha_2) + q\alpha_1\alpha_2\alpha_3\beta_2\beta_3 \bigr\}, \\
t_1 = -\frac{q}{\alpha_1^2\alpha_2^2\alpha_3}
\bigl\{ (1 + q)q\alpha_1^2\alpha_2^2\alpha_3^2 + q\alpha_1\alpha_2\alpha_3\gamma_1\gamma_2(\beta_1 + \beta_2 + \beta_3) \\
\phantom{t_1 =}{} + q\beta_1\beta_2\beta_3\gamma_1\gamma_2(\alpha_1 + \alpha_2) + \alpha_1^2\alpha_2^2\alpha_3(\gamma_1 + q\gamma_2) \bigr\}, \\
u_5 = q^3 \bigl(1+q+q^2\bigr) \frac{\beta_1^2\beta_2^2\beta_3^2\gamma_1^4\gamma_2^4}{\alpha_1^4\alpha_2^4\alpha_3}, \\
u_4 = -\frac{q^3\beta_1\beta_2\beta_3\gamma_1^3\gamma_2^3}{\alpha_1^4\alpha_2^4\alpha_3^2}
\bigl\{ \beta_1\beta_2\beta_3\gamma_1\gamma_2 + (1 + q)\alpha_3\beta_1\beta_2\beta_3(\gamma_1 + q\gamma_2) \\
\phantom{u_4 =}{} + (1 + q)\alpha_1\alpha_2\alpha_3(\beta_1\beta_2 + \beta_2\beta_3 + \beta_3\beta_1) + (1 + q)\alpha_1\alpha_2\alpha_3^2(\alpha_1 + \alpha_2) \bigr\}, \\
u_3 = \frac{q^2\gamma_1^2\gamma_2^2}{\alpha_1^4\alpha_2^4\alpha_3^2}
\big[ \alpha_1\alpha_2\alpha_3\beta_1\beta_2\beta_3(\gamma_1 + q\gamma_2)(\alpha_1\alpha_2 + q\alpha_2\alpha_3 + q\alpha_3\alpha_1 + q\beta_1\beta_2 + q\beta_2\beta_3 +q\beta_3\beta_1) \\
\phantom{u_3 =}{} + q\alpha_1\alpha_2\alpha_3\beta_1\beta_2\beta_3\gamma_1\gamma_2(\beta_1 + \beta_2 + \beta_3)
+ q\beta_1^2\beta_2^2\beta_3^2\gamma_1\gamma_2(\alpha_1 + \alpha_2 + q\alpha_3) \\
\phantom{u_3 =}{} + q\alpha_1\alpha_2\alpha_3\beta_1\beta_2\beta_3 \{\alpha_1\alpha_2(\beta_1 + \beta_2 + \beta_3) + (1 + q)\alpha_1\alpha_2\alpha_3 \} \\
\phantom{u_3 =}{} + q\alpha_1^2\alpha_2^2\alpha_3^2(\alpha_1 + \alpha_2)(\beta_1\beta_2 + \beta_2\beta_3 + \beta_3\beta_1) + q\alpha_1^3\alpha_2^3\alpha_3^3 \big], \\
u_2 = -\frac{q^2\gamma_1\gamma_2}{\alpha_1^3\alpha_2^3\alpha_3^2}
\bigl\{ q\alpha_3\beta_1\beta_2\beta_3\gamma_1\gamma_2(\alpha_1 + \alpha_2)(\beta_1 + \beta_2 + \beta_3) \\
\phantom{u_2 =}{} + q\alpha_1\alpha_2\alpha_3^2\gamma_1\gamma_2(\beta_1\beta_2 + \beta_2\beta_3 + \beta_3\beta_1) \\
\phantom{u_2 =}{} + (1 + q)\alpha_1\alpha_2\alpha_3\beta_1\beta_2\beta_3\gamma_1\gamma_2 + \alpha_1\alpha_2\alpha_3\beta_1\beta_2\beta_3(\gamma_1 + q\gamma_2)(\alpha_1 + \alpha_2 + q\alpha_3) \\
\phantom{u_2 =}{} + \alpha_1^2\alpha_2^2\alpha_3^2(\gamma_1 + q\gamma_2)(\beta_1 + \beta_2 + \beta_3) + \alpha_1^2\alpha_2^2\alpha_3^2(\alpha_1\alpha_2 + q\alpha_2\alpha_3 + q\alpha_3\alpha_1) \\
\phantom{u_2 =}{} + q\alpha_1^2\alpha_2^2\alpha_3^2(\beta_1\beta_2 + \beta_2\beta_3 + \beta_3\beta_1) + q\beta_1^2\beta_2^2\beta_3^2\gamma_1\gamma_2 \bigr\}, \\
u_1 = \frac{q}{\alpha_1^2\alpha_2^2\alpha_3}
\bigl\{ q^2\alpha_1^2\alpha_2^2\alpha_3^2 + (1 + q)q\alpha_1\alpha_2\alpha_3\gamma_1\gamma_2(\beta_1 + \beta_2 + \beta_3) \\
\phantom{u_1 =}{} + (1 + q)\alpha_1^2\alpha_2^2\alpha_3(\gamma_1 + q\gamma_2) + (1 + q)q\beta_1\beta_2\beta_3\gamma_1\gamma_2(\alpha_1 + \alpha_2) \bigr\}.
\end{gather*}
Each coefficient of $g_3 \bigl(q^j x\bigr)$ $(j=0,1,2,3)$ is a polynomial in $x$ of degree $5$.
Remark that the coefficients in equation \eqref{eq:3-111-111-21-3-g3} are expected to be characterized by some conditions
inherent in equation \eqref{eq:3-111-111-21-3}.

By the $q$-integral transformation induced from the $q$-convolution $c_{\lambda'}^q$ and the transformation by the matrix $P$, we obtain formal $q$-integral representations of solutions to equation \eqref{eq:3-111-111-21-3}, and they include the expression
\begin{align*}
g_3(x) ={}& \frac{(\alpha_3 - \alpha_1)(\alpha_1 - \alpha_2)(\gamma_2 - \alpha_3)}{(\alpha_1 - \beta_1)(\alpha_1 - \beta_2)(\alpha_1 - \beta_3)\alpha_2} x^{-\lambda'} \int_0^{\xi' \infty} \frac{\bigl(q^{\lambda'+1}s/x, \gamma_1 s, q\gamma_2 s; q\bigr)_{\infty}}{(qs/x, \alpha_1 s, \alpha_2 s; q)_{\infty}} \\
& \times \frac{s^{-\lambda}}{1 - \alpha_3 s} \int_0^{\xi \infty} \frac{\bigl(q^{\lambda+1}t/s, q\alpha_1 t, q\alpha_2 t, q\alpha_3 t; q\bigr)_{\infty}}{(qt/s, \beta_1 t, \beta_2 t, \beta_3 t; q)_{\infty}} (\eta_1 - \eta_2 t) {\rm d}_q t {\rm d}_q s,
\end{align*}
where
\begin{gather*}
q^{\lambda} = \beta_1\beta_2\beta_3/(\alpha_1\alpha_2\alpha_3),\qquad q^{\lambda'} = \alpha_1^2\alpha_2^2\alpha_3/(\beta_1\beta_2\beta_3\gamma_1\gamma_2), \\ \eta_1 = \alpha_1\alpha_2(\alpha_1 + \alpha_2) - \alpha_1\alpha_2(\beta_1 + \beta_2 + \beta_3) + \beta_1\beta_2\beta_3,\qquad \text{and}\\
\eta_2 = \alpha_1^2\alpha_2^2 - \alpha_1\alpha_2(\beta_1\beta_2 + \beta_2\beta_3 + \beta_3\beta_1) + (\alpha_1 + \alpha_2)\beta_1\beta_2\beta_3.
\end{gather*}
Since the assumption of Theorem~\ref{thm:qcint} is not satisfied in general for this case, the formal $q$\nobreakdash-integral representations may not satisfy the designated $q$-difference equation, and we need to discuss more to obtain actual solutions.

\subsection*{Acknowledgements}
The authors are grateful to the anonymous referees for their insightful comments, which helped to improve the manuscript.
The second author was supported by JSPS KAKENHI Grant Number JP22K03368.

\pdfbookmark[1]{References}{ref}
\LastPageEnding


\begin{thebibliography}{99}
\footnotesize\itemsep=0pt

\bibitem{Ar}
Arai Y.,
Solutions to $q$-hypergeometric equations associated with $q$-middle convolution,
 \href{http://arxiv.org/abs/2403.02662v2}{arXiv:2403.02662v2}.

\bibitem{AT}
Arai Y., Takemura K., On {$q$}-middle convolution and {$q$}-hypergeometric
 equations, \href{https://doi.org/10.3842/SIGMA.2023.037}{\textit{SIGMA}}
 \textbf{19} (2023), 037, 40~pages,
 \href{http://arxiv.org/abs/2209.02227}{arXiv:2209.02227}.

\bibitem{DR1}
Dettweiler M., Reiter S., An algorithm of {K}atz and its application to the
 inverse {G}alois problem, \href{https://doi.org/10.1006/jsco.2000.0382}{\textit{J.~Symbolic Comput.}} \textbf{30} (2000), 761--798.

\bibitem{DR2}
Dettweiler M., Reiter S., Middle convolution of {F}uchsian systems and the
 construction of rigid differential systems,
 \href{https://doi.org/10.1016/j.jalgebra.2007.08.029}{\textit{J.~Algebra}}
 \textbf{318} (2007), 1--24.

\bibitem{FN}
Fujii T., Nobukawa T., Hypergeometric solutions for variants of the
 $q$-hypergeometric equation,
 \href{http://arxiv.org/abs/2207.12777}{arXiv:2207.12777}.

\bibitem{GR}
Gasper G., Rahman M., Basic hypergeometric series, 2nd ed., \textit{Encyclopedia Math.
 Appl.}, Vol.~96,
 \href{https://doi.org/10.1017/CBO9780511526251}{Cambridge University Press},
 Cambridge, 2004.

\bibitem{Har}
Haraoka Y., Linear differential equations in the complex domain. {F}rom
 classical theory to forefront, \textit{Lecture Notes in Math.}, Vol.~2271,
 \href{https://doi.org/10.1007/978-3-030-54663-2}{Springer}, Cham, 2020.

\bibitem{HMST}
Hatano N., Matsunawa R., Sato T., Takemura K., Variants of {$q$}-hypergeometric
 equation, \href{https://doi.org/10.1619/fesi.65.159}{\textit{Funkcial.
 Ekvac.}} \textbf{65} (2022), 159--190,
 \href{http://arxiv.org/abs/1910.12560}{arXiv:1910.12560}.

\bibitem{JS}
Jimbo M., Sakai H., A {$q$}-analog of the sixth {P}ainlev\'e equation,
 \href{https://doi.org/10.1007/BF00398316}{\textit{Lett. Math. Phys.}}
 \textbf{38} (1996), 145--154,
 \href{http://arxiv.org/abs/chao-dyn/9507010}{arXiv:chao-dyn/9507010}.

\bibitem{KK}
Kakei S., Kikuchi T., A {$q$}-analogue of~{$\widehat{\mathfrak{gl}}_3$}
 hierarchy and {$q$}-{P}ainlev\'e~{VI},
 \href{https://doi.org/10.1088/0305-4470/39/39/S11}{\textit{J.~Phys.~A}}
 \textbf{39} (2006), 12179--12190,
 \href{http://arxiv.org/abs/nlin.SI/0605052}{arXiv:nlin.SI/0605052}.

\bibitem{Katz}
Katz N.M., Rigid local systems, \textit{Ann. of Math. Stud.}, Vol.~139,
 \href{https://doi.org/10.1515/9781400882595}{Princeton University Press},
 Princeton, NJ, 1996.

\bibitem{O}
Oshima T., Classification of {F}uchsian systems and their connection problem,
 in Exact {WKB} {A}nalysis and {M}icrolocal {A}nalysis, \textit{RIMS K\^{o}ky\^{u}roku
 Bessatsu}, Vol.~B37, Research Institute for Mathematical Sciences (RIMS), Kyoto,
 2013, 163--192, \href{http://arxiv.org/abs/0811.2916}{arXiv:0811.2916}.

\bibitem{SY}
Sakai H., Yamaguchi M., Spectral types of linear {$q$}-difference equations and
 {$q$}-analog of middle convolution,
 \href{https://doi.org/10.1093/imrn/rnw089}{\textit{Int. Math. Res. Not.}}
 \textbf{2017} (2017), 1975--2013,
 \href{http://arxiv.org/abs/1410.3674}{arXiv:1410.3674}.

\bibitem{STT}
Sasaki S., Takagi S., Takemura K., {$q$}-middle convolution and
 {$q$}-{P}ainlev\'e equation,
 \href{https://doi.org/10.3842/SIGMA.2022.056}{\textit{SIGMA}} \textbf{18}
 (2022), 056, 21~pages,
 \href{http://arxiv.org/abs/2201.03960}{arXiv:2201.03960}.

\bibitem{TakR}
Takemura K., Degenerations of {R}uijsenaars--van {D}iejen operator and
 {$q$}-{P}ainlev\'e equations,
 \href{https://doi.org/10.1093/integr/xyx008}{\textit{J.~Integrable Syst.}}
 \textbf{2} (2017), xyx008, 27~pages,
 \href{http://arxiv.org/abs/1608.07265}{arXiv:1608.07265}.

\bibitem{TakqH}
Takemura K., On {$q$}-deformations of the {H}eun equation,
 \href{https://doi.org/10.3842/SIGMA.2018.061}{\textit{SIGMA}} \textbf{14}
 (2018), 061, 16~pages,
 \href{http://arxiv.org/abs/1712.09564}{arXiv:1712.09564}.

\bibitem{Tki}
Takemura K., Kernel function, {$q$}-integral transformation and {$q$}-{H}eun
 equations, \href{https://doi.org/10.3842/SIGMA.2024.083}{\textit{SIGMA}}
 \textbf{20} (2024), 083, 22~pages,
 \href{http://arxiv.org/abs/2309.09341}{arXiv:2309.09341}.

\bibitem{Yg}
Yamaguchi M., The rigidity index of the linear $q$-difference equation and the
 $q$-middle convolution, {M}aster Thesis, {U}niversity of Tokyo, 2011.

\end{thebibliography}
\end{document}